\documentclass[12pt, final, 3p, times]{Qarticle}
\usepackage{epsfig}
\usepackage{amssymb}
\usepackage{float}
\usepackage{amsfonts}
\usepackage{mathrsfs}
\usepackage{amsthm}
\usepackage{amsmath}
\usepackage{amscd}
\usepackage{graphicx,pgfplots}
\usepackage{lineno}
\usepackage{mathrsfs}

\usepackage[colorlinks,citecolor=blue,linkcolor=blue,urlcolor=blue]{hyperref}
\usepackage{graphicx}
\usepackage{subfigure}

\usepackage{stmaryrd}
\usepackage{geometry}
\usepackage{fancyhdr}
\theoremstyle{plain}

\theoremstyle{definition}

\numberwithin{equation}{section}
\selectfont
\allowdisplaybreaks[4]

\begin{document}

\title{Eigenvalue inequalities and  three-term asymptotic formulas of the heat traces for the Lam\'{e} operator and Stokes operator}
\fancyhead[C]{\footnotesize\textit{G. Q. Liu \ /  \ Eigenvalue inequalities and  three-term asymptotic formulas  of the heat traces}}

\begin{frontmatter}

\author[1]{Genqian Liu\corref{c}}
\ead{liuswgq@163.com}
 \cortext[c]{Corresponding author.}
 \address[1]{School of Mathematics and Statistics, Beijing Institute of Technology, Beijing 100081,  P. R. China}

\date{}

\begin{abstract}   This paper is devoted to establish the most essential connections of the eigenvalue problems for the Laplace operator,   Lam\'{e} operator,  Stokes operator, buckling operator and clamped plate operator. We show that the $k$-th Stokes (respectively, Laplace)  eigenvalue is the limit of the $k$-th  Lam\'{e}  eigenvalue  for the Dirichlet or traction boundary condition as the Lam\'{e} coefficient $\lambda$  tends to $+\infty$ (respectively,  to $-\mu$).  
Furthermore, we establish the eigenvalue inequalities  and  three-term asymptotic formulas of the heat traces  for the Laplace operator, the Lam\'{e}  operator, the Stokes  operator  and buckling operator with the Dirichlet and traction  boundary conditions. 
\end{abstract}


\end{frontmatter}

\section{Introduction}

Revealing the relationship between the eigenvalues of the most basic physical operators,  establishing the heat trace asymptotics for partial differential operators have been the subject of extensive research for over a century. It has attracted the attention of many outstanding mathematicians and physicists. Beyond the  wonderful limit laws and  beautiful asymptotic formulas that are intimately related to the geometric properties of the domain and its boundary, a sustaining force has been its important role in mathematics, mechanics and theoretical physics (see, for example,  \cite{AsB}, \cite{BG}, \cite{Chav1}, \cite{Cha2}, \cite{CLN}, \cite{CH},  \cite{ES}, \cite{BG},  \cite{Gil2}, \cite{Gru}, \cite{Ho4}, \cite{Iv}, \cite{Kac}, \cite{Lap},   \cite{MS-67},  \cite{Pa1},  \cite{Ple1}, \cite{Ple2}, \cite{Po},  \cite{Sar},  \cite{ShY},  \cite{Ste}, \cite{Ta2},  \cite{Wei}, \cite{Wey1},  \cite{Wey2},  etc). 

  Let $(\Omega, g )$ be a smooth compact Riemannian
manifold of dimension $n\ge 2$ with smooth boundary $\partial \Omega$. 
The Lam\'{e} operator ${L}_{\lambda, \mu}$ of linear elasticity  acting on vector fields $\mathbf{u}=(u^1,\cdots, u^n)$ can be written as   (see \cite{Liu-19}, \cite{Liu-23} and \cite{Liu2})
 \begin{eqnarray} \label{1-1}  {L}_{\lambda,\mu} \mathbf{u}\!\!& \!=\!\!&\! 2\mu \,{\mbox{Def}}^*\, {\mbox{Def}}\; \mathbf{u}-\lambda\, \mbox{grad} \,\mbox{div}\; \mathbf{u} \;\;\;\,\;\;\mbox{in}\;\; \Omega,  
  \end{eqnarray}
  where  $\,\mbox{Def}\, \mathbf{u}= \frac{1}{2} (\nabla \mathbf{u} +(\nabla \mathbf{u})^T)$ (see, p.$\,$464 of \cite{Ta1}), $\,\,(\nabla \mathbf{u})^T$ is the transpose of $\nabla \mathbf{u}$,   $\;\mbox{Def}^*$ is the adjoint  of deformation operator $\mbox{Def}$ (see,  \S2 of \cite{Liu-19}, or   Exercise 16 on p$\,$153 of \cite{Ta1}),  and  $\mbox{div}$ and $\mbox{grad}$ are the usual divergence and gradient operators. Here  $\lambda$ and $\mu$  
 are real constants known as Lam\'{e}  parameters, assumed to satisfy the conditions
  \begin{eqnarray} \label{24.6.23-1}  \mu>0, \;  \, \lambda+2\mu>0, \end{eqnarray} 
which guarantee strong  ellipticity of ${L}_{\lambda,\mu}$, see, e.g., \cite{ADN2}, \cite{CiMa},  \cite{Hah}, \cite{MaHu}, \cite{McL},  \cite{MHNZ}, \cite{LiQin} or \cite{Liu-39}.  
  Equivalently, $L_{\lambda, \mu}$ can be written as 
  \begin{eqnarray} \label{24.4.12-9}  L_{\lambda, \mu} \mathbf{u}= 
    \mu \, \big( \nabla^* \nabla \mathbf{u}    - \mbox{Ric} (\mathbf{u})\big)
   -(\lambda+\mu)\, \mbox{grad}\;\mbox{div}\, \mathbf{u}, \end{eqnarray}   
   where $(\nabla^* \nabla \mathbf{u})^j:=- \nabla_k \nabla^k u^j $ is the Bochner Laplacian (see (2.11) of \cite{Liu2}),
      $\;\nabla_k u^j=u^j_{\;\;\;;\,k}\!:= \frac{\partial u^j}{\partial x_k} + \Gamma^{j}_{kl} u^l$
     are the components of the covariant derivative  of the vector field $\mathbf{u}$,   $\;\nabla^k= g^{kl}\nabla_l$, and    \begin{eqnarray} \label{18/12/22} \mbox{Ric} \,(\mathbf{u})= \big(R^{\,1}_{l} u^l,   R^{\,2}_{l} u^l, \cdots,  R^{\,n}_{l} u^l\big)\end{eqnarray} denotes the action of Ricci tensor $\mbox{R}_l^{\;j}:= R^{k\,\,j}_{\,lk}$ on $\mathbf{u}$. 
     Here and further on we adopt the Einstein summation convention over repeated indices.  
More precisely, ${L}_{\lambda,\mu}$ can be expressed as
 \begin{eqnarray} \label{24.4.9-1} \big ({L}_{\lambda,\mu} \mathbf{u}\big)^{j} :=  - \mu\, \Big(
\nabla_k\nabla^k u^j + \mbox{Ric}^{\,\,j}_{k} u^k\Big) -(\lambda+\mu)\nabla^j\nabla_k u^k,  \;\;\;\;\, j=1, \cdots, n.  \end{eqnarray}  
For  the Lam\'{e} operator ${L}_{\lambda,\mu}$,  the  Dirichlet boundary condition is  $\mathbf{u}\big|_{\partial \Omega}=0$, and  the traction (i.e., free or Neumann) boundary condition is  $\mathcal{T}_{\lambda,\mu} \mathbf{u}\big|_{\partial \Omega}=0$,  where $\mathcal{T}_{\lambda,\mu}$ is defined by 
\begin{eqnarray}  \label{24.4.9-2}  \mathcal{T}_{\lambda,\mu} \mathbf{u} :=\!\!\!&\!\!\!& \!\! -2\mu \,(\mbox{Def}\; \mathbf{u})^\# \boldsymbol{\nu}-\lambda\, (\mbox{div}\; \mathbf{u})\nu \;\;\; \mbox{on}\;\; \partial \Omega, \end{eqnarray} 
where $ \;\#$ is the sharp operator (for a tensor) by raising index,  and $\boldsymbol{\nu}=(\nu^1, \cdots, \nu^n)$ is the exterior unit normal vector to the boundary $\partial \Omega$.  
Clearly,  (\ref{24.4.9-2}) is just   
 \begin{eqnarray} \label{24.4.10-6} \big(\mathcal{T}_{\lambda,\mu} \mathbf{u}\big)^j = -\lambda \nu^j \nabla_k u^k -\mu \big( \nu^k \nabla_k u^j +\nu_k \nabla^j u^k\big) \;\;\;   \mbox{on} \;\; \partial \Omega, \,\;\;\; j=1,\cdots, n. \end{eqnarray} 
  In physics,   the elastic wave equations $\frac{\partial^2 \mathbf{u}}{\partial t^2}=-{L}_{\lambda,\mu} \mathbf{u}$  describe the propagations of waves in an isotropic homogeneous elastic medium,  and the constants $\sqrt{\lambda+2\mu}$  and  $\sqrt{\mu}\,$ are the velocities of longitudinal and transverse waves, respectively.

Set \begin{eqnarray} \label{24.4.11-1}  \mathcal{V} =\big\{\mathbf{u}\in C_0^\infty (\Omega, T\Omega)\,\big|\, \mbox{div}\,\mathbf{u}=0\big\}.\end{eqnarray} Then \begin{eqnarray} \label{24.4.11} J^{\,k} = \mbox{closure of } \mathcal{V} \;\, \mbox{in}\;\,  H^k (\Omega, T\Omega), \;\;\;\; k=0,1.\end{eqnarray} Let $P$ be the orthogonal projection of $L^2(\Omega,T\Omega)$ onto $J^0$.
The operator $S_\mu \mathbf{u}:=-\mu P\Delta_g  \mathbf{u}$ with domain $\mathcal{D} (S_\mu) = \{\mathbf{u}\in J^0 \cap H^2 (\Omega,T\Omega) \}$  is called the Stokes operator, where $\Delta_g$ is the Beltrami-Laplace operator and $\mu>0$ is the viscosity (constant) coefficient.
 The Stokes operator  $S_\mu$ can equivalently  be written as   (see,  Theorem 2.1 of \cite{Liu4}): 
  \begin{eqnarray} \label{24.4.10-4} \mathcal{S}_{\mu} \mathbf{u}:=- \mbox{div}\, (\sigma_\mu (\mathbf{u}, p)^\# ,\, \;\;\, \mbox{div}\, \mathbf{u}=0,\end{eqnarray} 
  where  $\mathbf{u} = (u^1,\cdots, u^n)$ is  the velocity vector field,  
  $\,\sigma_\mu (\mathbf{u},p)^\#= 2\mu \,\big(\mbox{Def}\, (\mathbf{u})\big)^\#-p \mathbf{I}_n $ is a tensor of field of type $(1,1)$,   and $p$ is the  pressure.  
 Precisely, the Stokes operator $S_\mu$ can be represented  as (see \cite{Liu4})
 \begin{eqnarray} \label{2022.6.24-1}  S_\mu\mathbf{u}= \mu \,\big( \nabla^*\nabla \mathbf{u} -\mbox{Ric}\, (\mathbf{u}) \big)  + \nabla_g p \;\;\, \mbox{with}\;\; \, \mbox{div}\, \mathbf{u}=0\,\;\;\, \mbox{in}\;\, \Omega, \end{eqnarray}
   where  $-\nabla^*\nabla \mathbf{u}\,$  and $\,\mbox{Ric}\,(\mathbf{u})$  are defined as before.
In physics,   a fluid flow  obeying  the stationary Stokes equations
   \begin{eqnarray} \label{2022.3.20-1}  \left\{\!\begin{array}{ll} - \mbox{div}\; (\sigma_\mu (\mathbf{u}, p))^\sharp =0 \;\; & \mbox{in} \;\; \Omega,\\
   \mbox{div}\; \mathbf{u}=0 \;\;& \mbox{in}\;\; \Omega,\end{array} \right. \end{eqnarray}
 is called the Stokes flow (i.e., creeping flow). For the Stokes operator, the Dirichlet boundary condition  is $\mathbf{u}\big|_{\partial \Omega}=0\,$ (i.e., $\mathcal{D}(S_\mu)=J^0\cap H_0^1(\Omega, T\Omega)\cap H^2(\Omega, T\Omega)$. Note that in this case, $J^1=\big\{\mathbf{u}\in H_0^1(\Omega) \,\big| \,\mbox{div}\, \mathbf{u}=0\;\,\mbox{in}\;\, \Omega\big\}$, see Lemma 5.5 on p.$\,$575 of \cite{Ta3}), and the Cauchy force boundary condition is $\big(\sigma_\mu (\mathbf{u}, p)\big)^\#\boldsymbol{\nu} \,\big|_{\partial \Omega}=0$.

 We consider the following classical eigenvalue problems:
 \begin{eqnarray}\label{4.11_1}\!\!\!\!\!\!\!\!\!\!\!\!\!\!\!\! \!\!\!  \!\!\!\! \!\!  \!\!\!\! (DL)  \quad \quad \quad \;\;\left\{ \begin{array}{ll}  L_{\lambda,\mu} \mathbf{u} = \tau^{\begin{scriptsize}\mbox{(D)}\end{scriptsize}} \,\mathbf{u}  \;\;\, &\mbox{in}\;\; \Omega,\\
   \mathbf{u}=0 \;\;\; & \mbox{on}\;\; \partial \Omega;\end{array} 
   \right.\; \quad \;\;\quad\end{eqnarray} 
  \begin{eqnarray}\label{4.11_2}\!\!\!\!\!\!\!\! \!\!\!\!\!\!\!\!\!\!\!\! \!\!\!\! \! \!\!\!\!\! (TL)\;\;\quad \quad \quad \left\{ \begin{array}{ll}  L_{\lambda,\mu} \mathbf{u} = \tau^{\begin{scriptsize}\mbox{(T)}\end{scriptsize}} \,\mathbf{u}  \;\;\, &\mbox{in}\;\; \Omega,\\
  \mathcal{T}_{\lambda,\mu}  \mathbf{u}=0 \;\;\; & \mbox{on}\;\; \partial \Omega;\end{array}  \right. \quad \;\; \quad\end{eqnarray} 
  \begin{eqnarray}\label{4.11_3} \!\!\!\!\!\!\!\!\!\!\!\!\!\!\!\! \!\!\!\! \!\!\!\! \!\!  \!\!\!\! (DS)\;\;\quad \quad \quad \left\{ \begin{array}{ll}  S_\mu \mathbf{v} = \varsigma^{\begin{scriptsize}\mbox{(D)}\end{scriptsize}}  \mathbf{v}  \;\;\, &\mbox{in}\;\; \Omega,\\
 \mbox{div}\; \mathbf{v}=0 \;\;\, &\mbox{in}\;\, \Omega,\\
 \mathbf{v}=0
\;\;\; & \mbox{on}\;\; \partial \Omega;\end{array}  \right. \,\quad\quad \;\; \,\end{eqnarray} 
 \begin{eqnarray}\label{4.11_4} \!\!\!\!\!\!\!\! \!\!\!\!\!\!\!\! \!\!\!\!\!\!\!\!\! \! \!\!\!\! \!  (CS)\;\;\quad \quad \quad \left\{ \begin{array}{ll}  S_\mu \mathbf{v} = \varsigma^{\begin{scriptsize}\mbox{(C)}\end{scriptsize}}  \mathbf{v}  \;\;\, &\mbox{in}\;\; \Omega,\\
 \mbox{div}\; \mathbf{v}=0 \;\;\, &\mbox{in}\;\, \Omega,\\
      \big( \sigma_\mu (\mathbf{v}, p)\big)^\#\,\boldsymbol{\nu} =0
    \;\;\; & \mbox{on}\;\; \partial \Omega;\end{array}  \right. \;\, \end{eqnarray} 
 \begin{eqnarray}\label{4.11_5} (DB)\;\;\quad \quad \quad \left\{ \begin{array}{ll}  \mu\big(\nabla^*\nabla \mathbf{w} - \mbox{Ric} (\mathbf{w}) \big)  = \theta^{\,\begin{scriptsize}\mbox{(D)}\end{scriptsize}} \mathbf{w}  \;\;\, &\mbox{in}\;\; \Omega,\\
    \mathbf{w}=0 \;\;\; & \mbox{on}\;\; \partial \Omega;\end{array}  \right. \quad\, \;\; \end{eqnarray} 
 \begin{eqnarray} \label{4.11_6}  (TB)\;\;\quad \quad \quad \left\{ \begin{array}{ll}  \mu\big(\nabla^*\nabla \mathbf{w} - \mbox{Ric} (\mathbf{w}) \big)   = \theta^{\,\begin{scriptsize}\mbox{(T)}\end{scriptsize}} \mathbf{w}  \;\;\, &\mbox{in}\;\; \Omega,\\
   \mathcal{T}_{-\mu,\mu} \mathbf{w}=0 \;\;\; & \mbox{on}\;\; \partial \Omega;\end{array}  \right. \;\quad \;\;\end{eqnarray}

(\ref{4.11_1})  (respectively, (\ref{4.11_2})) is the Lam\'{e} eigenvalue problem  with Dirichlet (respectively, traction) boundary condition; 
 (\ref{4.11_3})  (respectively, (\ref{4.11_4})) is the Stokes eigenvalue problem with Dirichlet (respectively, Cauchy force) boundary condition;  (\ref{4.11_5})  (respectively, (\ref{4.11_6})) is the Laplace eigenvalue problem  with Dirichlet (respectively, traction) boundary condition.  In each of these cases, the spectrum is discrete and we arrange the eigenvalues in
non-decreasing order (repeated according to multiplicity)
  \begin{eqnarray*} && 0< \tau_1^{\begin{scriptsize} \mbox{(D)}\end{scriptsize} } \le \tau_2^{\begin{scriptsize} \mbox{(D)}\end{scriptsize}} \le \cdots \le  \tau_k^{\begin{scriptsize} \mbox{(D)}\end{scriptsize}}\le \cdots;\\
&&  0\le  \tau_1^{\begin{scriptsize} \mbox{(T)}\end{scriptsize} } \le \tau_2^{\begin{scriptsize} \mbox{(T)}\end{scriptsize} } \le \cdots \le  \tau_k^{\begin{scriptsize} \mbox{(T)}\end{scriptsize} } \le \cdots;\\
&&  0<  \varsigma_1^{\begin{scriptsize} \mbox{(D)}\end{scriptsize} }  \le \varsigma_2^{\begin{scriptsize} \mbox{(D)}\end{scriptsize} } \le \cdots \le  \varsigma_k^{\begin{scriptsize} \mbox{(D)}\end{scriptsize} } \le \cdots;\\
&&  0\le \varsigma_1^{\begin{scriptsize} \mbox{(C)}\end{scriptsize} }  \le \varsigma_2^{\begin{scriptsize} \mbox{(C)}\end{scriptsize} }  \le \cdots \le  \varsigma_k^{\begin{scriptsize} \mbox{(C)}\end{scriptsize} } \le \cdots;\\
&&  0<  \theta_1^{\begin{scriptsize} \mbox{(D)}\end{scriptsize} } \le \theta_2^{\begin{scriptsize} \mbox{(D)}\end{scriptsize} } \le \cdots \le  \theta_k^{\begin{scriptsize} \mbox{(D)}\end{scriptsize} }\le \cdots;\\
&&   0\le \theta_1^{\begin{scriptsize} \mbox{(T)}\end{scriptsize} } \le \theta_2^{\begin{scriptsize} \mbox{(T)}\end{scriptsize} } \le \cdots \le  \theta_k^{\begin{scriptsize} \mbox{(T)}\end{scriptsize} }\le \cdots.
   \end{eqnarray*}  
The corresponding eigenfunctions are expressed as $\mathbf{u}_1^{\begin{scriptsize} \mbox{(D)}\end{scriptsize} }, \mathbf{u}_2^{\begin{scriptsize} \mbox{(D)}\end{scriptsize} }, \mathbf{u}_3^{\begin{scriptsize} \mbox{(D)}\end{scriptsize} }, \cdots$;
$\;\,\mathbf{u}_1^{\begin{scriptsize} \mbox{(T)}\end{scriptsize} }, \mathbf{u}_2^{\begin{scriptsize} \mbox{(T)}\end{scriptsize} }, \mathbf{u}_3^{\begin{scriptsize} \mbox{(T)}\end{scriptsize} }, \cdots$;
   $\;\,\mathbf{v}_1^{\begin{scriptsize} \mbox{(D)}\end{scriptsize} }, \mathbf{v}_2^{\begin{scriptsize} \mbox{(D)}\end{scriptsize} }, \mathbf{v}_3^{\begin{scriptsize} \mbox{(D)}\end{scriptsize} }, \cdots$; 
    $\,\;\mathbf{v}_1^{\begin{scriptsize} \mbox{(C)}\end{scriptsize} }, \mathbf{v}_2^{\begin{scriptsize} \mbox{(C)}\end{scriptsize} }, \mathbf{v}_3^{\begin{scriptsize} \mbox{(C)}\end{scriptsize} }, \cdots$; 
      $\,\;\mathbf{w}_1^{\begin{scriptsize} \mbox{(D)}\end{scriptsize} }, \mathbf{w}_2^{\begin{scriptsize} \mbox{(D)}\end{scriptsize} }, \mathbf{w}_3^{\begin{scriptsize} \mbox{(D)}\end{scriptsize} }, \cdots$;  
  $\,\;\mathbf{w}_1^{\begin{scriptsize} \mbox{(T)}\end{scriptsize} }, \mathbf{w}_2^{\begin{scriptsize} \mbox{(T)}\end{scriptsize} }, \mathbf{w}_3^{\begin{scriptsize} \mbox{(T)}\end{scriptsize} }, \cdots$.

\vskip 0.38 true cm

Our main results in this paper are following:

 \vskip 0.22 true cm

 \noindent{\bf Theorem 1.1.} \ {\it  Let $(\Omega,g)$ be an $n$-dimensional smooth compact  Riemannian manifold with smooth boundary $\partial \Omega$, $\,(n\ge 2)$.   
 Then for given  Lam\'{e} coefficient  $\mu>0$,
 the spectrum $\sigma_D (L_{\lambda, \mu})$ (respectively, $\sigma_T(L_{\lambda, \mu})$)  of $L_{\lambda, \mu}$ with Dirichlet (respectively, traction) boundary condition is discrete for any $\lambda\in (-2\mu, +\infty)$. \ The eigenvalues $\tau^{\begin{scriptsize}\mbox{(D)}\end{scriptsize}}_k (\lambda)$ (respectively,  $\tau_k^{\begin{scriptsize}\mbox{(T)}\end{scriptsize}} (\lambda)$) of  
 $L_{\lambda,\mu}$ ordered in a non-decreasing sequence, is continuous  non-decreasing functions of $\lambda$. }

 \vskip 0.39 true cm 
 
 \noindent{\bf Theorem 1.2.} \ {\it  Let $(\Omega,g)$ be an $n$-dimensional smooth compact  Riemannian manifold with smooth boundary $\partial \Omega$, $\,(n\ge 2)$.   
 Then for given  Lam\'{e} coefficient  $\mu>0$ and any $k\ge 1$, 
  the Lam\'{e} eigenvalues $\tau_k^{\begin{tiny}(D)/(T)\end{tiny}}(\lambda)\,$   tend to  the Stokes eigenvalues $\varsigma_k^{\begin{tiny}(D)/(C)\end{tiny}}$ (respectively, the Laplace eigenvalues $\theta_k^{\begin{tiny}(D)/(T)\end{tiny}}$)  as the Lam\'{e} constant $\lambda\to +\infty$ (respectively, $\lambda\to -\tau$),  and the corresponding  Lam\'{e} eigenspaces tend to the corresponding Stokes eigenspaces (respectively, Laplace eigenspaces).}
 
 \vskip 0.38  true cm

 Combining the results of Theorem 1.1 and Theorem 1.2  with the mini-max formulas for the eigenvalues,
we have  the basic eigenvalue inequalities:

\vskip 0.29 true cm 

 \noindent{\bf Corollary 1.3}.  {\it Let $(\Omega,g)$ be an $n$-dimensional smooth compact  Riemannian manifold with smooth boundary $\partial \Omega$, $\,(n\ge 2)$.  Let $\tau_k^{\begin{tiny}(D)/(T)\end{tiny}}$,  $\varsigma_k^{\begin{tiny}(D)/(C)\end{tiny}}$, $\theta_k^{\begin{tiny}(D)/(T)\end{tiny}}$  be the
$k$-th eigenvalue with repetition according to the multiplicity of the self-adjoint operators
 $L_{\lambda,\mu}$,  $S_\mu$ and $B_\mu:=L_{-\mu,\mu}$ with corresponding boundary conditions, respectively. Then for any $\lambda\in (-\mu, +\infty)$ there holds
 \begin{eqnarray} \label{24.4.14-6}  \theta_k^{\begin{tiny}(D)\end{tiny}} \le \tau_k^{\begin{tiny}(D)\end{tiny}}(\lambda) \le \varsigma_k^{\begin{tiny}(D)\end{tiny}}, \;\;\;\,\;
 \theta_k^{\begin{tiny}(T)\end{tiny}} \le \tau_k^{\begin{tiny}(T)\end{tiny}}(\lambda) \le \varsigma_k^{\begin{tiny}(T)\end{tiny}}, \;\;\; \, k=1,2,3,\cdots.\end{eqnarray}  }
  
\vskip 0.26  true cm

We introduce the {\it partition function}, or  the {\it trace of the heat semigroup} for the Lam\'{e} operator with  Dirichlet (respectively,  traction) boundary condition, by
$\mathcal{Z}^{\begin{tiny}(D)\end{tiny}} (t) := \mbox{Tr}\; e^{-t L_{\lambda,\mu}^{\begin{tiny}(D)\end{tiny}}} = \sum_{k=1}^\infty e^{-t \tau_k^{\begin{tiny}(D)\end{tiny}}}$ (respectively, $\mathcal{Z}^{\begin{tiny}(T)\end{tiny}} (t) := \mbox{Tr}\; e^{-t L_{\lambda,\mu}^{\begin{tiny}(T)\end{tiny}}} = \sum_{k=1}^\infty e^{-t \tau_k^{\begin{tiny}(T)\end{tiny}}}$) defined for $t>0$ and monotone decreasing in $t$.

\vskip 0.36 true cm

\noindent{\bf Theorem 1.4.} \ {\it Let $(\Omega,g)$ be an $n$-dimensional compact smooth  Riemannian manifold with smooth boundary $\partial \Omega$.  Assume that  the Lam\'{e} coefficients  satisfy $\mu>0$ and $\lambda+2\mu> 0$.  Let $0< \tau_1^{\begin{tiny}(D)\end{tiny}}< \tau_2^{\begin{tiny}(D)\end{tiny}} \le \tau^{\begin{tiny}(D)\end{tiny}}_3\le \cdots \le \tau_k^{\begin{tiny}(D)\end{tiny}} \le \cdots$ (respectively, $0\le \tau_1^{\begin{tiny}(T)\end{tiny}} \le \tau_2^{\begin{tiny}(T)\end{tiny}} \le \tau_3^{\begin{tiny}(T)\end{tiny}} \le \cdots \le \tau_k^{\begin{tiny}(T)\end{tiny}} \le \cdots $) be the eigenvalues of the Lam\'{e} operator $L_{\lambda,\mu}$  with respect to the  Dirichlet (respectively, traction) boundary condition. Then
\begin{eqnarray}  \left. \begin{array}{ll} & \mathcal{Z}^{\begin{tiny}(D)/(T)\end{tiny}}(t) =  \sum_{k=1}^\infty e^{-t \tau_k^{\begin{tiny}(D)/(T)\end{tiny}}} = a_0 \,t^{-n/2} {\mbox{Vol}_n}(\Omega) + a_1^{\begin{tiny}(D)/(T)\end{tiny}}\, t^{(1-n)/2} {\mbox{Vol}_{n-1}}(\partial\Omega)\\
[3mm]&\qquad \quad \;\; \quad\, \;\; +a_2\, t^{(2-n)/2} +o(t^{(2-n)/2})  \quad \mbox{as}\;\; t\to 0^+, \;\; \;\,   \end{array} \right. \;\;\quad \;\; \mbox{for}\;\; n\ge 2,\nonumber \end{eqnarray}
where \begin{eqnarray*} \!\!\!\!\! &\!\!\!\!\! &\!\!\!\!\!\!\!\!\!\!\!\!\!\!a_0= \frac{n-1}{(4\pi \mu )^{n/2}}\!
 + \! \frac{1}{(4\pi (\lambda+2\mu) )^{n/2}}, \;\;\;\, a_1^{\begin{tiny}(D)\end{tiny}}= -\frac{1}{4} \bigg[  \frac{n-1}{(4\pi \mu\,)^{(n-1)/2}}
 +  \frac{1}{(4\pi (\lambda+2\mu) )^{(n-1)/2}}\bigg],\\
\!\!\!\!\! &\!\!\!\!\! &\!\!\!\!\!\!\!\!\!\!\!\!\!\! a_1^{\begin{tiny}(T)\end{tiny}}=  \frac{1}{4} \bigg[  \frac{n-1}{(4\pi \mu\,)^{(n-1)/2}}
 +  \frac{1}{(4\pi (\lambda+2\mu) )^{(n-1)/2}}\bigg],\\
 \!\!\!\!\! &\!\!\!\!\!& \!\!\!\!\!\!\!\!\! \!\!\!\! a_2=\frac{1}{6 (4\pi)^{n/2}} \bigg[\bigg( \frac{1}{ (\lambda+2\mu)^{(n-2)/2}} +\frac{n-7}{ \mu^{(n-2)/2}} + \frac{12\mu}{n} \Big( \frac{1}{(\lambda+2\mu)^{n/2}} +\frac{n-1}{\mu^{n/2}} \Big) \bigg)\! \int_{\Omega} \! R(x)\, dV\\  && \qquad \;\;\;\quad + 2 \,\bigg(\frac{1}{ (\lambda+2\mu)^{(n-2)/2}} +\frac{n-7}{ \mu^{(n-2)/2}}\bigg) \int_{\partial \Omega}H(x) \,ds\bigg] .\nonumber\end{eqnarray*} 
 Here ${\mbox{Vol}}_{n}(\Omega)$ denotes the $n$-dimensional volume of $\, \Omega$,  ${\mbox{Vol}}_{n-1} (\partial\Omega)$  is the $(n-1)$-dimensional volume of $\partial \Omega$,  $R(x)$ denotes the scalar curvature at $x\in \Omega$, and $H(x)$  is the mean curvature of the boundary $\partial \Omega$ at $x\in \partial \Omega$. }

\vskip 0.36  true cm

From Theorem 1.3 and Theorem 1.4 (i.e., as $\lambda\to +\infty$), we immediately obtain: 

\vskip 0.22  true cm 

\noindent{\bf  Corollary 1.5.} \ {\it Let $(\Omega,g)$ be an $n$-dimensional compact smooth  Riemannian manifold with smooth boundary $\partial \Omega$.  Assume that  the viscosity coefficient  satisfy $\mu>0$.  Let $0< \varsigma_1^{(D)}< \varsigma_2^{(D)} \le \varsigma^{(D)}_3\le \cdots \le \varsigma_k^{(D)} \le \cdots$ (respectively, $0\le \varsigma_1^{(C)} < \varsigma_2^{C)} \le \varsigma_3^{(C)}\le \cdots \le \varsigma_k^{(C)} \le \cdots $) be the eigenvalues of the Stokes operator $S_{\mu}$  with respect to the  Dirichlet (respectively, Cauchy force) boundary condition. Then
\begin{eqnarray}  \left. \begin{array}{ll} & \mathcal{Z}^{(D)/(C)}(t) =  \sum_{k=1}^\infty e^{-t\, \varsigma_k^{\begin{scriptsize}(D)/(C)\end{scriptsize}}} = b_0 \,t^{-n/2} {\mbox{Vol}_n}(\Omega) + b_1^{\begin{scriptsize}(D)/(C)\end{scriptsize}}\, t^{(1-n)/2} {\mbox{Vol}_{n-1}}(\partial\Omega)\\
[3mm]&\qquad \quad \;\; +b_2 t^{(2-n)/2} +O(t^{(3-n)/2})  \quad \mbox{as}\;\; t\to 0^+, \;\; \;\,   \end{array} \right. \;\;\quad \;\; \mbox{for}\;\; n\ge 2,\nonumber \end{eqnarray}
where \begin{eqnarray*} \!\!\!\!\! &\!\!\!\!\! &\!\!\!\!\!\!\!\!\!\!\!\!\!\! b_0= \frac{n-1}{(4\pi \mu )^{n/2}}, \;\; \;\;\;\, \, b_1^{(D)}= -   \frac{n-1}{4(4\pi \mu\,)^{(n-1)/2}}, \;\;\, \;\;\, b_1^{(C)}=    \frac{n-1}{4(4\pi \mu\,)^{(n-1)/2}}\\
\!\!\!\!\! &\!\!\!\!\!& \!\!\!\!\!\!\!\!\! \!\!\!\! b_2=\frac{1}{6 (4\pi)^{n/2}} \bigg[\bigg(\beta_n +\frac{n-7}{ \mu^{(n-2)/2}} + \frac{12 (n-1)}{n\,\mu^{(n-2)/2}} \bigg)\! \int_{\Omega} \! R(x)\, dV +
2 \,\bigg(\,\beta_n +\frac{n-7}{\mu^{(n-2)/2}}\bigg) \int_{\partial \Omega} H(x)\, ds\bigg],\nonumber\end{eqnarray*} 
 and \begin{eqnarray*} \left\{ \begin{array}{ll} \beta_n=0 \;\; &\mbox{when}\;\; n\ge 3,\\
 \beta_n =1 \;\; & \mbox{when}\;\;  n=2.\end{array}\right.\end{eqnarray*}}

\vskip 0.29  true cm 

\noindent{\bf  Remark 1.6.} \ {\it  For the Stokes eigenvalues with Dirichlet boundary condition in a bounded domain $\Omega\subset \mathbb{R}^n$,  ($n\ge 2$),  the corresponding two-term asymptotic formula have been established  in \cite{Liu3}  by a completely different method, in which we used the celebrated semigroup theory  by   Abe and Giga for the Stokes equations (see, \cite{AGi} and \cite{Giga}) and a (pseudodifferential operator) representation of inverse of the Stokes operator by Kozhevnikov \cite{Koz}. }

\vskip 0.46  true cm

\noindent{\bf Corollary 1.7.} \ {\it Let $(\Omega,g)$ be an $n$-dimensional compact smooth  Riemannian manifold with smooth boundary $\partial \Omega$. Let $0< \theta_1^{(D)}< \theta_2^{(D)} \le  \cdots \le \theta_k^{(D)} \le \cdots$ (respectively, $0\le \theta_1^{(T)} \le \theta_2^{(T)} \le\cdots \le \theta_k^{(T)} \le \cdots $) be the eigenvalues of the Laplace operator $  \mu \, \big( \nabla^* \nabla \mathbf{u}    - \mbox{Ric} (\mathbf{u})\big)$  with respect to the  Dirichlet (respectively, traction) boundary condition. Then
\begin{eqnarray}  \left. \begin{array}{ll} & \mathcal{Z}^{(D)/(T)}(t) =  \sum_{k=1}^\infty e^{-t \theta_k^{(D)/(T)}} = c_0 \,t^{-n/2} {\mbox{Vol}_n}(\Omega) + c_1^{(D)/(T)}\, t^{(1-n)/2} {\mbox{Vol}_{n-1}}(\partial\Omega)\\
[3mm]&\qquad \quad \;\; +c_2\, t^{(2-n)/2} +O(t^{(3-n)/2})  \quad \mbox{as}\;\; t\to 0^+, \;\; \;\,   \end{array} \right. \;\;\quad \;\; \mbox{for}\;\; n\ge 2,\nonumber \end{eqnarray}
where \begin{eqnarray*} \!\!\!\!\! &\!\!\!\!\! &\!\!\!\!\!\!\!\!\!\!\!\!\!\! c_0= \frac{n}{(4\pi \mu )^{n/2}}, \;\;\;\, c_1^{(D)}= - \frac{n}{4(4\pi \mu\,)^{(n-1)/2}}, \;\;\; \;\, c_1^{(T)}=  \frac{n}{4(4\pi \mu\,)^{(n-1)/2}},\\
\!\!\!\!\! &\!\!\!\!\!& \!\!\!\!\!\!\!\!\! \!\!\!\! c_2=\frac{1}{6 (4\pi)^{n/2}} \bigg[\frac{n+6}{ \mu^{(n-2)/2}}  \int_{\Omega} \! R(x)\, dV+
\frac{2(n-6)}{\mu^{(n-2)/2}} \int_{\partial \Omega} H(x)\,ds \bigg].\nonumber\end{eqnarray*} 
 }

\vskip 0.32  true cm

Let us remark that all the above  conclusions hold if $\Omega$ is a bounded domain with smooth boundary in the Euclidean space $\mathbb{R}^n$.   Now,  let $\Omega$ be a bounded domain with smooth boundary in the Euclidean plane $\mathbb{R}^2$, and let $\mu$ be a positive constant.  Denote by $\Lambda^{(D)}$ the buckling eigenvalues with Dirichlet boundary condition
 \begin{eqnarray} \label{24.6.7-2}  \left\{ \begin{array}{ll}  \mu \Delta^2 \psi +\Lambda \Delta \psi =0 \;\;\;& \mbox{in}\;\;\, \Omega,\\
 u=\frac{\partial u}{\partial \nu}=0 \;\;\, & \mbox{on}\;\; \, \partial \Omega,  \end{array}\right. \end{eqnarray}  
and   by $\Lambda^{(C)}$ the buckling eigenvalue with Cauchy force  boundary
 \begin{eqnarray} \left\{ \begin{array}{ll}  \mu \Delta^2 \psi +\Lambda \Delta \psi =0 \;\;\;&\mbox{in}\;\;\, \Omega,\\
\psi=0 \;\;\, &\mbox{on}\;\; \, \partial \Omega,\\
(B{\psi})\boldsymbol{\nu}=0\;\,  &\mbox{on} \;\;\partial \Omega, \end{array} \right.\end{eqnarray} 
where \begin{eqnarray*} (B\psi) \boldsymbol{\nu}=\left\{ \mu \begin{pmatrix}  -2 \partial_1\partial_2 \psi_k  &  \partial_1^2 \psi_k -\partial_2^2 \psi_k\\
\partial_1^2 \psi_k -\partial_2^2 \psi_k &  2\partial_1\partial_2 \psi_k \end{pmatrix} - \mu \sqrt{ (\partial_1^2 \psi_k -\partial_2^2\psi_k)^2 +(\partial_1\partial_2 \psi_k)^2} \;\, \mathbf{I}_2 \right\}\boldsymbol{\nu} \;\;\;\, \mbox{on}\;\;\partial \Omega. \end{eqnarray*}
 It is easy to verify that each of spectrum is discrete with 
 \begin{eqnarray*}  &&0< \Lambda_1^{(D)}< \Lambda_2^{(D)} \le \cdots \le \Lambda_k^{(D)} \le \cdots, \\
 && 0< \Lambda_1^{(C)}< \Lambda_2^{(C)} \le \cdots \le \Lambda_k^{(C)} \le \cdots, \end{eqnarray*}
each eigenvalue repeated according to its multiplicity.

\vskip 0.48  true cm  

\noindent{\bf  Theorem 1.8.} \ {\it  Let $\Omega\subset \mathbb{R}^2$ be a bounded domain with smooth boundary $\partial \Omega$. 
 Let $\Lambda_k^{(D)}$ (respectively, $\Lambda_k^{(C)}$) be the $k$-th  buckling eigenvalue with the Dirichlet (respectively, Cauchy force) boundary condition. Then the set $\{\Lambda_k^{(D)}\}_{k=1}^\infty$ (respectively, $\{\Lambda_k^{(C)}\}_{k=1}^\infty$) of all the buckling eigenvalues with Dirichlet (respectively, Cauchy force) boundary condition is the same as the set $\{\varsigma_k^{(D)}\}_{k=1}^\infty$ (respectively,  $\{\varsigma_k^{(C)}\}_{k=1}^\infty$) of all the Stokes eigenvalues with Dirichlet (respectively, Cauchy force) boundary condition. That is, 
 \begin{eqnarray} \label{24.4.23-1}  \Lambda_k^{(D)} =\varsigma_k^{(D)}, \;\;\;\;\, \Lambda_{k}^{(C)} =\varsigma_{k+1}^{(C)}, \;\;\; \;\; k=1,2,3, \cdots.\end{eqnarray}
  } 
 
 According to  Corollary 1.5 and $\varsigma_1=0$, we  get the following

\vskip 0.22  true cm 

\noindent{\bf  Corollary 1.9.} \ {\it Let $\Omega$ be a bounded domain with smooth boundary in the Euclidean plane $\mathbb{R}^2$.  Then, for the buckling eigenvalues $\{\Lambda_k^{(D)}\}_{k=1}^\infty$ and $\{\Lambda_k^{(C)}\}_{k=1}^\infty$, 
\begin{eqnarray}  \left. \begin{array}{ll} & \mathcal{Z}^{(D)/(C)}(t) =  \sum_{k=1}^\infty e^{-t \Lambda_k^{(D)/(C)}} = \frac{|\Omega|}{4\pi \mu \,t}  +d_1^{(D)/(C)}\,t^{-1/2}  +d_2    +O(1)  \quad \mbox{as}\;\; t\to 0^+,   \end{array} \right.  \nonumber \end{eqnarray}
where
  \begin{eqnarray*}   \!\!\!\!\! &\!\!\!\!\!& \!\!\!\!\!\!\!\!\! \!\!\!\! d_1^{(D)}=  -\, \frac{|\partial \Omega|}{4\sqrt{4\pi \mu}},  \;\;\, \;d_1^{(C)}=   \frac{|\partial \Omega|}{4\sqrt{4\pi \mu}}, \\
 \!\!\!\!\! &\!\!\!\!\!& \!\!\!\!\!\!\!\!\! \!\!\!\! d_2=-1-\,\frac{1}{3 \pi } \int_{\partial \Omega}H(x) \,ds.\nonumber\end{eqnarray*}  }

\vskip 0.22  true cm

Let $\Omega$ be a bounded domain with smooth boundary in the Euclidean space $\mathbb{R}^2$, and let $\mu$ be a positive constant.  Denote by $\{\Gamma_k^2, \phi_k)\}_{k=1}^\infty$ the clamped plate eigenvalues 
 \begin{eqnarray} \label{24.6.7-1}  \left\{ \begin{array}{ll} \mu^2 \Delta^2 \phi_k -\Gamma_k^2 \phi_k =0 \;\;\;& \mbox{in}\;\;\, \Omega,\\
\phi_k=\frac{\partial \phi_k}{\partial \nu}=0 \;\;\, & \mbox{on}\;\; \, \partial \Omega.  \end{array}\right. \end{eqnarray}

 For the Neumann eigenvalues $\{\Theta_k\}_{k=1}^\infty$ of the Laplace operator $\mu \,\Delta$,  the Dirichlet eigenvalues $\{\Xi_k\}_{k=1}^\infty$ of the Laplace operator $\mu \,\Delta$, the clamped plate eigenvalues $\{\Gamma_k^2\}_{k=1}^\infty$ and the buckling eigenvalues $\{\Lambda_k^{(D)}\}_{k=1}^\infty$ on a smooth compact Riemannian manifold $(\Omega, g)$ with smooth boundary $\partial \Omega$,   the present author in \cite{Liu-11} proved the following basic  inequalities:
\begin{eqnarray} \label{24.5.1-1} \Theta_k <\Xi_k < \Gamma_k <\Lambda_k^{(D)}, \;\;\;\;\; k=1,2,3,\cdots.\end{eqnarray} 
  The inequalities (\ref{24.5.1-1}) are sharp and they can not been further improved in a Riemannian manifold, see \cite{Liu-11}  for the detailed examples.  Thus, from the asymptotic formula of the Dirichlet eigenvalues  of the Laplace operator $\mu\, \Delta$ (see, \cite{MS-67}), Corollary 1.9 and the inequalities 
  $$\sum_{k=1}^\infty e^{-t \Lambda_k^{(D)}}< \sum_{k=1}^\infty e^{-t \Gamma_k}<\sum_{k=1}^\infty e^{-t \Xi_k},$$  we can obtain the following two-term asymptotic formulas:
  
  \vskip 0.2  true cm 

   \noindent{\bf  Corollary 1.10.} \ {\it Let $\Omega$ be a bounded domain with smooth boundary in the Euclidean plane $\mathbb{R}^2$.  Then, for the Dirichlet eigenvalues $\{\Theta_k\}_{k=1}^\infty$ of the  Laplace operator $\mu\, \Delta$, the buckling eigenvalues $\{\Lambda_k^{\begin{tiny}(D)\end{tiny}}\}_{k=1}^\infty$ (in  (\ref{24.6.7-2})) and the clamped plate eigenvalues $\{\Gamma_k\}_{k=1}^\infty$ (in (\ref{24.6.7-1})),  the following two-term asymptotic formulas holds:
\begin{eqnarray} \label{24.5.23-1} \left\{ \begin{array}{ll}    \sum_{k=1}^\infty e^{-t \Xi_k} \\
 \sum_{k=1}^\infty e^{-t \Gamma_k}  \\  \sum_{k=1}^\infty e^{-t \Lambda_k^{(D)}}\end{array}\right\} \;\large\sim  
 \frac{|\Omega|}{4\pi \mu \,t} -   \frac{|\partial \Omega|}{4\sqrt{4\pi \mu\,t}}   +O(1)  \quad \mbox{as}\;\; t\to 0^+.   \end{eqnarray}
 }

\vskip 0.46  true cm

The remainder of this paper is organized as follows.  In Section 2, we give the proof of Theorem 1.1 by using eigenvalue variational  principle. In Section 3, by a construct technique and by applying G{\aa}rding inequality  we prove Theorem 1.2.  Section 4  shows that the elastic Lam\'{e} operator is exactly the generalized Ahlfors Laplacian, and  provides the third term coefficient of its asymptotic expansion from a result of \cite{BGOP}.  In Section 5, we prove Theorem 1.4 by combining double operator (defined on the double closed  manifold), pseudodifferential operator technique and the result of Section 4.  In Section 6, we  discuss the  buckling and clamped plate eigenvalue problems and further establish the corresponding eigenvalue inequalities and asymptotic formulas
 from property of the Stokes operator in a  two-dimensional Euclidean bounded domain.  

\vskip 0.78 true cm 

\section{Monotonicity of the Lam\'{e} eigenvalues with respect to parameter}

\vskip 0.38 true cm

  \noindent  {\it Proof of Theorem 1.1.} \  By the maximum-minimum property of the Lam\'{e} eigenvalues, we have  that for any $k\ge 1$,
\begin{eqnarray*}  \label{24.4.20-1}\tau_k^{\begin{scriptsize}\mbox{(D)}\end{scriptsize}} (\lambda) \,=\, \max_{\substack{E\subset H_0^1(\Omega, T\Omega) \\  \mbox{\begin{footnotesize}dim\end{footnotesize}}\, E\,=\,k-1}} \, \;\, \min_{\substack{\boldsymbol{\phi}\in E^{\perp}}}\, \, \frac{ \int_{\Omega}  \Big[2\mu \big( \mbox{Def}\; \boldsymbol{\phi}, \mbox{Def}\; \boldsymbol{\phi}\big) + \lambda\,\big(\mbox{div}\; \boldsymbol{\phi}\big)^2 \Big] \,dV }{\int_\Omega \boldsymbol{\phi} \cdot  \boldsymbol{\phi} \; dV}\end{eqnarray*}
and  \begin{eqnarray*} \label{24.4.20-2} \tau_k^{\begin{scriptsize}\mbox{(T)}\end{scriptsize}}(\lambda) \,= \,\max_{\substack{F\subset H^1(\Omega, T\Omega) \\  \mbox{\begin{footnotesize}dim\end{footnotesize}}\, F\,=\,k-1}} \, \;\, \min_{\substack{\boldsymbol{\phi}\in F^{\perp}}}\, \, \frac{ \int_{\Omega}  \Big[2\mu \big( \mbox{Def}\; \boldsymbol{\phi},  \mbox{Def}\; \boldsymbol{\phi}\big) + \lambda\,\big(\mbox{div}\; \boldsymbol{\phi}\big)^2 \Big] \,dV }{\int_\Omega \boldsymbol{\phi}\cdot  \boldsymbol{\phi} \; dV},\end{eqnarray*}
where $E^{\perp}$  (respectively, $F^{\perp}$) is the orthogonal complement space  of $E$ (respectively, $F$) in $H_0^1(\Omega,T\Omega)$ (respectively, $H^1(\Omega, T\Omega)$) under the $L^2$-norm.
Or equivalently, there exists  a  sequence of eigenpairs  $\{(\tau_k^{\begin{scriptsize}\mbox{(D)}\end{scriptsize}}(\lambda), \mathbf{u}_k^{\begin{scriptsize}\mbox{(D)}\end{scriptsize}}(\lambda))\}_{k=1}^\infty$ (respectively, $\{(\tau_k^{\begin{scriptsize}\mbox{(T)}\end{scriptsize}}(\lambda), \mathbf{u}_k^{\begin{scriptsize}\mbox{(T)}\end{scriptsize}}(\lambda))\}_{k=1}^\infty$ such that 
\begin{eqnarray} \label{24.4.27-1} & & \tau_1^{\begin{scriptsize}\mbox{(D)}\end{scriptsize}}(\lambda) \, = \frac{ \int_{\Omega}  \Big[2\mu \big( \mbox{Def}\, \mathbf{u}^{\begin{scriptsize}\mbox{(D)}\end{scriptsize}}_1(\lambda),\mbox{Def}\, \mathbf{u}_1^{\begin{scriptsize}\mbox{(D)}\end{scriptsize}} (\lambda)\big)+ \lambda\,\big(\mbox{div}\, \mathbf{u}_1^{\begin{scriptsize}\mbox{(D)}\end{scriptsize}}(\lambda)\big)^2 \Big] \,dV }{\int_\Omega
\mathbf{u}_1^{\begin{scriptsize}\mbox{(D)}\end{scriptsize}}(\lambda)\cdot  \mathbf{u}_1^{\begin{scriptsize}\mbox{(D)}\end{scriptsize}}(\lambda)\; dV}\\
  & & \qquad \;\;\,\;\, =  \min_{\substack{\boldsymbol{\phi} \in H_0^1(\Omega, T\Omega)}} \frac{ \int_{\Omega}  \Big[2\mu \big( \mbox{Def}\, \boldsymbol{\phi},  \mbox{Def}\, \boldsymbol{\phi}\big)+ \lambda\,\big(\mbox{div}\, \boldsymbol{\phi}\big)^2 \Big] \,dV }{\int_\Omega \boldsymbol{\phi}\cdot \boldsymbol{\phi}\;  dV},\nonumber \\
 [4mm]&& \quad  \quad \cdots  \cdots \cdots \cdots \nonumber\\
 [4mm] \label{24.4.27-2}  &&   \tau_{k}^{\begin{scriptsize}\mbox{(D)}\end{scriptsize}}(\lambda) =\frac{ \int_{\Omega}  \Big[2\mu \big(\mbox{Def}\; \mathbf{u}^{\begin{scriptsize}\mbox{(D)}\end{scriptsize}}_k(\lambda),\mbox{Def}\, \mathbf{u}_k^{\begin{scriptsize}\mbox{(D)}\end{scriptsize}} (\lambda)\big)+ \lambda\,\big(\mbox{div}\; \mathbf{u}_k^{\begin{scriptsize}\mbox{(D)}\end{scriptsize}}(\lambda)\big)^2 \Big] \,dV }{\int_\Omega \mathbf{u}_k^{\begin{scriptsize}\mbox{(D)}\end{scriptsize}}(\lambda)\cdot \mathbf{u}_k^{\begin{scriptsize}\mbox{(D)}\end{scriptsize}}(\lambda)  \; dV} \\
 [1mm]  &&\qquad \;\;\, \; =  \min_{\substack{\boldsymbol{\phi} \in H_0^1(\Omega, T\Omega)\\
\int_\Omega  \boldsymbol{\phi} \,\cdot \,\mathbf{u}_j^{\begin{tiny}\mbox{(D)}\end{tiny}}(\lambda)  \; dV=0, \; \;j =1,\cdots, k-1}} \!\!\frac{ \int_{\Omega}  \Big[2\mu \big( \mbox{Def}\; \boldsymbol{\phi}, \mbox{Def}\; \boldsymbol{\phi}\big)+ \lambda\,\big(\mbox{div}\; \boldsymbol{\phi}\big)^2 \Big] \,dV }{\int_\Omega \boldsymbol{\phi}\cdot  \boldsymbol{\phi} \; dV},\nonumber \\
 [6mm] &&  \quad \quad  \cdots  \cdots \cdots \cdots,
 \nonumber \end{eqnarray} 
  and 
  \begin{eqnarray}&&\label{24.4.27-3}  \!\!\!\tau_1^{\begin{scriptsize}\mbox{(T)}\end{scriptsize}}(\lambda) \, = \frac{ \int_{\Omega}  \Big[2\mu \big( \mbox{Def}\, \mathbf{u}^{\begin{scriptsize}\mbox{(T)}\end{scriptsize}}_1(\lambda),\mbox{Def}\, \mathbf{u}_1^{\begin{scriptsize}\mbox{(T)}\end{scriptsize}} (\lambda)\big)+ \lambda\,\big(\mbox{div}\, \mathbf{u}_1^{\begin{scriptsize}\mbox{(T)}\end{scriptsize}}(\lambda)\big)^2 \Big] \,dV }{\int_\Omega \mathbf{u}_1^{\begin{scriptsize}\mbox{(T)}\end{scriptsize}}(\lambda) \cdot  \mathbf{u}_1^{\begin{scriptsize}\mbox{(T)}\end{scriptsize}}(\lambda)\; dV}\\
 [1mm] &&  \qquad \;\; = \min_{\substack{\boldsymbol{\phi} \in H^1(\Omega, T\Omega)}} \frac{ \int_{\Omega}  \Big[2\mu \big(\mbox{Def}\, \boldsymbol{\phi},\mbox{Def}\, \boldsymbol{\phi}\big)+ \lambda\,\big(\mbox{div}\, \boldsymbol{\phi}\big)^2 \Big] \,dV }{\int_\Omega \boldsymbol{\phi}\cdot  \boldsymbol{\phi} \; dV}, \qquad \nonumber \\
 [6mm] && \quad \;   \cdots  \cdots \cdots \cdots \nonumber  \\
 [4mm] \label{24.4.2-10} && \tau_{k}^{\begin{scriptsize}\mbox{(T)}\end{scriptsize}}(\lambda) =\frac{ \int_{\Omega}  \Big[2\mu \big(\mbox{Def}\, \mathbf{u}^{\begin{scriptsize}\mbox{(T)}\end{scriptsize}}_k(\lambda),\mbox{Def}\, \mathbf{u}_k^{\begin{scriptsize}\mbox{(T)}\end{scriptsize}} (\lambda)\big)+ \lambda\,\big(\mbox{div}\, \mathbf{u}_k^{\begin{scriptsize}\mbox{(T)}\end{scriptsize}}(\lambda)\big)^2 \Big] \,dV }{\int_\Omega \mathbf{u}_k^{\begin{scriptsize}\mbox{(T)}\end{scriptsize}}(\lambda)\cdot  \mathbf{u}_k^{\begin{scriptsize}\mbox{(T)}\end{scriptsize}}(\lambda) \; dV}\\
[1mm] && \qquad\;\,\;\, =  \min_{\substack{\boldsymbol{\phi} \in H^1(\Omega, T\Omega) \\
\int_\Omega    \boldsymbol{\phi}\,\cdot \, \mathbf{u}_j^{\begin{tiny}\mbox{(T)}\end{tiny}}(\lambda)  \; dV=0, \; \;j =1,\cdots, k-1}}\!\! \frac{ \int_{\Omega}  \Big[2\mu \big( \mbox{Def}\, \boldsymbol{\phi},\mbox{Def}\, \boldsymbol{\phi}\big)+ \lambda\,\big(\mbox{div}\, \boldsymbol{\phi}\big)^2 \Big] \,dV }{\int_\Omega \boldsymbol{\phi}\cdot  \boldsymbol{\phi}\; dV}, \nonumber \\
 [4mm]  && \quad    \cdots  \cdots \cdots \cdots. \nonumber 
  \end{eqnarray} 
Thus, we have to show that the solutions of the above variational problems are also eigenvectors for the Lam\'{e} boundary value problems, and they  furnish {\it all} the eigenvectors.     For $\boldsymbol{\phi} \in H^1_0(\Omega,T\Omega)$ or $\boldsymbol{\phi}\in H^1(\Omega, T\Omega)$,  we  denote the  quadratic functionals for the variational eigenvalue problems: 
 \begin{eqnarray*}   A[ \boldsymbol{\phi}] = \int_{\Omega}  \Big(2\mu \big(\mbox{Def}\, \boldsymbol{\phi},\mbox{Def}\, \boldsymbol{\phi}\big)+ \lambda\,\big(\mbox{div}\, \boldsymbol{\phi}\big)^2 \Big) \,dV\end{eqnarray*}
and \begin{eqnarray*} B[\boldsymbol{\phi}] =\int_\Omega  \boldsymbol{\phi}\cdot  \boldsymbol{\phi} \; dV.\end{eqnarray*}
 It is easy to verify that  the associated polar forms
 \begin{eqnarray*}
  &&   A[ \boldsymbol{\phi}, \boldsymbol{\psi}] = \int_{\Omega}  \Big(2\mu \big(\mbox{Def}\, \boldsymbol{\phi},\mbox{Def}\, \boldsymbol{\psi}\big)+ \lambda\,\big(\mbox{div}\, \boldsymbol{\phi}\big) \big(\mbox{div}\, \boldsymbol{\psi}\big) \Big) \,dV,\\
& &  B[\boldsymbol{\phi}, \boldsymbol{\psi}] =\int_{\Omega}  \boldsymbol{\phi}\cdot \boldsymbol{\psi}\;dV
 \end{eqnarray*}
 satisfy the relations: 
 \begin{eqnarray*} &&{A} [\boldsymbol{\phi} +\boldsymbol{\psi} ] ={ A} [\boldsymbol{\phi}]+ 2{A} [\boldsymbol{\phi}, \boldsymbol{\psi}] +{A} [\boldsymbol{\psi}],\\
 && B[\boldsymbol{\phi}+\boldsymbol{\psi}]= B[\boldsymbol{\phi}] +2B[\boldsymbol{\phi}, \boldsymbol{\psi}]+B[\boldsymbol{\psi}].\end{eqnarray*}

Without loss of generality, we only discuss the case of the  variational problem with traction boundary condition  in $H^1(\Omega, T\Omega)$ (since,  for the case  of the variational problem with Dirichlet boundary condition in $H_0^1(\Omega, T\Omega)$, it is quite similar and easier). We  assume that its solution $\mathbf{u}_1^{\begin{scriptsize}\mbox{(T)}\end{scriptsize}}(\lambda)$ satisfies $B[\mathbf{u}_1^{\begin{scriptsize}\mbox{(T)}\end{scriptsize}}(\lambda)]=1$. 
 If $\mathbf{w}$ is any vector field  in $H^1(\Omega, T\Omega)$, and if $\epsilon$ is an arbitrary constant, then for every value of $\epsilon$ and for $\mathbf{u}= \mathbf{u}_1^{\begin{scriptsize}\mbox{(T)}\end{scriptsize}}(\lambda)$, $\tau =\tau_1^{\begin{scriptsize}\mbox{(T)}\end{scriptsize}}(\lambda)$, we have 
 \begin{eqnarray*} {A} [ \mathbf{u} +\epsilon \mathbf{w}] \ge \tau B[\mathbf{u}+\epsilon \mathbf{w}],\end{eqnarray*}    
i.e.,  \begin{eqnarray}   \label{24.4.20-09} \epsilon\Big\{{A} [\mathbf{u}, \mathbf{w}] -\tau B[\mathbf{u}, \mathbf{w}] +\frac{\epsilon}{2} \big( {A}[ \mathbf{w}] -\tau B[\mathbf{w})\big) \Big\} \ge 0\end{eqnarray}  
in view of the relation ${A}[\mathbf{u}] =\tau B[\mathbf{u}]$. 
The inequality (\ref{24.4.20-09})
 can be valid for arbitrary values of $\epsilon$  only if equality 
 \begin{eqnarray} \label{24.4.20-02} {A}[\mathbf{u}, \mathbf{w}] -\tau B[\mathbf{u}, \mathbf{w}] =0\end{eqnarray}
 holds, i.e., if   the expression ${A} -\tau B$   vanishes  (corresponding to the first variational minimum problem (\ref{24.4.27-3})). 
 By applying Green's formula (see, Section 2 of \cite{Liu-19}),  we get 
  \begin{eqnarray} \label{24.5.16-1} {A} [\mathbf{u}, \mathbf{w}] = \int_\Omega \big( L_{\mu, \lambda} \mathbf{u}\big) \cdot  \mathbf{w}
  \; dV + \int_{\partial \Omega} \Big(2 \mu \,\big( \mbox{Def}\, \mathbf{u}\big)\boldsymbol{\nu}  +\lambda\, \big(\mbox{div}\, \mathbf{u}\big) \boldsymbol{\nu}\Big)\cdot  \mathbf{w}\; ds. \end{eqnarray}   
 Since the vector field $\mathbf{w}\in H^1(\Omega, T\Omega)$ is arbitrary, by (\ref{24.4.20-02}) -- (\ref{24.5.16-1})  we immediately obtain equations 
 \begin{eqnarray} \label{24.4.20-05} \left\{ \begin{array}{ll}  L_{\lambda, \mu} \mathbf{u}= \tau \mathbf{u} \;\;\; &\mbox{in}\;\; \Omega, \\
 2 \mu \, \big(\mbox{Def}\, \mathbf{u}\big) \,\boldsymbol{\nu} +\lambda \,\big( \mbox{div}\, \mathbf{u}\big) \,\boldsymbol{\nu} =0 \;\;\, & \mbox{on}\;\;\partial \Omega \end{array} \right.\end{eqnarray}  
  for   $\mathbf{u}=\mathbf{u}_1^{\begin{scriptsize}\mbox{(T)}\end{scriptsize}} (\lambda)$ and $\tau=\tau_1^{\begin{scriptsize}\mbox{(T)}\end{scriptsize}}(\lambda)$. 
   Considering the second minimum problem (in (\ref{24.4.2-10}))  with the additional condition $B[\boldsymbol{\phi}, \mathbf{u}_1^{\begin{scriptsize}\mbox{(T)}\end{scriptsize}}(\lambda)]=0$, where $\boldsymbol{\phi}\in H^1(\Omega, T\Omega)$, we obtain equation (\ref{24.4.20-02})  for $\mathbf{u}=\mathbf{u}_2^{\begin{scriptsize}\mbox{(T)}\end{scriptsize}}(\lambda)$ and $\tau= \tau_2^{\begin{scriptsize}\mbox{(T)}\end{scriptsize}}(\lambda)$ at first only under the assumption that $\mathbf{w}\in H^1(\Omega, T\Omega)$ satisfies relation
   \begin{eqnarray} \label{24.4.20-06} B[\mathbf{w}, \mathbf{u}_1^{\begin{scriptsize}\mbox{(T)}\end{scriptsize}}(\lambda)]=0.\end{eqnarray}   
  Now if $\boldsymbol{\eta}$  is any vector field in $H^1(\Omega, T\Omega)$, we can determine a number $t$ in such a way that the function $\mathbf{w}=\boldsymbol{\eta} +t \,\mathbf{u}_1$ satisfies condition (\ref{24.4.20-06}) (this can be realized by setting $t=- B[\mathbf{u}_1^{(T)}(\lambda), \boldsymbol{\eta}]$).      
 Particularly,  by substituting the vector field $\mathbf{w} =\mathbf{u}_2^{\begin{scriptsize}\mbox{(T)}\end{scriptsize}}(\lambda)$ in equation (\ref{24.4.20-02}) with $\mathbf{u} =\mathbf{u}_1^{\begin{scriptsize}\mbox{(T)}\end{scriptsize}}(\lambda)$ and $\tau= \tau_1^{\begin{scriptsize}\mbox{(T)}\end{scriptsize}}(\lambda)$,   we get  \begin{eqnarray} \label{24.4.20-08} {A}[\mathbf{u}_2^{\begin{scriptsize}\mbox{(T)}\end{scriptsize}}(\lambda). \mathbf{u}_1^{\begin{scriptsize}\mbox{(T)}\end{scriptsize}}(\lambda)]=0\end{eqnarray}
 since  $\mathbf{u}_2^{\begin{scriptsize}\mbox{(T)}\end{scriptsize}} (\lambda)$ satisfies the additional condition 
  \begin{eqnarray} \label{24.4.20-07}  B[\mathbf{u}_2^{\begin{scriptsize}\mbox{(T)}\end{scriptsize}}(\lambda), \mathbf{u}_1^{\begin{scriptsize}\mbox{(T)}\end{scriptsize}}(\lambda)]=0. \end{eqnarray} 
  If we substitute our vector field $\mathbf{w} =\boldsymbol{\eta} + t\, \mathbf{u}_1^{\begin{scriptsize}\mbox{(T)}\end{scriptsize}} (\lambda) $ in equation (\ref{24.4.20-02}), writing  $\mathbf{u}=\mathbf{u}_2^{\begin{scriptsize}\mbox{(T)}\end{scriptsize}}(\lambda)$, $\tau=\tau_2^{\begin{scriptsize}\mbox{(T)}\end{scriptsize}}(\lambda)$, we find 
  \begin{eqnarray} {A}[ \mathbf{u}, \boldsymbol{\eta}] -\tau B[\mathbf{u}, \boldsymbol{\eta}] +t\, \big( {A} [\mathbf{u}, \mathbf{u}_1^{\begin{scriptsize}\mbox{(T)}\end{scriptsize}}(\lambda)] -\tau B[\mathbf{u}, \mathbf{u}_1^{\begin{scriptsize}\mbox{(T)}\end{scriptsize}}(\lambda)]\big)=0\end{eqnarray} 
  or, taking into consideration equations (\ref{24.4.20-08} ) and (\ref{24.4.20-07}), 
  \begin{eqnarray} \label{24.4.4.20-09}  {A} [\mathbf{u}, \boldsymbol{\eta}]-\tau B[\mathbf{u}, \boldsymbol{\eta}]=0.\end{eqnarray}
  That is, equation (\ref{24.4.20-02}) holds also for arbitrary vector fields $\boldsymbol{\eta}$ or $\mathbf{w}$
   without regard to the auxiliary  condition (\ref{24.4.20-06}). From
 this it follows directly, as above, that equation (\ref{24.4.20-05})  is also valid for $\mathbf{u}= \mathbf{u}_2^{\begin{scriptsize}\mbox{(T)}\end{scriptsize}}(\lambda)$ and $\tau =\tau_2^{\begin{scriptsize}\mbox{(T)}\end{scriptsize}}(\lambda)$.     
 Continuing in this process, we conclusion that the eigenvalue equation (\ref{24.4.4.20-09}) holds generally for the solutions $\mathbf{u}_k^{\begin{scriptsize}\mbox{(T)}\end{scriptsize}}(\lambda)$ and the minimum values $\tau_k^{\begin{scriptsize}\mbox{(T)}\end{scriptsize}}(\lambda)$.    
 For the solutions of the problem, we have the relations 
 \begin{eqnarray} \label{24.4.20-10} \left. \begin{array} {ll} {A}[\mathbf{u}_k^{\begin{scriptsize}\mbox{(T)}\end{scriptsize}}(\lambda)]=\tau_k^{\begin{scriptsize}\mbox{(T)}\end{scriptsize}},\;\; \;\,\;\; {A}[\mathbf{u}_k^{\begin{scriptsize}\mbox{(T)}\end{scriptsize}}(\lambda), \mathbf{u}_j^{\begin{scriptsize}\mbox{(T)}\end{scriptsize}}(\lambda)]=0\\
  B[\mathbf{u}_k^{\begin{scriptsize}\mbox{(T)}\end{scriptsize}}(\lambda)]=1, \;\;\;\,\;\;\;\;\,  B[\mathbf{u}_k^{\begin{scriptsize}\mbox{(T)}\end{scriptsize}}(\lambda), \mathbf{u}_j^{\begin{scriptsize}\mbox{(T)}\end{scriptsize}}(\lambda)] =0 \end{array}\right. \;\;\; \;\;\; (k\ne j).\end{eqnarray} 
   Our variational problems produce an infinite sequence of eigenpairs of the associated differential equations problem with traction boundary condition. Conversely, it is also true.  That is,  solutions of the variational eigenvalues problem (in $H^1(\Omega, T\Omega)$)  are equivalent to the solutions of Lam\'{e} eigenvalue problem with traction boundary condition.  Furthermore, the regularity of equations (\ref{24.4.20-05}) implies that all weak solutions $\{ \mathbf{u}_k^{\mbox{\begin{scriptsize}(T)\end{scriptsize}}}(\lambda)\}_{k=1}^\infty$ belong to $C^\infty (\Omega, T\Omega)$.

  As pointed out before,  the variational eigenvalue problem (\ref{24.4.27-1})--(\ref{24.4.27-2}) (which is considered in $H_0^1(\Omega, T\Omega)$) produce an infinite sequence of eigenpairs $\{(\tau^{\begin{scriptsize}\mbox{(D)}\end{scriptsize}}_k(\lambda), \mathbf{u}_k^{\begin{scriptsize}\mbox{(D)}\end{scriptsize}}(\lambda))\}$ with Dirichlet boundary condition. In addition, it can easily  be shown that $\mathbf{u}_k^{\mbox{\begin{tiny}(D)\end{tiny}}}$ belong to $C^\infty (\Omega, T\Omega)$.
 Finally, if  $\lambda', \lambda''\in (-2\mu, +\infty)$ satisfy $\lambda'<\lambda''$, then  for any  $\boldsymbol{\phi}\in H_0^1(\Omega, T\Omega)$ (in the case of Dirichlet boundary condition) or  $ \boldsymbol{\phi}\in H^1(\Omega, T\Omega)$ (in the case of traction boundary condition), 
\begin{eqnarray}\label{24.4.20-3}  && \frac{ \int_{\Omega}  \Big[2\mu \big(\mbox{Def}\, \boldsymbol{\phi}, \mbox{Def}\, \boldsymbol{\phi}\big) + \lambda'\,\big(\mbox{div}\, \boldsymbol{\phi}\big)^2 \Big] \,dV }{\int_\Omega  \boldsymbol{\phi}\cdot \boldsymbol{\phi} \; dV} \le    \frac{ \int_{\Omega}  \Big[2\mu \big\langle\mbox{Def}\, \boldsymbol{\phi}, \mbox{Def}\; \boldsymbol{\phi}\big)+ \lambda''\,\big(\mbox{div}\, \boldsymbol{\phi}\big)^2 \Big] \,dV }{\int_\Omega \boldsymbol{\phi}\cdot  \boldsymbol{\phi} \; dV}. \qquad \qquad\end{eqnarray}
Combining (\ref{24.4.27-1}) -- (\ref{24.4.2-10}) and (\ref{24.4.20-3}), we get that $\tau^{(D)}_k(\lambda)\,$ (respectively, $\tau_k^{(T)}(\lambda)$) is continuous in $\lambda$, and for any $k\ge 1$,   
\begin{eqnarray} \left. \begin{array}{ll} & \tau_k^{\begin{scriptsize}\mbox{(D)}\end{scriptsize}} (\lambda')\le \tau_k^{\begin{scriptsize}\mbox{(D)}\end{scriptsize}} (\lambda''), \\
[1mm]  & \tau_k^{\begin{scriptsize}\mbox{(T)}\end{scriptsize}}(\lambda')\le \tau_k^{\begin{scriptsize}\mbox{(T)}\end{scriptsize}} (\lambda'')\end{array} \right. \;\; \;\;\mbox{when}\;\; -2\mu <\lambda'<\lambda''<+\infty.  \end{eqnarray}
\qed

  \vskip 0.78 true cm 

\section{Relationship of eigenvalues for the Lam\'{e} operator, Stokes operator and Laplace operator}

\vskip 0.38 true cm 

 For the Laplace operator  (i.e., the Lam\'{e} operator when $\lambda+\mu=0$) with Dirichlet and traction boundary conditions,  we need the following:   

\vskip 0.39 true cm

 \noindent{\bf Remark 3.1.} \     {\it  (i) \  Let $(\Omega,g)$ be a smooth compact Riemann manifold with smooth boundary. It is well-known (see,  (4.26) and (4.28) on p.$\,$305 of \cite{Ta2}) that 
 \begin{eqnarray} \label{24.6.7-6} \mbox{Def}^* \mathbf{u}=- \mbox{div}\; \mathbf{u},\end{eqnarray}
 where $\big(\mbox{div}\; \mathbf{u}\big)^j=u^{jk}_{\;\, \;\; ;k}$, and the following Weitzenbock formula holds on $(\Omega,g)$:
 \begin{eqnarray} \label{24.5.25-4}  2\,\mbox{div}\; \mbox{Def}\; \mathbf{u}=-\nabla^*\nabla \mathbf{u} +\mbox{grad}\; \mbox{div}\; \mathbf{u} + \mbox{Ric}\,(\mathbf{u}).\end{eqnarray} 
 Note that (\ref{24.5.25-4})  can also be obtained from (\ref{1-1}), (\ref{24.4.12-9}),  (\ref{24.6.7-6}) and $\mu=-\lambda =1$. 
 Moreover,  from  (\ref{24.5.25-4}) and Green's formula for the Lam\'{e} operator (see,  \S2 of \cite{Liu-19}),  we have  that \begin{eqnarray}  2 \int_{\Omega} \big(\mbox{Def}\; \mathbf{u}, \mbox{Def}\; \mathbf{u}\big)\, dV - \| \mbox{div}\; \mathbf{u}\|_{L^2}^2   = \big( \nabla^*\nabla \mathbf{u}, \mathbf{u}\big)_{L^2} -(\mbox{Ric}\, (\mathbf{u}), \mathbf{u})_{L^2} \end{eqnarray} 
 for any $\mathbf{u}\in H^1_0(\Omega, T\Omega)\,$ or $\,\mathbf{u}\in \big\{\mathbf{v} \in H^1(\Omega, T\Omega) \,\big| \, \mathcal{T}_{-1, 1}\mathbf{v}=0\big\},$
 where $\mathcal{T}_{\lambda, \mu}\mathbf{v}= -2\mu (\mbox{Def}\; \mathbf{v})^\# \boldsymbol{\nu} -\lambda \, (\mbox{div}\; \mathbf{v})\boldsymbol{\nu}$.  It follows that  for  all $\mathbf{u}\in H^1_0(\Omega, T\Omega)\,$ or $\,\mathbf{u}\in \big\{\mathbf{v}\in H^1(\Omega, T\Omega) \,\big|\, (T_{-\mu, \mu} \mathbf{v}\big)\big|_{\partial \Omega}=0\big\}$, 
  \begin{eqnarray}\label{24.6.4-1} \mu\! \int_{\Omega} \!2  (\mbox{Def}\, \mathbf{u}, \mbox{Def}\, \mathbf{u}) \, dV\! -\! \mu\! \int_{\Omega} \!(\mbox{div}\, \mathbf{u})(\mbox{div}\, \mathbf{u}) \,dV \! =\! -\mu \!\int_{\Omega} \! (\Delta_B \mathbf{u})\cdot \mathbf{u}\;dV-\mu \big(\mbox{Ric} (\mathbf{u}), \mathbf{u}\big)_{L^2},\;\end{eqnarray}  
 where $\Delta_B:=-\nabla^*\nabla$.
 
 (ii) \   If $\Omega\subset \mathbb{R}^n$ is a bounded domain with smooth boundary $\partial \Omega$,  then the corresponding result  (\ref{24.6.4-1}) can  directly be verified as follows. 
 
 Let us recall that $(\mbox{Def}\, \mathbf{u}, \mbox{Def}\, \mathbf{w})=\sum_{j,k=1}^n \big(\frac{1}{2} (\partial_{x_j} u_k+\partial_{x_{{}_k}} u_j) \big)  
\big(\frac{1}{2} (\partial_{x_j} w_k+\partial_{x_{{}_k}} w_j) \big)$. Similarly,  $\,\mbox{grad}\, \mathbf{u}\! :\!\mbox{grad}\, \mathbf{w}\,$ denots $\,\sum_{j,k=1}^n (\partial_{x_j} u_k)(\partial_{x_j} u_k)\,$ and we write $\,\mbox{grad}\, \mathbf{u}\!:\! 
(\mbox{grad}\, \mathbf{w})^T\,$ for $\,\sum_{j,k=1}^n  (\partial_{x_j} u_k)(\partial_{x_{{}_k}} w_j)$.  Then it is easy to verify that  for any vector fields $\mathbf{u}$ and $\mathbf{v}$ in $[H^1(\Omega)]^n$:
\begin{eqnarray} \label{24.5.19-1}  2\, (\mbox{Def}\, \mathbf{u}, \mbox{Def}\, \mathbf{w}) =\mbox{grad}\, \mathbf{u}: \mbox{grad}\, \mathbf{w} + \mbox{grad}\, \mathbf{u}: (\mbox{grad}\, \mathbf{w})^T.  \end{eqnarray} 
 Applying the integration by parts we have 
 \begin{eqnarray}   \label{24.5.19-2}
\!\!\! &&\int_{\Omega} (\mbox{grad}\, \mathbf{u}:\mbox{grad}\, \mathbf{w} )\, dx = \int_{\Omega}  \sum_{j,k=1}^n (\partial_{x_j} u_k )
 (\partial _{x_j} w_k) \\
 \!\!\!&& \;\;\;= \int_{\Omega} \bigg[ \sum_{k=1}^n \mbox{div}\,\Big( w_k\partial_{x_1} u_k, \cdots,    w_k \partial_{x_n} u_k\Big) \; dx - \int_{\Omega}  \sum_{k=1}^n   w_k  (\Delta u_k) \, dx \nonumber \\
 \!\!\! &&  \;\;\; = \! \int_{\partial \Omega}  \sum_{j,k=1}^n  \big( w_k  \nu_j \partial_{x_j} u_k \big)
ds - \int_{\Omega} (\Delta \mathbf{u}) \cdot \mathbf{w}  \, dx, \; \;\;\, \;\;\; \forall \,\mathbf{u}, \mathbf{w}\in [H^2(\Omega)]^n.\quad  \nonumber \end{eqnarray}  
 Similarly,  two integrations by parts we get 
 \begin{eqnarray}\!\! \!\!\! \!\!\!\!\!\!\!\! &&  \int_{\Omega}  \mbox{grad} \, \mathbf{u} : (\mbox{grad}\, \mathbf{w})^T dx = -\int_{\Omega} \big(\mbox{grad}\;\mbox{div}\, \mathbf{u} \big)\cdot \mathbf{w} \; dx +\int_{\partial \Omega} \sum_{j,k=1}^n\! w_j \nu_k \partial_{x_j} u_k \; ds \;\;\;\;\,  \\
 \!\!\!\!\! \!\!\! \!\! \!\!\!\! &&\, \quad = \! \int_{\Omega} (\mbox{div}\; \mathbf{u})(\mbox{div}\, \mathbf{w})\, dx 
\! -\! \int_{\partial \Omega} \!(\mbox{div}\, \mathbf{u}) (\boldsymbol{\nu} \cdot \mathbf{w}) \,ds  
 +\!\int_{\partial \Omega} \!
 \sum_{j,k=1}^n   w_k \nu_j\partial_{x_{{}_k}} u_j \; ds, \;\;\, \; \forall \mathbf{u}, \mathbf{w}\in [H^2(\Omega)]^n.\quad  \nonumber 
 \end{eqnarray}
If follows from (\ref{24.5.19-1})--(\ref{24.5.19-3}) that  for any $\mathbf{u}, \mathbf{w}\in [H^2(\Omega)]^n$, 
\begin{eqnarray} \label{24.5.19-3} \!\!\!&&  \int_{\Omega} 2\, (\mbox{Def}\, \mathbf{u}, \mbox{Def}\, \mathbf{w})\,dx  =
\int_\Omega \Big[ -(\Delta \mathbf{u})\cdot \mathbf{w}   + (\mbox{div}\, \mathbf{u}) (\mbox{div}\, \mathbf{w})\Big]dx \\
\!\!\!&& \qquad\;\;\; + \int_{\partial  \Omega} \bigg[ \sum_{j,k=1}^n 
 \Big(w_k\nu_j  \partial_{x_j} u_k +  w_k\nu_j \partial_{x_{{}_k}} u_j \Big) + (\mbox{div}\, \mathbf{u})(\boldsymbol{\nu}\cdot  \mathbf{w})\bigg] ds \nonumber  \\
\!\!\! &&  \qquad = \int_{\Omega} \Big(-(\Delta \mathbf{u}) \cdot \mathbf{w}\,  + (\mbox{div}\, \mathbf{u})(\mbox{div}\; \mathbf{w}) \Big) \,dx +\int_{\partial \Omega}  \Big( 2 \, (\mbox{Def}\; \mathbf{u} )\boldsymbol{\nu} + (\mbox{div}\, \mathbf{u})\boldsymbol{\nu}\Big) \cdot \mathbf{w} \; ds.   \nonumber 
 \end{eqnarray} 
Thus, in the case of $\lambda+\mu=0$,  we have 
\begin{eqnarray} \!\!\!\! &&\mu \int_{\Omega} 2 \, (\mbox{Def}\, \mathbf{u}, \mbox{Def}\, \mathbf{w}) \, dx - \mu \int_{\Omega} (\mbox{div}\, \mathbf{u})(\mbox{div}\, \mathbf{w}) \,dx   \\
 \!\!\! \!&& \quad \;  = -\mu \int_{\Omega}  (\Delta \mathbf{u})\cdot \mathbf{w}\;dx+ \mu 
 \int_{\partial \Omega} \Big( 2\, (\mbox{Def}\, \mathbf{u}) \boldsymbol{\nu} +(\mbox{div}\, \mathbf{u}) \boldsymbol{\nu} \Big) \cdot \mathbf{w}\; ds, \quad \forall \, \mathbf{u}, \mathbf{w}\in [H^2(\Omega)]^n. \quad \;\quad   \nonumber\end{eqnarray} 
 This implies that   for any  $\,\mathbf{u}, \mathbf{w}\in [H^1_0(\Omega)]^n\,$ or $\;\mathbf{u}, \mathbf{w}\in \!\big\{\mathbf{v} \in [H^1(\Omega)]^n \,\big|\, -2\mu (\mbox{Def}\, \mathbf{v})^\# \boldsymbol{\nu} + \mu (\mbox{div}\, \mathbf{v})\boldsymbol{\nu}=0\big\},$
 \begin{eqnarray} \!\!\!\! &&\mu \int_{\Omega} 2 \, (\mbox{Def}\, \mathbf{u}, \mbox{Def}\, \mathbf{w}) \, dx - \mu \int_{\Omega} (\mbox{div}\, \mathbf{u})(\mbox{div}\, \mathbf{w}) \,dx  = -\mu \int_{\Omega}  (\Delta \mathbf{u})\cdot \mathbf{w}\;dx.\end{eqnarray} 
  
 Therefore, for the Dirichlet eigenvalue problem of the Laplacian  (i.e., the Lam\'{e} operator when $\lambda+\mu=0$),  our definition  is the same as the classical one. But,  for the Laplacian   with traction boundary condition,  the eigenvalue problem has a bit difference with classical Neumann Laplacian (the Rayleigh quotient is different). Their relation can be seen in (\ref{24.6.4-2}) of the proof of Theorem 1.4.   In physics, the traction Laplacian eigenvalues have significant application, it describes the vibrational case of an elastic body when the velocities of the longitudinal  and transverse waves  are equal.  }

\vskip 0.52 true cm

The following Lemma plays a key role for the Lam\'{e} operator $L_{\lambda, \mu}$:

\vskip 0.15 true cm

\noindent {\bf  Lemma 3.2.}  \  {\it Let $(\Omega, g)$ be a smooth compact Riemannian manifold with smooth boundary $\partial \Omega$, and let the Lam\'{e} coefficients satisfy $\mu>0$ and $\lambda+2\mu>0$. 
Then the following G{\aa}rding's inequality holds:  
  there are constants $\tilde{c}>0$ and $\tilde{d}>0$ such that for all $\mathbf{u}\in  H^{(D)/(T)}_\partial$, 
  \begin{eqnarray} \label{24.12.01-04}  \int_\Omega \Big(2\mu\, \big( \mbox{Def}\; \mathbf{u}, \mbox{Def}\; \mathbf{u}\big)
+\lambda \,(\mbox{div}\, \mathbf{u})\,(\mbox{div}\, \mathbf{u}) \Big)\,dV\ge \tilde{c}\,\|\mathbf{u}\|_{H^1(\Omega, T\Omega)}^2 -\tilde{d}\, \|\mathbf{u}\|_{L^2(\Omega, T\Omega)}^2,\end{eqnarray}
 where $H^{(D)}_\partial:= H^1_0(\Omega, T\Omega)\cap H^2(\Omega, T\Omega)$ and  $H^{(T)}_{\partial}:=
  \big\{\mathbf{w}\in H^2(\Omega, T\Omega) \,\big| \, \mathcal{T}_{\lambda, \mu} \mathbf{w}=0\;\, \mbox{on}\;\, \partial \Omega\big\}$.}

\vskip 0.42 true cm

   \noindent  {\it Proof.}  \      It is easy to verify that the Lam\'{e} operator \  is  strongly elliptic,  \  and the Lam\'{e} operator $L_{\lambda,\mu}$ with Dirichlet (respectively, traction) boundary condition is symmetric (i.e.,  $(L_{\lambda, \mu}^{\begin{tiny} (D)/(T)\end{tiny}} \mathbf{u}, \mathbf{v})_{L^2(\Omega, T\Omega)} =(\mathbf{u}, L_{\lambda, \mu}^{(D)/(T)} 
\mathbf{v})_{L^2(\Omega, T\Omega)}$ for all $\mathbf{u}, \mathbf{v}\in H_\partial^{(D)/(T)}$, see also \cite{Liu-39}).  Combing this and  Proposition 1.5 on p.$\,$318 of \cite{MaHu},  we  immediately  obtain  the desired result (\ref{24.12.01-04}).  \qed

\vskip 0.52 true cm

   \noindent  {\it Proof of Theorem 1.2.}  \     In order to keep the proof as simple as possible, we assume that the Stokes eigenvalues are simple and that each Lam\'{e} eigenvalue   is simple for $\lambda$ large enough.  The proof is divided into two cases according to the boundary  conditions.

i) \ \  Let $(\varsigma^{[0]}, \mathbf{u}^{[0]})$ be a Stokes eigenpair with Dirichlet boundary condition (assume that it is simple). Then there exists a unique 
 $p^{[0]} \in L^2(\Omega)$ with  $\; \int_{\Omega} p^{[0]}\,dV=0$ such that 
\begin{eqnarray*} \left\{  \begin{array}{ll} \mu \big(\nabla^*\nabla \mathbf{u}^{[0]}-\mbox{Ric}\,(\mathbf{u}^{[0]})\big) +\nabla_g p^{[0]} =\varsigma^{[0]} \mathbf{u}^{[0]}\;\;  &\mbox{in}\;\; \Omega,\\
\mbox{div}\, \mathbf{u}^{[0]}=0 \;\; &\mbox{in}\;\; \Omega,\\
\mathbf{u}^{[0]}=0 \;\; & \mbox{on}\;\; \partial \Omega, \end{array} \right. 
\end{eqnarray*} 
since $p^{[0]}$ is unique up to an additive arbitrary constant.  Clearly, by multiplying a constant, it can be assume that $\|\mathbf{u}^{[0]}\|_{H_0^1(\Omega, T\Omega)}< \frac{1}{2}$.  We can construct a triplet $(\varsigma^{[1]}, \mathbf{u}^{[1]}, p^{[1]})$ with $\varsigma^{[1]}\in \mathbb{R}$, $\mathbf{u}^{[1]}\in J^{0}\cap H_0^1(\Omega, T\Omega)\cap H^2(\Omega, T\Omega)\,$, and $p^{[1]}\in L^2 (\Omega)$,  $\;\int_\Omega p^{[1]}\; dV=0$  such that 
\begin{eqnarray} \label{24.4.29-3} \left\{ \begin{array}{ll}  \mu\big(\nabla^*\nabla \mathbf{u}^{[1]} -\mbox{Ric}\,(\mathbf{u}^{[1]}) \big) +\nabla_g p^{[1]} = \varsigma^{[0]} \mathbf{u}^{[1]}+ \varsigma^{[1]} \mathbf{u}^{[0]} \;\; &\mbox{in}\;\; \Omega,\\
 \mbox{div}\, \mathbf{u}^{[1]}=-p^{[0]}\;\; &\mbox{in}\;\;\Omega\\
\mathbf{u}^{[1]}=0 \; \; & \mbox{on}\;\; \partial \Omega.\end{array} \right. \end{eqnarray} 
In fact, there exists a unique solution $\mathbf{\tilde{u}}$ satisfying the following  boundary value problem \begin{eqnarray} \label{24.4.29-1}
\left\{ \begin{array}{ll}  \mbox{div}\, \mathbf{\tilde{u}} =-p^{[0]} \;\; &\mbox{in}\;\;\Omega,\\
\mathbf{\tilde{u}}=0\;\; &\mbox{on}\;\; \partial \Omega.\end{array} \right.\end{eqnarray} 
It follows from (\ref{24.4.29-3})  and (\ref{24.4.29-1}) that 
\begin{eqnarray}\label{24.4.29-7} 
\! \left\{ \!\begin{array}{ll} \mu \Big( \nabla^*\nabla \big(\mathbf{u}^{[1]}\!-\! \mathbf{\tilde{u}}\big) -\mbox{Ric}\, \big(\mathbf{u}^{[1]} \!-\!\mathbf{\tilde{u}}\big) \Big)+ \nabla_g p^{[1]} =  \varsigma^{[0]} \big(\mathbf{u}^{[1]}\!-\! \mathbf{\tilde{u}}\big)+ \varsigma^{[1]}\mathbf{u}^{[0]}\\
\qquad\qquad \quad \quad \qquad \qquad  \qquad \qquad \quad \qquad -\mu\big( \nabla^*\nabla \mathbf{\tilde{u}} -\mbox{Ric}\, \mathbf{\tilde{u}} \big) + \varsigma^{[0]} \mathbf{\tilde{u}}\;\;\;&\mbox{in}\;\; \Omega,  \\ 
 \mbox{div}\, \big(\mathbf{u}^{[1]}-\mathbf{\tilde{u}}\big)=0\;\;\; &\mbox{in}\;\;\Omega\\
\mathbf{u}^{[1]}-\mathbf{\tilde{u}}=0  \;\;\; &\mbox{on}\;\;\partial \Omega.\end{array} \right. \quad \end{eqnarray} 
Put $\mathbf{\tilde{\tilde{u}}} :=\mathbf{u}^{[1]}-\mathbf{\tilde{u}}$ and $\mathbf{f}^{[1]}:= -\mu\big(\nabla^*\nabla \mathbf{\tilde{u}} -\mbox{Ric}\,(\mathbf{\tilde{u}})\big)+\varsigma^{[0]} \mathbf{\tilde{u}}$.  Then (\ref{24.4.29-7}) becomes 
\begin{eqnarray}\label{24.4.29-11} \left\{ \begin{array}{ll} \big(S_\mu -\varsigma^{[0]} \big) \mathbf{\tilde{\tilde{u}}} =\mathbf{f}^{[1]} +\varsigma^{[1]}\mathbf{u}^{[0]}\;\;\, & \mbox{in}\;\; \Omega,\\
\mbox{div}\; \mathbf{\tilde{\tilde{u}}}=0\;\;\,& \mbox{in}\;\; \Omega,\\
\mathbf{\tilde{\tilde{u}}}=0 \;\;& \mbox{on}\;\; \partial \Omega.\end{array} \right.  \end{eqnarray} 
  which is  
 solvable if its right hand side is orthogonal to $\mathbf{u}^{[0]}$:  this can be realized by the choice of $\varsigma^{[1]}$: 
\begin{eqnarray*} \varsigma^{[1]} =- \,\frac{\big( \mathbf{f}^{[1]}, \mathbf{u}^{[0]}\big)_{L^2}}{ \big(\mathbf{u}^{[0]}, \mathbf{u}^{[0]}\big)_{L^2}}. \end{eqnarray*} 
By solving (\ref{24.4.29-11}) (i.e., (\ref{24.4.29-7})),  we  get $\mathbf{\tilde{\tilde{u}}}$ and $p^{[1]}$, so that  we obtain the required $(\varsigma^{[1]}$, $\mathbf{u}^{[1]}, p^{[1]})$ (since $\tilde{\mathbf{u}}$ is known by (\ref{24.4.29-1})).   
 In addition, it can be assume $\|\mathbf{u}^{[1]}\|_{H_0^1(\Omega, T\Omega)}< \frac{1}{2^2}$. 
Generally,  using the same way we can construct a sequence $(\varsigma^{[j]}, \mathbf{u}^{[j]}, p^{[j]})$, $\, j\ge 1$, with $\varsigma^{[j]}\in \mathbb{R}$, $\,\mathbf{u}^{[j]}\in J^0\cap H^1_0(\Omega, T\Omega) \cap H^2(\Omega,T\Omega)$ 
and $p^{[j]}\in L^2 (\Omega)$, $\; \int_\Omega p^{[j]}\, dV =0$ such that 
\begin{eqnarray}  \label{24.4.30-01}\!\!\! \!\!\!\! \left\{ \! \begin{array}{ll} \mu\big(\nabla^*\nabla \mathbf{u}^{[j]} \!-\! \mbox{Ric}\, (\mathbf{u}^{[j]}) \big) \!+\!\nabla_g p^{[j]} \!=\! \varsigma^{[0]} \mathbf{u}^{[j]} +\varsigma^{[1]} \mathbf{u}^{[j-1]} +\cdots + \varsigma^{[j]} \mathbf{u}^{[0]}\;\; &\mbox{in}\;\,\Omega,\\
 \mbox{div}\, \mathbf{u}^{[j]} = -p^{[j-1]} \;\,&\mbox{in}\;\;\Omega,\\
 \mathbf{u}^{[j]}=0 \;\; &\mbox{on}\;\, \partial \Omega.\end{array} \right. \end{eqnarray}
Indeed, as the above discussion, after subtraction of a solution $\mathbf{\tilde{w}} \in  H^1_0 (\Omega, T\Omega)\cap H^2(\Omega, T\Omega)$ of the equation 
\begin{eqnarray} \left\{ \begin{array}{ll} \mbox{div}\, \mathbf{\tilde{w}}=-p^{[j-1]} \;\, &\mbox{in}\;\; \Omega,\\
\mathbf{\tilde{w}}=0 \;\, &\mbox{on}\;\; \partial \Omega,\end{array} \right. \end{eqnarray} 
(\ref{24.4.30-01})  has the form 
\begin{eqnarray} \label{24.4.29-13} \big(S_\mu -\varsigma^{[0]}) \,\mathbf{\tilde{\tilde{w}}} =\mathbf{f}^{[j]} +\varsigma^{[j]} \mathbf{u}^{[0]}\end{eqnarray}  by setting $\mathbf{\tilde{\tilde{w}}}:= \mathbf{u}^{[j]} - \mathbf{\tilde{w}}$ and 
$\mathbf{f}^{[j]} := -\mu \big( \nabla^* \nabla \mathbf{\tilde{w}} -\mbox{Ric}\,(\mathbf{\tilde{w}})\big) + \varsigma^{[1]} \mathbf{u}^{[j-1]} +\cdots + \varsigma^{[j-1]} \mathbf{u}^{[1]}$,   which is solvable if its right hand side is orthogonal to $\mathbf{u}^{[0]}$:  this can also  be achieved by the choice of $\varsigma^{[j]}$: 
\begin{eqnarray*} \varsigma^{[j]} = -\,\frac{\big( \mathbf{f}^{[j]}, \mathbf{u}^{[0]}\big)_{L^2}}{ \big(\mathbf{u}^{[0]}, \mathbf{u}^{[0]}\big)_{L^2}}. \end{eqnarray*} 
By solving (\ref{24.4.29-13})  we then get the desired result  (i.e., $(\varsigma^{[j]}, \mathbf{u}^{[j]}, p^{[j]})$ satisfying  (\ref{24.4.30-01})) and $\|\mathbf{u}^{[j]}\|_{H_0^1(\Omega, T\Omega)}<\frac{1}{2^{j+1}}$.

With $\epsilon=1/ (\lambda+\mu)$, we see that the power series in $\epsilon$, $\;\big(\!\sum_{j\ge 0} \epsilon^j \varsigma^{[j]}, \sum_{j\ge 0} \epsilon^j \mathbf{u}^{[j]}\big)$ is a formal eigenpair of $L_{\lambda, \mu}$ with Dirichlet boundary condition, since 
 \begin{eqnarray*} L_{\lambda,\mu} \mathbf{u} = \mu \big( \nabla^*\nabla \mathbf{u} -\mbox{Ric}\,(\mathbf{u})\big) - (\lambda+\mu) \,\mbox{grad}\, \mbox{div}\, \mathbf{u}.\end{eqnarray*}
Let us show now that this also holds in the sense of asymptotic expansions: setting for any
 $m\ge 1$,
\begin{eqnarray*} \underline{\varsigma}^{(m)} =\sum_{j=0}^m \epsilon ^j \varsigma^{[j]} \;\;\;\mbox{and}\;\;\, \underline{\mathbf{u}}^{(m)} =\sum_{j=0}^m 
\epsilon^j \mathbf{u}^{[j]}\end{eqnarray*}
 we find by $\lambda+\mu=1/\epsilon$,  $\,\mbox{div}\; \mathbf{u}^{[j]}=-p^{[j-1]}$ and (\ref{24.4.30-01}) that
  \begin{eqnarray*}\!\!\!\!\! \!\!\!\! \! & & \!\! L_{\lambda,\mu} \underline{\mathbf{u}}^{(m)} -\underline{\varsigma}^{(m)} \underline{\mathbf{u}}^{(m)} \\
\!\!\!\! \!\! \!\!\! & & \!\! = \sum_{j=0}^{m} \epsilon^j  \Big[ \mu \big( \nabla^*\nabla \mathbf{u}^{[j]} -\mbox{Ric}\, \big( \mathbf{u}^{[j]}\big) \big) \Big] 
   - (\lambda+\mu) \,\mbox{grad}\;  \Big(\sum_{j=0}^m  \epsilon^{j} \,\mbox{div}\; \mathbf{u}^{[j]}\Big) - \Big( \sum_{j=0}^m \epsilon^j \varsigma^{[j]} \Big)\Big( \sum_{j=0}^m \epsilon^j \mathbf{u}^{[j]}\Big)    \\  
   \!\!\!\!\! \! \!\!\! & &\!\!  = \epsilon^{m} \Big[ \mu\Big( \nabla^*\nabla \mathbf{u}^{[m]}\! -\!\mbox{Ric}\, \big(\mathbf{u}^{[m]}\big)  \Big) \Big] 
  \! + \!   \sum_{j=0}^{m-1}\! \epsilon^{j} \Big[ \mu\Big( \nabla^*\nabla \mathbf{u}^{[j]} \!-\!\mbox{Ric}\, \big(\mathbf{u}^{[j]}\big)\Big) \!+\!\mbox{grad}\; p^{[j]} \Big]   \\
 \!\!\!  &&\!\!   \;\; \; \, - \,\bigg[ \sum_{j=0}^{m-1} \epsilon^j\, \Big( \sum_{k=0}^j  \varsigma^{[k]} \mathbf{u}^{[j-k]}  \Big)  +
   \sum_{j=0}^m \epsilon^{m+j} \Big( \sum_{k=j}^m \varsigma^{[k]} \mathbf{u}^{[m+j-k]}  \bigg]   \\
   \!\!\!\!\! \!\! \! && = \epsilon^{m} \Big[ \mu\Big( \nabla^*\nabla \mathbf{u}^{[m]}\! -\!\mbox{Ric}\, \big(\mathbf{u}^{[m]}\big)  \Big) \Big] 
  \! + \!   \sum_{j=0}^{m-1}\! \epsilon^{j}  \sum_{k=0}^j \varsigma^{[k]} \mathbf{u}^{[j-k]}\\
   \!\!\!  && \!  \;\; \; \,  - \,\bigg[ \sum_{j=0}^{m-1} \epsilon^j\, \Big( \sum_{k=0}^j  \varsigma^{[k]} \mathbf{u}^{[j-k]}  \Big)  +
   \sum_{j=0}^m \epsilon^{m+j} \Big( \sum_{k=j}^m \varsigma^{[k]} \mathbf{u}^{[m+j-k]} \Big) \bigg]  
   \\   \!\!\! \!\!\!\! \!\!  &&\!\!  =\epsilon^m \Big[ \mu\Big( \nabla^*\nabla \mathbf{u}^{[m]} -\mbox{Ric} \,\big(\mathbf{u}^{[m]}\big)\Big)\Big] -  \sum_{j=0}^m \epsilon^{m+j} \Big( \sum_{k=j}^m \varsigma^{[k]} \mathbf{u}^{[m+j-k]}\Big).\end{eqnarray*}
  Whence the uniform estimate for $\epsilon$
   small enough (i.e., $\lambda$  large enough)
    \begin{eqnarray*}    \| L_{\lambda, \mu}\underline{\mathbf{u}}^{(m)} -\underline{\varsigma}^{(m)} \underline{\mathbf{u}}^{(m)} 
\|_{(H_0^1(\Omega,T\Omega))' } \le C \epsilon^m \| \underline{\mathbf{u}}^{(m)}\|_{H_0^1(\Omega,T\Omega)},\end{eqnarray*}
where $(H_0^1(\Omega,T\Omega))'$ is the dual space of $H_0^1(\Omega,T\Omega)$.
However, with $\mathbf{u}_k^{\begin{tiny}(D)\end{tiny}}(\lambda)$  the normalized (in $L^2(\Omega, T\Omega)$) 
 eigenvector associated to the eigenvalue $\tau_k^{\begin{tiny}(D)\end{tiny}}(\lambda)$ of the Lam\'{e} operator $L_{\lambda,\mu}$ with Dirichlet boundary condition,    for any $\varsigma \in \mathbb{R}$  and $\mathbf{u}\in H^1_0(\Omega, T\Omega)\cap H^2(\Omega, T\Omega)$  there holds 
  \begin{eqnarray} \label{24.4.16-02}    \| L_{\lambda,\mu} \mathbf{u} -\varsigma \,\mathbf{u} \|_{(H^1_0(\Omega,T\Omega))'}^2  
  =\!&&\!\!\!\! \!\sum_{k\ge 1}  ( \tau_k^{\begin{tiny}(D)\end{tiny}}(\lambda) -\varsigma)^2 (\mathbf{u}, \mathbf{u}_k^{\begin{tiny}(D)\end{tiny}}(\lambda))^2 \|\mathbf{u}_k^{\begin{tiny}(D)\end{tiny}}(\lambda)\|_{(H_0^1 (\Omega, T\Omega))'}^2.\end{eqnarray}
  It follows from G{\aa}rding's inequality (\ref{24.12.01-04})  that for any $\mathbf{u}\in H_\partial^{(D)}$, 
      \begin{eqnarray*} \int_\Omega \Big[ 2\mu\, \big(\mbox{Def}\; \mathbf{u}, \mbox{Def}\; \mathbf{u}\big) +\lambda\,(\mbox{div}\, \mathbf{u})(\mbox{dix}\, \mathbf{u}) +\tilde{d} \, |\mathbf{u}|^2\Big] \,dV \ge \tilde{c} \,
   \|\mathbf{u} \|^2_{H^1_0(\Omega, T\Omega)},\end{eqnarray*}
      where $\tilde{c}$ and $\tilde{d}$ are the positive constants in G{\aa}rding's inequality. That is,   
      \begin{eqnarray}\label{24.5.27-11}  \big((L_{\lambda,\mu} +\tilde{d}) \mathbf{u}, \mathbf{u}\big)_{L^2(\Omega, T\Omega)} \ge \tilde{c} \,
   \|\mathbf{u} \|^2_{H^1_0(\Omega, T\Omega)},\;\;\;  \forall\, \mathbf{u}\in H_\partial^{(D)}.\end{eqnarray}
Noting  that $(L_{\lambda,\mu} +\tilde{d}) \,\mathbf{u}_k^{\begin{tiny}(D)\end{tiny}}(\lambda) = (\tau^{\begin{tiny}(D)\end{tiny}}_k+\tilde{d}) \,\mathbf{u}_k^{\begin{tiny}(D)\end{tiny}}(\lambda) $,  by (\ref{24.5.27-11}) we have 
\begin{eqnarray} \label{24.6.4-6}  \big((\tau^{\begin{tiny}(D)\end{tiny}}_k+\tilde{d}) \,\mathbf{u}_k^{\begin{tiny}(D)\end{tiny}}(\lambda), \mathbf{u}_k^{(D)}(\lambda)\big)_{L^2(\Omega, T\Omega)} \ge \tilde{c} \|\mathbf{u}_k^{(D)}(\lambda)\|^2_{H_0^1(\Omega, T\Omega)}.\end{eqnarray}  
 Since   
   \begin{eqnarray*} \| \mathbf{u}_k^{(D)} (\lambda) \|_{(H_0^{1} (\Omega, T\Omega))'} =\sup\limits_{\mathbf{w}\in H_0^1(\Omega, T\Omega)} \frac{\big| \big(\mathbf{u}_k^{\begin{tiny}(D)\end{tiny}}(\lambda) , \mathbf{w}\big)_{L^2(\Omega, T\Omega)} \big|}{ \|\mathbf{w}\|_{H_0^1(\Omega, T\Omega)} } \ge   \frac{ \big|\big(\mathbf{u}_k^{\begin{tiny}(D)\end{tiny}}(\lambda) , 
   \mathbf{u}_k^{\begin{tiny}(D)\end{tiny}}(\lambda)  \big)_{L^2(\Omega, T\Omega)} \big|}{ \|\mathbf{u}_k^{\begin{tiny}(D)\end{tiny}}(\lambda) \|_{H_0^1(\Omega, T\Omega)} }, \end{eqnarray*} 
  we find from this and  (\ref{24.6.4-6}) that  
    \begin{eqnarray*} (\tau_k^{\begin{tiny}(D)\end{tiny}} \!+\!\tilde{d} )\, \| \mathbf{u}_k^{\begin{tiny}(D)\end{tiny}}(\lambda)  \|^2_{(H_0^{1} (\Omega, T\Omega))'}\ge  \Bigg( \frac{  \big((\tau_k^{\begin{tiny}(D)\end{tiny}}(\lambda) \! +\!\tilde{d} ) \mathbf{u}_k^{\begin{tiny}(D)\end{tiny}}(\lambda), 
   \mathbf{u}_k^{\begin{tiny}(D)\end{tiny}}(\lambda)  \big)_{L^2(\Omega, T\Omega)}} { \|\mathbf{u}_k^{\begin{tiny}(D)\end{tiny}}(\lambda) \|^2_{H_0^1(\Omega, T\Omega)} }
   \Bigg)    \,  \| \mathbf{u}_k^{\begin{tiny}(D)\end{tiny}}(\lambda) \|^2_{L^2(\Omega, T\Omega)}   
  \ge \tilde{c}, \end{eqnarray*} 
 so  \begin{eqnarray} \label{24.5.27-12}\|\mathbf{u}_k^{\begin{tiny}(D)\end{tiny}}(\lambda) \|^2_{(H_0^1(\Omega, T\Omega))'} \ge \frac{\tilde{c}}{\tau_k^{\begin{tiny}(D)\end{tiny}}(\lambda)  +\tilde{d}}\end{eqnarray} 
 (Note that  $\tau_k^{\begin{tiny}(D)\end{tiny}}(\lambda) +\tilde{d} \ge \tau_1^{\begin{tiny}(D)\end{tiny}}(\lambda)+\tilde{d}\ge \tau_1^{\begin{tiny}(D)\end{tiny}}(-\mu) +\tilde{d}>0$ for all $k\ge 1$).
Combining this and  (\ref{24.4.16-02}), we find that for  any $\varsigma \in \mathbb{R}$  and $\mathbf{u}\in H^1_0(\Omega, T\Omega)\cap H^2(\Omega, T\Omega)$,
      \begin{eqnarray} \label{24.4.16-102}   &&  \| L_{\lambda,\mu} \mathbf{u} -\varsigma \,\mathbf{u} \|_{(H^1_0(\Omega,T\Omega))'}^2    \ge  C\, \sum_{k\ge 1} \frac{(\tau_k^{\begin{tiny}(D)\end{tiny}}(\lambda) -\varsigma)^2}{\tau_k^{\begin{tiny}(D)\end{tiny}}(\lambda)+\tilde{d} } (\mathbf{u}, \mathbf{u}_k^{\begin{tiny}(D)\end{tiny}}(\lambda))^2. \nonumber
  \end{eqnarray} 
We see that 
\begin{eqnarray}   \label{24.4.16-03}\sum_{k\ge 1} \frac{(\tau_k^{\begin{tiny}(D)\end{tiny}}(\lambda) -{\underline{\varsigma}}^{(m)})^2}{\tau_k^{\begin{tiny}(D)\end{tiny}}(\lambda)+\tilde{d}} (\underline{\mathbf{u}}^{(m)}, \mathbf{u}_k^{\begin{tiny}(D)\end{tiny}}(\lambda))^2\le C \epsilon^{2m}.\end{eqnarray}
Hence there exists $k_0$ such that $\lim_{\lambda \to \infty} \tau_{k_0}^{\begin{tiny}(D)\end{tiny}}(\lambda)$  is equal to $\lim_{\epsilon\to 0} \underline{\varsigma}^{(m)}=\varsigma^{[0]}$, 
 and 
 \begin{eqnarray}   \label{24.4.16-04} \exists \delta>0, \, \;\;\,  \forall k\ne k_0, \;\;\, \frac{(\tau_k^{\begin{tiny}(D)\end{tiny}}(\lambda) -\underline{\varsigma}^{(m)} )^2}{\tau_k^{\begin{tiny}(D)\end{tiny}}(\lambda)+\tilde{d}}\ge \delta,\end{eqnarray}  
thus
\begin{eqnarray*}  \sum_{k\ne k_0} (\underline{\mathbf{u}}^{(m)}, \mathbf{u}_k^{\begin{tiny}(D)\end{tiny}}(\lambda))^2\le C \epsilon^{2m}.\end{eqnarray*} 
On the other hand, for $ k = k_0$ ,  we must have $|\underline{\varsigma}^{(m)} -\tau_k^{\begin{tiny}(D)\end{tiny}}(\lambda)| \le C \epsilon^m$. 
 This proves that $\sum_{l\ge 0} \epsilon^l \varsigma^{[l]} $  (respectively,   $\sum_{l\ge 0} \epsilon^l \mathbf{u}^{[l]}$) 
is the asymptotic development of $\tau_{k_0}^{\mbox{\begin{tiny}(D)\end{tiny}}}(\lambda)$  (respectively,  $\mathbf{u}_{k_0}^{\mbox{\begin{tiny}(D) \end{tiny}}}(\lambda)$) as $\lambda\to +\infty$.

  Conversely, let us fix  $\tau^{\begin{tiny}(D)\end{tiny}}(\lambda) =\tau_k^{\begin{tiny}(D)\end{tiny}}(\lambda)$ the $k$-th eigenvalue of $L_{\lambda,\mu}$ with Dirichlet boundary condition. From Theorem 1.1 and the previous 
proof we deduce that  $\tau^{\begin{tiny}(D)\end{tiny}}(\lambda) \le \varsigma_k^{\begin{tiny}(D)\end{tiny}}$, where $\varsigma_k^{\begin{tiny}(D)\end{tiny}}$ is the $k$-th Stokes eigenvalue with Dirichlet boundary condition. Thus $\tau_k^{\begin{tiny}(D)\end{tiny}}(\lambda)$  is bounded as $\lambda\to \infty$ and has a limit 
$\varsigma^{[0]}$. 
The corresponding normalized eigenvectors $\mathbf{u} (\lambda) = \mathbf{u}_k^{\begin{tiny}(D)\end{tiny}}(\lambda)$ are thus bounded in the domain
of any power of $L_{\lambda,\mu}$, thus in  $H^{1+\delta}(\Omega, T\Omega)$ for $\delta>0$ small enough. Therefore $\mathbf{u}(\lambda)$ has a limit $\mathbf{u}^{[0]}$ in $J^0\cap H_0^1(\Omega, T\Omega)\cap H^2(\Omega,T\Omega)$. Going back to the equations satisfied by $\mathbf{u}(\lambda)$ we find that
as $\lambda\to +\infty$,  \begin{eqnarray*}  \mbox{div}\, \mathbf{u}^{[0]}=0 \;\; \mbox{and} \;\;  (\lambda+\mu) \,\mbox{grad}\, \mbox{div}\, \mathbf{u}(\lambda) \to \mu \big(\nabla^*\nabla \mathbf{u^{[0]}}\! -\!\mbox{Ric}(\mathbf{u}^{[0]})\big) -\varsigma^{[0]} \mathbf{u}^{[0]}\;\;\, \mbox{in}\;\; J^0\cap H_0^1\cap  H^2.\end{eqnarray*} 
Setting $p(\lambda) = (\lambda+\mu)\,\mbox{div}\, \mathbf{u}(\lambda)$, we obtain that it converges in $L^2 (\Omega)$ to a limit  $-p^{[0]}$ as $\lambda\to +\infty$.  Thus $(\varsigma^{[0]}, \mathbf{u}^{[0]})$ is a Stokes eigenpair. 
 
 ii) \ \ For the Stokes operator with Cauchy force boundary condition,   the proof  is almost the same as in i)  only 
 some little  modifications are done:  
 from (\ref{24.4.16-02}) on,  we have that 
    for any $\varsigma \in \mathbb{R}$  and $\mathbf{u}\in H^1(\Omega, T\Omega)\cap H^2(\Omega, T\Omega)$ 
    \begin{eqnarray} \label{24.4.106-002}    \| L_{\lambda,\mu} \mathbf{u} -\varsigma \,\mathbf{u} \|_{(H^1(\Omega,T\Omega))'}^2 \ge  C \sum_{k\ge 1}(\tau_k^{\begin{tiny}(T)\end{tiny}}(\lambda) -\varsigma)^2 (\mathbf{u}, \mathbf{u}_k^{\begin{tiny}(T)\end{tiny}}(\lambda))^2 \|\mathbf{u}_k^{\begin{tiny}(T)\end{tiny}}(\lambda)\|_{(H^1 (\Omega, T\Omega))'}^2, 
  \end{eqnarray} 
  where $\{(\tau_k^{\begin{tiny}(T)\end{tiny}}(\lambda), \mathbf{u}_k^{\begin{tiny}(T)\end{tiny}}(\lambda))\}_{k=1}^\infty$ are eigenpairs of $L_{\lambda, \mu}$ with traction boundary condition (note that $\tau_1^{\begin{tiny}(T)\end{tiny}}(\lambda)=0$ and $\mathbf{u}_1^{\begin{tiny}(T)\end{tiny}}(\lambda)=\mbox{constant}$).  
 Note that $\tau_k^{\begin{tiny}(T)\end{tiny}}(\lambda) +\tilde{d} \ge \tau_1^{\begin{tiny}(T)\end{tiny}}(\lambda)+\tilde{d}\ge \tau_1^{\begin{tiny}(T)\end{tiny}}(-\mu)+\tilde{d}\ge \tilde{d} >0$ for all $k\ge 1$.
 It is similar to the discussion of (i)  by applying G{\aa}rding's inequality that 
    \begin{eqnarray} \label{24.5.27-111}\|\mathbf{u}_k^{\begin{tiny}(T)\end{tiny}}(\lambda) \|^2_{(H_0^1(\Omega, T\Omega))'} \ge \frac{\tilde{c}}{\tau_k^{\begin{tiny}(T)\end{tiny}}(\lambda)  +\tilde{d}}.\end{eqnarray} 
  From this and (\ref{24.4.106-002}) we get that for  any $\varsigma \in \mathbb{R}$  and $\mathbf{u}\in H^{(T)}_\partial$,
      \begin{eqnarray} \label{24.4.16-102}   &&  \| L_{\lambda,\mu} \mathbf{u} -\varsigma \,\mathbf{u} \|_{(H^1_0(\Omega,T\Omega))'}^2    \ge  C\, \sum_{k\ge 1} \frac{(\tau_k^{\begin{tiny}(T)\end{tiny}}(\lambda) -\varsigma)^2}{\tau_k^{\begin{tiny}(T)\end{tiny}}(\lambda)+\tilde{d} } (\mathbf{u}, \mathbf{u}_k^{\begin{tiny}(T)\end{tiny}}(\lambda))^2. \nonumber
  \end{eqnarray} 
  Therefore 
\begin{eqnarray}   \label{24.4.16-03}\sum_{k\ge 1} \frac{(\tau_k^{\begin{tiny}(T)\end{tiny}} (\lambda)-{\underline{\varsigma}}^{(m)})^2}{\tau_k^{\begin{tiny}(T)\end{tiny}}(\lambda)+\tilde{d}} (\underline{\mathbf{u}}^{(m)}, \mathbf{u}_k^{\begin{tiny}(T)\end{tiny}}(\lambda))^2\le C \epsilon^{2m}.\end{eqnarray}
Hence there exists $k_0$ such that $\lim_{\lambda \to \infty} \tau_{k_0}^{\begin{tiny}(T)\end{tiny}}(\lambda)$  is equal to $\lim_{\epsilon\to 0} \underline{\varsigma}^{(m)}=\varsigma^{[0]}$, 
 and 
 \begin{eqnarray}   \label{24.4.16-04} \exists \delta>0, \, \;\;\,  \forall k\ne k_0, \;\;\, \frac{(\tau_k^{\begin{tiny}(T)\end{tiny}}(\lambda) -\underline{\varsigma}^{(m)} )^2}{\tau_k^{\begin{tiny}(T)\end{tiny}}(\lambda)+\tilde{d}}\ge \delta,\end{eqnarray}  
thus
\begin{eqnarray*}  \sum_{k\ne k_0} (\underline{\mathbf{u}}^{(m)}, \mathbf{u}_k^{\begin{tiny}(T)\end{tiny}}(\lambda))^2\le C \epsilon^{2m},\end{eqnarray*} 
and for $ k = k_0$  we must have $|\underline{\varsigma}^{(m)} -\tau_k^{\begin{tiny}(T)\end{tiny}}(\lambda)| \le C \epsilon^m$.
Others are completely similar to proof of (i).   Therefore, the desired result still holds.
  \vskip 0.28 true cm 
  Finally, since each of the Lam\'{e} eigenvalues is continuous with respect to the Lam\'{e} coefficients,  by noting Remark 3.1 we immediately obtain the corresponding conclusions   for the Laplace-type  operator $-\mu \big( \nabla^*\nabla -\mbox{Ric}\big)$ with giving boundary conditions as $\lambda\to -\mu$. 
  \qed

\vskip 0.8 true cm

\section{Generalized Ahlfors Laplacian and the third coefficient of the heat trace asymptotic of the Lam\'{e} operator}
\vskip 0.36 true cm

In order to exactly give the third coefficient of  the heat trace asymptotic expansion for the Lam\'{e} operator, we need to use a result of  T. Branson, P.  Gilkey, B. {\O}rsted,  A. Pierzchalski  in  \cite{BGOP} for the generalized Ahlfors Laplacian. Thus, we must show that the $1$-form representation of the Lam\'{e} operator $L_{\lambda,\mu}$ (defined on a Riemannian manifold $(\Omega, g)$) is just the generalized Ahlfors Laplacian.  Let us remark that the classical Ahlfors Laplacian originated from the conformal geometry and was introduced by Ahlfors in 1974 and 1976 (see \cite{Ahl-74} and \cite{Ahl-76}).     

 Let $T\Omega$ and $T^*\Omega$ be the tangent and cotangent bundles of the Riemannian $n$-manifold $(\Omega,g)$, respectively. The space of all $C^\infty$ vector fields will be denoted by $\mathscr{X}$. 
Recall that in terms of vector fields, the Lam\'{e} operator $ L_{\lambda, \mu}$ defined on $(\Omega, g)$  can be written as (cf. (\ref{1-1})): \begin{eqnarray} \label{23.9.22-2} L_{\lambda,\mu} \mathbf{u}:\!\!\!\!\!\!&& \!\!\!\!= \mu \big(\nabla^* \nabla \mathbf{u}\big)  -(\mu +\lambda) \,\mbox{grad}\; \mbox{div}\, \mathbf{u} -\mu \, \mbox{Ric}(\mathbf{u})\\
\!\!\!\!\!\!&&\!\!\!\!= \Big(-\mu\nabla_k\nabla^k u^j -(\lambda+\mu) \nabla^j\nabla_k u^k - \mu\,\mbox{Ric}^{j}_{\;k}\, u^k\Big) \frac{\partial }{\partial x_j}\nonumber \\
\!\!\!\!\!\!&&\!\!\!\!=\Big( - \mu \nabla_k \nabla^k  u^j -\lambda \nabla^j \nabla_ku^k -\mu \nabla_k \nabla^j u^k\Big)\,\frac{\partial}{\partial x_j}, \;\;\;\;\; \mathbf{u}= u^j \frac{\partial}{\partial x_j} \in  \mathscr{X},\nonumber\end{eqnarray}
 where $\nabla^*\nabla \mathbf{u}$ is given by (\ref{24.4.12-9}), and $\big\{ \frac{\partial }{\partial x_j} \big\}_{j=1}^n $ are the coordinate basis.  
 On the other hand,  the  Lam\'{e} operator as well as the associated boundary value problems can equivalently be discussed in terms of $1$-form.  
  If   $\alpha$ is  the $1$-form dual to the vector field $\mathbf{u}= u^j \frac{\partial}{\partial x_j}$ in the sense that 
 \begin{eqnarray*} \alpha(\mathbf{X})=g(\mathbf{u}, \mathbf{X}), \;\; \mathbf{X}\in \mathscr{X}, \end{eqnarray*}
 then  \begin{eqnarray} \label{23.9.22-5} \alpha=  u_j dx_j \;\;\, \mbox{and}\;\; u_j= g_{jl}u^l. \end{eqnarray}
   Noting  that  $\nabla g=0$,  one has,  in  index notation, 
   \begin{eqnarray*}&& \nabla_k \nabla^k u^j=  \nabla_k\nabla^k g^{jl} u_l 
  =
    g^{jl} g^{km} \nabla_k \nabla_m u_l
   = g^{jl}\nabla^m \nabla_m u_l ,\\
 &&  \nabla^j \nabla_ku^k=  \nabla^j \nabla_k  g^{kl} u_l= \nabla^j \nabla^l u_l =  g^{jm}\nabla_m \nabla^l u_l,\\
 &&  \nabla_k \nabla^j u^k   = \nabla_k \nabla^j g^{kl} u_l =  g^{kl} \nabla_k g^{jm} \nabla_m u_l= g^{jm} \nabla^l \nabla_m u_l.\end{eqnarray*} 
Combining the above facts and the last line of (\ref{23.9.22-2}) we find   that for $\mathbf{u}\!=\!u^j \frac{\partial}{\partial x_j} \in  \mathscr{X}$, 
\begin{eqnarray} \label{23.9.21-1} && \; L_{\lambda,\mu} \mathbf{u} = \Big(\!- \mu  g^{jl}\nabla^m \nabla_m u_l  -\lambda g^{jm}\nabla_m \nabla^l u_l-\mu g^{jm} \nabla^l \nabla_m u_l\Big) \frac{\partial }{\partial x_j}\\
&& \quad \quad \; =g^{jm}\Big(\!- \mu  \nabla^l \nabla_l u_m  -\lambda \nabla_m \nabla^l u_l-\mu \nabla^l \nabla_m u_l\Big) \frac{\partial }{\partial x_j}\nonumber\\
&& \quad \quad \;:= \phi^j \frac{\partial}{\partial x_j}, \nonumber\end{eqnarray}
where $$\phi^j:= g^{jm}\Big(\!- \mu  \nabla^l \nabla_l u_m  -\lambda \nabla_m \nabla^l u_l-\mu \nabla^l \nabla_m u_l\Big).$$
 Therefore,  we have  that 
  \begin{eqnarray} \label{23.9.22-3}  && L_{\lambda,\mu}^{\flat} (\alpha) : =\big( L_{\lambda, \mu}\mathbf{u}\big)^\flat= \psi_k \,dx_k, \end{eqnarray}
  where 
  \begin{eqnarray} &&\psi_k = g_{kj}\phi^j =     g_{kj} g^{jm}\Big(\!- \mu  \nabla^l \nabla_l u_m  -\lambda \nabla_m \nabla^l u_l-\mu \nabla^l \nabla_m u_l\Big) \\
 && \quad\;  = - \mu  \nabla^l \nabla_l u_k  -\lambda \nabla_k \nabla^l u_l-\mu \nabla^l \nabla_k u_l \nonumber\\
  &&\quad\;  =  -\mu \nabla^l \nabla_l u_k -(\lambda +2\mu) \nabla_k \nabla^l u_l +\mu \nabla^l \nabla_k u_l -2\mu \,\mbox{Ric}^l_{\, k}u_l,\nonumber
   \end{eqnarray}
    $\alpha$ is given by (\ref{23.9.22-5}),  and  $\flat$  is the flat operator (for a vector field) by lowering an index. 
 Hence \begin{eqnarray} \label{23.9.23-10} L_{\lambda,\mu}^\flat (\alpha) = \Big( -(\lambda +2\mu) \nabla_k \nabla^l u_l +\mu \nabla^l \nabla_k u_l -\mu \nabla^l \nabla_l u_k -2\mu \,\mbox{Ric}^l_{\, k}u_l\Big)\, dx_k.\end{eqnarray} 
   Let $d: \Lambda  T^*\Omega\to \Lambda T^*\Omega$ be the exterior differential operator, where $\Lambda T^*\Omega=\bigoplus_{p=0}^n \Lambda^p T^*\Omega$.    The adjoint operator $\delta$ of  $d$ acting on a $p$-form $\alpha$ is defined in terms of $d$ and the Hodge star operator by formula 
 \begin{eqnarray*} \delta \alpha =(-1)^{np+n+1} * d* \alpha,\end{eqnarray*}  and the Hodge star operator $*: \Lambda^p T^*\Omega \to \Lambda^{n-p}T^* \Omega,\;  p=0, \cdots, n$, is defined by \begin{eqnarray*}  \langle \gamma, \eta\rangle =\gamma \wedge  *\eta \end{eqnarray*} 
 for any $\gamma,\eta\in \Lambda^pT^*\Omega$.
  It is well known (see, for example, p.$\,$16--17 of \cite{CLN})  that $\nabla^l u_l= \mbox{div}\, \alpha= -\delta \alpha$ 
   for the $1$-form $\alpha= u_j dx_j$.
   Furthermore,  from Exercise 5 on p.$\,$561 in \cite{Ta3} we know that $\delta d \alpha =  \big(\nabla^l \nabla_k u_l-  \nabla^l\nabla_l u_k\big) dx_k$ for the above $1$-form $\alpha$.  Combining these facts and   (\ref{23.9.23-10}) we  obtain that
 \begin{eqnarray} \label{23.9.21-2} L_{\lambda,\mu}^{\flat} (\alpha) = (\lambda+2\mu) d\delta \alpha + \mu \delta d \alpha - 2\mu \, \mbox{Ric}\, (\alpha).\end{eqnarray}
 Here $\mbox{Ric\,}(\alpha)$ denotes the Ricci action on $1$-forms $\alpha$:
\begin{eqnarray*} \mbox{Ric} \,(\alpha) =\mbox{Ric}(\mathbf{u}, \cdot)= \big( \mbox{Ric}^l_{\,k} u_l \big)\, dx_k,\end{eqnarray*} 
where $\mathbf{u}$ is the vector field dual to $\alpha$.
(\ref{23.9.21-2}) is the $1$-form representation of the Lam\'{e} operator.

  The {\it generalized Ahlfors Laplacian}  (see,  \cite{BGOP})  is defined as the operator $P=ad\delta +b\delta d -\epsilon \rho$ on $\Lambda^1 T^*\Omega$, where $a$ and $b$ are positive constants and where $\epsilon \rho$ is an arbitrary constant multiple of the Ricci tensor.  
 Clearly, from (\ref{23.9.21-2}) we see that  the $1$-form representation $L_{\lambda,\mu}^\flat$ of the elastic Lam\'{e} operator $L_{\lambda,\mu}$  is just the {\it generalized Ahlfors Laplacian} on $\Lambda^1T^* \Omega$ with $\lambda+2\mu=a>0$, $\mu=b>0$ and $2\mu=\epsilon$.

 Combining the above discussion and Branson-Gilkey-{\O}rsted-Pierzchalski's theorem (see, \cite{BGOP}), we get the following:
  
  \vskip 0.25 true cm 
  
   \noindent{\bf  Lemma 4.1.} \ {\it Let   $L_{\lambda,\mu}$ be the Lam\'{e} operator defined on an  $n$-dimensional compact Riemannian manifold $(\Omega, g)$ with smooth boundary $\partial \Omega$. Then, for the Dirichlet (or traction)  boundary conditions,  there is an asymptotic expansion of the form:
\begin{eqnarray} \label{4.4.2-8}  \mbox{Tr}\; e^{-t L_{\lambda, \mu}} \sim \sum_{k=0}^\infty a_k(L_{\lambda,\mu}) \,t^{(k-n)/2} \;\;\, \mbox{as} \;\; t\to 0^{+},\end{eqnarray} 
where  the first coefficient $a_0 (L_{\lambda, \mu}) $ and the third coefficient $a_2(L_{\lambda,\mu})$ are  given by 
\begin{eqnarray}   \label{24.4.23-4} && a_0 (L_{\lambda,\mu}) = (4\pi)^{-n/2} \big( (\lambda+2\mu)^{-n/2} +(n-1) \mu^{-n/2}\big) \mbox{Vol}\, (\Omega),\end{eqnarray} 
\begin{eqnarray} \label{24.5.22-1}&a_2 (L_{\lambda,\mu}) = 
\frac{1}{6 (4\pi)^{n/2}} \Big[\Big( \frac{1}{ (\lambda+2\mu)^{(n-2)/2}} +\frac{n-7}{ \mu^{(n-2)/2}} + \frac{12\mu}{n} \big( \frac{1}{(\lambda+2\mu)^{n/2}} +\frac{n-1}{\mu^{n/2}} \big) \Big)\! \int_{\Omega} \! R(x)\, dV\\ 
 &\quad  + 2 \,\Big(\frac{1}{ (\lambda+2\mu)^{(n-2)/2}} +\frac{n-7}{ \mu^{(n-2)/2}}\Big) \int_{\partial \Omega}H(x) \,ds\Big]. \nonumber \end{eqnarray} 
 Here $R(x)$ is the scalar curvature at $x\in \Omega$,  $H(x):=L_{\alpha\alpha}(x)$ be the mean curvature of $\Omega$ at $x\in \partial \Omega$, and $L_{\alpha\beta}$ be the second fundamental form. }

\vskip 0.39 true cm 

 \noindent  {\it Proof.} \  From the above discussion, we see that the $1$-form representation of  the Lam\'{e} operator $L_{\lambda,\mu}$ is exactly the generalized Ahlfors Laplacian $a\,d \delta + b \,\delta d -\epsilon \rho$  with $a=\lambda+2\mu$, $b=\mu$ and $\epsilon =2\mu$.  It is completely similar to the proof of Theorem 4.2 (or  proof of Theorem 4.3)  of \cite{BGOP}) (i.e., by  calculating the coefficient $a_2$ on $S^{n-1}\times [0,\pi]$ and by applying  
 the invariance theory of Gilkey \cite{Gil2})  that 
   the first coefficient $a_0$ and the third coefficient $a_2$
 are obtained.  Actually, the desired conclusion  is the same as the result in (c) of Theorem 4.2 (or (c) of Theorem 4.3) on p.$\,$5 of \cite{BGOP})  with $a:=\lambda+2\mu$, $\,b:=\mu$ and $\epsilon=2\mu$. 
\qed

\section{Asymptotics formulas of heat traces for the Lam\'{e} operators}

\vskip 0.32 true cm
 
 \noindent  {\it Proof of Theorem 1.4.} \  Part proof (particularly, for two-term asymptotic expansions) is similar to \cite{Liu2} (see also \cite{Liu-23}).  For the  sake of completeness, we shall give a detailed proof up to three-term asymptotics.  From the theory of elliptic operators (see \cite{GiTr},  \cite{Mo3}, \cite{Pa},  \cite{Stew}), we see that the Lam\'{e} operator $-L_{\lambda, \mu}$ can generate  strongly continuous semigroup $(e^{-t L_{\lambda,\mu}^{\begin{tiny}(D)\end{tiny}}})_{t\ge 0}$  (respectively,  $(e^{-t L_{\lambda,\mu}^{\begin{tiny}(T)\end{tiny}}})_{t\ge 0}$)  with respect to the  Dirichlet (respectively,  traction boundary) condition,  in suitable spaces of vector-valued functions (for example, in $C_0(\Omega, T\Omega)$ (by Stewart \cite{Ste}) or in $L^2(\Omega, T\Omega)$ (by Browder \cite{Brow})), or in $L^p(\Omega, T\Omega)$ (by Friedman \cite{Fri}). Furthermore,
   there exist  matrix-valued functions ${\mathbf{K}}^{\begin{tiny}(D)\end{tiny}} (t, x, y)$ (respectively,  ${\mathbf{K}}^{\begin{tiny}(T)\end{tiny}} (t, x, y)$), which is called the integral kernel, such that (see \cite{Brow} or p.$\,$4 of \cite{Fri})
        \begin{eqnarray*}  e^{-t L_{\lambda,\mu}^{\begin{tiny}(D)/(T)\end{tiny}}}{\mathbf{w}}_0(x)=\int_\Omega {\mathbf{K}}^{\begin{tiny}(D)/(T)\end{tiny}}(t, x,y) {\mathbf{w}}_0(y)\,dV_y, \quad \,
        {\mathbf{w}}_0\in  L^2(\Omega, T\Omega).\end{eqnarray*}

Let $\big\{{\mathbf{u}}_k^{\begin{tiny}(D)\end{tiny}}\big\}_{k=1}^\infty$ (respectively,  $\big\{{\mathbf{u}}_k^{\begin{tiny}(T)\end{tiny}}\big\}_{k=1}^\infty$) be the orthonormal eigenvectors of the elastic operator $L_{\lambda,\mu}^{\begin{tiny}(D)\end{tiny}}$ (respectively, $L_{\lambda,\mu}^{\begin{tiny}(T)\end{tiny}}$) corresponding to the eigenvalues $\big\{\tau_k^{\begin{tiny}(D)\end{tiny}}\big\}_{k=1}^\infty$ (respectively, $L_{\lambda,\mu}^{\begin{tiny}(T)\end{tiny}}$), then the integral kernels  ${\mathbf{K}}^{\begin{tiny}(D)/(T)\end{tiny}}(t, x, y)=e^{-t L_{\lambda, \mu}^{\begin{tiny}(D)/(T)\end{tiny}}} 
\boldsymbol{\delta}(x-y)$ are given by \begin{eqnarray} \label{18/12/18} {\mathbf{K}}^{\begin{tiny}(D)/(T)\end{tiny}}(t,x,y) =\sum_{k=1}^\infty e^{-t \tau_k^{\begin{tiny}(D)/(T)\end{tiny}}} {\mathbf{u}}_k^{\begin{tiny}(D)/(T)\end{tiny}}(x)\otimes {\mathbf{u}}_k^{\begin{tiny}(D)/(T)\end{tiny}}(y).\end{eqnarray}
This implies that the integrals of the traces of ${\mathbf{K}}^{\begin{tiny}(D)/(T)\end{tiny}}(t,x,y)$ are actually  spectral invariants:
\begin{eqnarray} \label{1-0a-2}\int_{\Omega} \mbox{Tr}\,\big({\mathbf{K}}^{\begin{tiny}(D)/(T)\end{tiny}}(t,x,x)\big) \,dV_x=\sum_{k=1}^\infty e^{-t \tau_k^{\begin{tiny}(D)/(T)\end{tiny}}}.\end{eqnarray}

 We will combine calculus of symbols (see \cite{Gru}) and ``method of images'' to deal with asymptotic expansions for the trace integrals of integral kernels.
Let $\mathcal{M}=\Omega \cup (\partial \Omega)\cup \Omega^*$ be the (closed) double of $\Omega$, and $\mathcal{L}_{\lambda, \mu}$ the double to $\mathcal{M}$ of
 the  operator $L_{\lambda, \mu}$ on $\Omega$.

Let us explain the double Riemannian manifold $\mathcal{M}$ and the double differential operator $\mathcal{L}_{\lambda,\mu}$ more precisely, and introduce how to get them from the given Riemannian manifold $\Omega$ and the Lam\'{e} operator $L_{\lambda,\mu}$.
The double of $\Omega$ is the manifold $\Omega \cup_{\mbox{Id}} \Omega$, where $\mbox{Id}: \partial \Omega\to \partial \Omega$ is the identity map of $\partial \Omega$; it is obtained from $\Omega \sqcup \Omega$ by identifying each boundary point in one  copy of $\Omega$
with same boundary point in the other. It is a smooth manifold without boundary, and contains two regular domains diffeomorphic to $\Omega$  (see, p.$\,$226 of \cite{Lee}). When considering the double differential system $\mathcal{L}_{\lambda,\mu}$ crossing the boundary, we make use of the coordinates as follows.  Let $x'=(x_1, \cdots, x_{n-1})$ be any local coordinates for $\partial \Omega$.  For each point $(x',0)\in \partial \Omega$, let $x_{n}$ denote the parameter along the unit-speed geodesic starting at $(x',0)$ with initial direction given by the inward boundary normal to $\partial \Omega$ (Clearly, $x_n$ is the geodesic distance from the point $(x',0)$ to the point $(x',x_n)$). In such coordinates $x_{n}>0$ in
 $\Omega$, and $\partial \Omega$ is locally characterized by $x_{n}=0$ (see, \cite{LU} or \cite{Ta2}).
 Since the Lam\'{e} operator is a linear differential operator defined on $\Omega$, it can be further denoted as (see  (\ref{24.4.12-9})) $L_{\lambda,\mu}:=L_{\lambda,\mu}\big(x, \{g^{jk}(x)\}_{1\le j,k\le n}$,
$\{\Gamma^j_{kl}(x)\}_{1\le j,k,l\le n}$, $\{\frac{\partial \Gamma^s_{jk}(x)}{\partial  x_l}\}_{1\le s,j,k,l\le n}$, $\{R^j_k(x)\}_{1\le j,k\le n}$, $\frac{\partial}{\partial x_1}, \cdots, \frac{\partial}{\partial x_{n-1}}$,$\frac{\partial }{\partial x_n}\big)$.
  Let  $\varpi: (x_1, \cdots, x_{n-1}, x_n)\mapsto (x_1, \cdots, x_{n-1}, -x_n)$ be the reflection with respect to the boundary $\partial \Omega$ in $\mathcal {M}$ (here we always assume $x_n\ge 0$).  Then  we can get the $\Omega^*$ from the given $\Omega$ and $\varpi$.
    Now, we discuss the change of the metric $g$ from $\Omega$ to $\Omega^*$ by $\varpi$. Recall that the Riemannian metric $(g_{ij})$ is given in the local coordinates $x_1, \cdots, x_n$, i.e., $g_{ij}(x_1,\cdots, x_n)$. In terms of the new coordinates $z_1,\cdots, z_n$, with  $x_i=x_i(z_1,\cdots, z_n), \,\, i=1,\cdots, n,$  the same metric is given by the functions $\tilde{g}_{ij} =\tilde{g}_{ij} (z_1, \cdots, z_n)$, where
\begin{eqnarray} \label{2022.10-2} \tilde{g}_{ij}= \sum_{k,l=1}^n \frac{\partial x_k}{\partial z_i} g_{kl} \frac{\partial x_l}{\partial z_j}.\end{eqnarray}
  If $\varpi$ is a coordinate change in a neighborhood intersecting with $\partial \Omega$
  \begin{eqnarray}\label{2022.10.28-1} \left\{ \begin{array}{ll} x_1= z_1, \\
                   \cdots \cdots\\
                   x_{n-1}=z_{n-1},\\
                   x_n=-z_n,\end{array}\right.\end{eqnarray}
                                       then its Jacobian matrix is
\begin{eqnarray} \label{2022.11.8-1} J:=\begin{pmatrix} \frac{\partial x_1}{\partial z_1}& \cdots &   \frac{\partial x_1}{\partial z_{n-1}} & \frac{\partial x_1}{\partial z_{n}}\\
   \vdots & \ddots & \vdots & \vdots\\
  \frac{\partial x_{n-1}}{\partial z_{1}} & \cdots & \frac{\partial x_{n-1}}{\partial z_{n-1}} &  \frac{\partial x_{n-1}}{\partial z_{n}}\\
   \frac{\partial x_n}{\partial z_{1}}& \cdots & \frac{\partial x_n}{\partial z_{n-1}}&\frac{\partial x_n}{\partial z_{n}}\end{pmatrix} =   \begin{pmatrix} 1& \cdots &   0 & 0\\
   \vdots & \ddots & \vdots & \vdots\\
  0 & \cdots & 1 &  0\\
   0& \cdots & 0&-1\end{pmatrix}.\end{eqnarray}
Using this and (\ref{2022.10-2}), we immediately obtain the corresponding metric on the $\Omega^*$:
  (see \cite{MS-67},  or p.\,10169, p.\,10183 and p.\,10187 of \cite{Liu2})
\begin{eqnarray} \label{2021.2.6-3}  g_{jk} (\overset{*}{x})\!\!\!&\!=\!&\!\!\!- g_{jk} (x) \quad \, \mbox{for}\;\;
  j<k=n \;\;\mbox{or}\;\; k<j=n,\\ g(\overset{*}{x}) \!\!\!&\!=\!&\!\!\! g_{jk} (x)\;\; \;\;\mbox{for}\;\; j,k<n \;\;\mbox{or}\;\; j=k=n,\\
  \label{2021.2.6-4}  g_{jk}(x)\!\!\!&\!=\!&\!\!\! 0 \;\; \;\mbox{for}\;\; j<k=n \;\;\mbox{or}\;\; k<j=n \;\;\mbox{on}\;\; \partial \Omega,\end{eqnarray} where $x_n(\overset{*}{x})= -x_n (x)$.
  We denote such a new (isometric) metric on $\Omega^*$ as $g^*$.
  It is easy to verify that \begin{eqnarray*} \begin{bmatrix} g_{11}(x) & \cdots & g_{1,n-1}(x)& -g_{1n}(x)\\
  \vdots & \ddots & \vdots & \vdots\\
  g_{n-1,1}(x) & \cdots & g_{n-1,n-1}(x) &- g_{n-1,n}(x)\\
 - g_{n1}(x) & \cdots & -g_{n,n-1}(x) & g_{nn}(x)\end{bmatrix}^{-1}=\begin{bmatrix} g^{11}(x) & \cdots & g^{1,n-1}(x)& -g^{1n}(x)\\
  \vdots & \ddots & \vdots & \vdots\\
  g^{n-1,1}(x) & \cdots & g^{n-1,n-1}(x) &- g^{n-1,n}(x)\\
 - g^{n1}(x) & \cdots & -g^{n,n-1}(x) & g^{nn}(x)\end{bmatrix}, \end{eqnarray*}
  where $[g^{jk}(x)]_{n\times n}$ is the inverse of $[g_{jk}(x)]_{n\times n}$.
    In addition, by this reflection $\varpi$,
the differential operators $\frac{\partial }{\partial x_1}$, $\cdots$, $\frac{\partial}{\partial x_{n-1}}$, $\frac{\partial }{\partial x_n}$ (defined on  $\Omega$) are changed to $\frac{\partial }{\partial x_1}$, $\cdots$, $\frac{\partial}{\partial x_{n-1}}$, $-\frac{\partial }{\partial x_n}$ (defined on  $\Omega^*$), respectively. 
 It is easy to verify that
\begin{eqnarray}  \label{2023.2.22-8}\Gamma^{j}_{kl}(\overset{*}{x})=a_{jkl} \Gamma^{j}_{kl}(x),\end{eqnarray}
where \begin{eqnarray*} a_{jkl} =\left\{ \begin{array}{ll} 1 \;\;\; & \mbox{if there is no}\;\; n \;\; 
\mbox{among} \;\; j,k,l,\\
-1 \;\;\;\; & \mbox{if there is an}\;\; n \;\; 
\mbox{among} \;\; j,k,l,\\
1 \;\;\;\; & \mbox{if there are two}\;\; n \;\; 
\mbox{among} \;\; j,k,l,\\
-1 \;\;\; \; & \mbox{if there are three}\;\; n \;\; 
\mbox{among} \;\; j,k,l, \end{array}\right. \end{eqnarray*}
 \begin{eqnarray}\label{2022.5.22-9}
\frac{\partial \Gamma^s_{jk}}
{\partial {x}_l}(\overset{*}{x}) =b_{sjkl}\frac{\partial \Gamma^s_{jk}}{\partial x_l}(x),\end{eqnarray}
where \begin{eqnarray*} b_{sjkl}=\left\{ \begin{array} {ll} 1 \;\;\; &\mbox{if there is no}\;\; n \;\; \mbox{among}\;\; s,j,k,l,\\
-1\;\;\;& \mbox{if there is an}\;\; n\;\;  \mbox{among}\;\; s,j,k,l,\\
1\;\;  \;\; &\mbox{if there are two}\;\; n\;\; \mbox{among}\;\; s,j,k,l,\\
-1 \;\;\;& \mbox{if there are three}\;\; n \;\; \mbox{among}\;\; s,j,k,l,\\
1 \;\;\;& \mbox{if there are  four}\;\; n \;\;\mbox{among}\;\; s,j,k,l,\end{array}\right. \end{eqnarray*}
and \begin{eqnarray}\label{2023.5.22-10} R^j_k(\overset{*}{x}) =c_{jk} R^j_k(x),\end{eqnarray}
where  \begin{eqnarray*}
c_{jk}= \left\{ \begin{array} {ll}1 \;\;\; &\mbox{if there no} \;\; n \;\; \mbox{among}\;\; j,k, \\
- 1 \;\;\; &\mbox{if there is an} \;\; n\;\; \mbox{among}\;\; j,k, \\
1 \;\;\; &\mbox{if there are two} \;\; n  \;\;\mbox{among}\;\; j,k.\end{array}\right.\end{eqnarray*}
 In what follows, we will let Greek indices run from $1$ to $n-1$, Roman indices 
from $1$ to $n$.   
  We define \begin{eqnarray} \label{2022.10.18-2} \mathcal{L}_{\lambda,\mu}=\left\{\begin{array}{ll} \! L_{\lambda,\mu} \;\;\; \;\;\;\,\mbox{on} \;\, \Omega\\
 \!L^\star_{\lambda,\mu} \;\;\;\; \;\;\, \mbox{on} \;\, \Omega^*, \end{array} \right.\end{eqnarray}
where \begin{eqnarray} \label{2022.10.6-8}&& L_{\lambda,\mu}^\star:=L_{\lambda,\mu}\Big({g}^{\alpha\beta}(\overset {*}{x}), - {g}^{\alpha n}(\overset{*}{x}),- {g}^{n\beta}(\overset{*}{x}), {g}^{nn}(\overset{*}{x}),a_{jkl}\Gamma^{j}_{kl}(\overset{*}{x}), b_{sjkl}\frac{\partial \Gamma^s_{jk}}
{\partial{x}_l}(\overset{*}{x}),\\
&& \qquad \quad \;\quad c_{jk}R^j_k(\overset{*}{x}), \frac{\partial }{\partial x_1}, \cdots, \frac{\partial }{\partial x_{n-1}},- \frac{\partial }{\partial x_n}\Big),\nonumber \end{eqnarray}
and $\overset{*}{x}=(x',-x_n)\in \Omega^*$.  
 Roughly speaking, $L^\star_{\lambda,\mu}$ is obtained from the expression of $L_{\lambda,\mu}$ by replacing $\frac{\partial}{\partial x_n}$ by $-\frac{\partial}{\partial x_n}$. But we must rewrite such a $L^\star_{\lambda,\mu}$ in the terms  of 
the corresponding metric, Christoffel symbols and Ricci curvatures in $\Omega^*$. 
That is, $L^\star_{\lambda,\mu}$ is got if we replace $g^{\alpha\beta}(x)$, $g^{\alpha n}(x)$, $g^{n\beta} (x)$, $g^{nn}(x)$,
$\{\Gamma^j_{kl}(x)\}_{1\le j,k,l\le n}$, $\{\frac{\partial \Gamma^s_{jk}}{\partial  x_l}(x)\}_{1\le s,j,k,l\le n}$, $\{R^j_k(x)\}_{1\le j,k\le n}$,
 $\frac{\partial }{\partial x_n}$ by ${g}^{\alpha\beta}(\overset{*}{x})$, $-{g}^{\alpha n}(\overset{*}{x})$, $-{g}^{n\beta} (\overset{*}{x})$, ${g}^{nn}(\overset{*}{x})$,
$\;a_{jkl}\Gamma^{j}_{kl}(\overset{*}{x})$, $\,b_{sjkl}\frac{\partial \Gamma^s_{jk}}
{\partial {x}_l}(\overset{*}{x})$, $c_{jk}R^j_k(\overset{*}{x})$, 
 $-\frac{\partial }{\partial x_n}$ in $L_{\lambda,\mu}=L_{\lambda,\mu}\Big( g^{\alpha\beta} (x)$, $g^{\alpha n}(x)$, $g^{n\beta}(x)$, $g^{nn}(x)$,
$\{\Gamma^j_{kl}(x)\}_{1\le j,k,l\le n}$, $\{\frac{\partial \Gamma^s_{jk}}{\partial  x_l}(x)\}_{1\le s,j,k,l\le n}$, $\{R^j_k(x)\}_{1\le j,k\le n}$,
 $\frac{\partial }{\partial x_1}$, $\cdots$, $\frac{\partial }{\partial x_{n-1}}$, $\frac{\partial }{\partial x_n}\Big)$, respectively. 
Note that we have used the relations (\ref{2021.2.6-3})--(\ref{2023.5.22-10}).
In view of the metric matrices $g$ and $g^*$ have the same order principal minor determinants, we see that $\mathcal{L}_{\lambda,\mu}$ is still a linear elliptic differential operator on closed Riemannian manifold $\mathcal{M}$.

Let $\mathbf{K}(t,x,y)$ be the fundamental solution of the parabolic  system
 \begin{eqnarray*} \left\{ \begin{array}{ll} \frac{\partial \mathbf{u}}{\partial t} + \mathcal{L}_{\lambda,\mu}\mathbf{u}=0 \;\; &\mbox{in}\;\, (0,+\infty)\times \mathcal{M},\\
  \mathbf{u}=\boldsymbol{\phi} \;\; &\mbox{on}\;\; \{0\}\times \mathcal{M}.\end{array}\right.\end{eqnarray*}
  That is, for any $t\ge 0$ and $x,y\in \mathcal{M}$,
\begin{eqnarray}\label{2020.10.28-10}\left\{\begin{array}{ll}
    \frac{\partial \mathbf{K}(t,x,y)}{\partial t} + \mathcal{L}_{\lambda,\mu}\mathbf{K}(t,x,y)=0 \;\;  &\mbox{for}\;\, t>0, \, x, y\in \mathcal{M},\\
      \mathbf{K}(0,x,y)=\boldsymbol{\delta}(x-y) \;\; &\mbox{for}\;\;  x,y \in \mathcal{M}.\end{array}\right.\end{eqnarray}
   Here the operator $\mathcal{L}_{\lambda,\mu}$ is acted in the third argument $y$ of $\mathbf{K}(t,x,y)$.

 Clearly, the coefficients occurring in $\mathcal{L}_{\lambda,\mu}$ jump as $x$ crosses the $\partial \Omega$ (since the extended metric $g$ is $C^0$-smooth on whole $\mathcal{M}$ and $C^\infty$-smooth in $\mathcal{M}\setminus \partial \Omega$), but $\frac{\partial \mathbf{u}}{\partial t}+\mathcal{L}_{\lambda,\mu} \mathbf{u}=0$ with $\mathbf{u}(0,x)=\boldsymbol{\phi}(x)$ still has a nice fundamental solution $\mathbf{K}$  of class $[C^1((0,+\infty)\times \mathcal{M}\times  \mathcal{M})]_{n\times n} \cap [ C^\infty((0,+\infty)\times (\mathcal{M} \setminus \partial \Omega) \times (\mathcal{M}\setminus \partial \Omega))]_{n\times n}$, approximable even on $\partial \Omega$  by Levi's sum (see \cite{Liu2}, or another proof below). Now, let us restrict $x,y\in \Omega$.
 Recall that for an (elastic) vector field $\mathbf{u}$ defined in $\Omega$, the boundary traction operator can also  be equivalently written as 
 \begin{eqnarray} \label{23.12.12-1} \;\;\;\quad \; (\mathcal{T}_{\lambda,\mu} \mathbf{u})^k :=\mu \sum_{l=1}^n\Big( \nu^{\,l} \nabla_l u^k + \nu_l \nabla^k  u^l\Big)+  \lambda\, \boldsymbol{\nu}^k
\sum_{l=1}^n \nabla_l u^l \; \;\mbox{on} \;\, \partial \Omega, \,\;\; k=1,\cdots, n.\;\;\end{eqnarray}  
  Write  $\mathbf{K} (t,x,y)=\big(K^{jk} (t,x,y)\big)_{n\times n}$,  $ \;\mathbf{K}_j(t, x, y)= (K^{j1}(t,x,y)$, $\cdots$, $K^{jn} (t,x,y)$,  $\;(j=1,\cdots, n)$,   and denote   \begin{eqnarray}\label{23.12.13-2} \mathcal{T}_{\lambda,\mu} \mathbf{K} (t,x,y) =\big( \mathcal{T}_{\lambda,\mu}  \mathbf{K}_1 (t,x,y), \cdots,  \mathcal{T}_{\lambda,\mu}  \mathbf{K}_n(t,x,y)\big)  \end{eqnarray}  
 for all $t>0, x\in \Omega, y\in \partial \Omega$.  More precisely,  for each $j=1,\cdots, n$, the $k$-th component of $\mathcal{T}_{\lambda,\mu}  \mathbf{K}_j$ is 
 \begin{eqnarray} \label{23.12.13-6}  (\mathcal{T}_{\lambda,\mu} \mathbf{K}_j)^k = \mu \sum_{l=1}^n ( \nu^l \nabla_l K^{jk} +\nu_l \nabla^k K^{jl})
  +\lambda \nu^k \sum_{l=1}^n\nabla_l K^{jl}.\end{eqnarray}
 It is easy to see that $\frac{\partial \big(\mathbf{K}(t,x,y)+\mathbf{K}(t,x,\overset{*}{y})\big)}{\partial \nu_y}\big|_{\partial \Omega}=0$
for  all $t>0$,  $x\in \Omega$ and $y\in \partial \Omega$.  Changing all  terms $\frac{\partial K^{jk}(t,x,y)}{\partial \nu_y}\big|_{\partial \Omega}$ into $0$ in the expression  of $2\mathcal{T}_{\lambda,\mu} \mathbf{K}(t,x,y)$ (i.e., replacing  all terms $\frac{\partial K^{jk}(t,x,y)}{\partial \nu_y}\big|_{\partial \Omega}= \sum_{m=1}^n \frac{\partial K^{jk}(t,x,y)}{\partial y_m} \, \nu_m$ by $0$ in the expression of   $2\mathcal{T}_{\lambda,\mu} \mathbf{K}(t,x,y)$),  we obtain  a matrix-valued function $\boldsymbol{\Upsilon}(t,x,y)$  for $t>0$, $x\in \Omega$ and $y\in \partial \Omega$.  This implies that  $\boldsymbol{\Upsilon}(t,x,y)$ only contains the (boundary) tangent derivatives of $K^{jk}(t, x, y)$ with respect to $y\in \partial\Omega$ (without normal derivative of $K^{jk}(t,x,y)$) in the local expression). That is,  in local boundary normal coordinates (the inner normal $\boldsymbol{\nu}$ of $\partial \Omega$ is in the direction of  $x_n$-axis), \begin{eqnarray*}  \boldsymbol{\Upsilon} (t,x,y)= \Big( \boldsymbol{\Upsilon}_1 (t,x,y), \cdots,   \boldsymbol{\Upsilon}_n (t,x,y)\Big), \end{eqnarray*}
   \begin{eqnarray*}\!\!\!\!\!\!\! \!&\!\!\!&\!\!\! \boldsymbol{\Upsilon}_j (t,x,y) \\
\!\!\!\!\!\! \!\!\! &\!\!\!\!\!&\!\!\! \! =\!
  2\mu \!\begin{small}   \begin{pmatrix} \! 2\frac{\partial K^{j1}}{\partial x_1}+2\!\sum\limits_{m=1}^n \!\Gamma^{1}_{1m} K^{jm}  &\! \!\!\cdots\!  \!  & \!\!\!\frac{\partial K^{j1}}{\partial x_{n\!-\!1}} \!\!+\!\!\frac{\partial K^{\!j,n\!-\!1}}{\partial x_1} \!\!+\!\!\! \sum\limits_{m=1}^n\!\! \big( \Gamma^1_{\!n\!-\!1, m}\! \!+\!\!\Gamma^{n\!-\!1}_{\!1m} \big) K^{jm} 
            &   \frac{\partial K^{jn}}{\partial x_1}\!\!+\!\!\!\sum\limits_{m=1}^n\!\! \big( \Gamma^1_{\!nm} \!\!+\!\!\Gamma^n_{1m}\big) K^{jm}  
      \\
\!\!\! \!\cdots \!& \!\cdots \! \!&\!\!\cdots  & \cdots \\
 \!  \! \frac{\partial K^{\!j,n\!-\!1}}{\partial x_{1}}\!\!+\! \!\frac{\partial K^{j1}} {\partial x_{n\!-\!1}}\! \!+\! \!\!\sum\limits_{m=1}^n \!\! \big( \Gamma^{n\!-\!1}_{\!1m}\!\! +\!\! \Gamma^1_{\!n\!-\!1, m} \big)  K^{jm}
 \! \!&\! \cdots    &\!\!\!2 \frac{\partial K^{j,n\!-\!1}}{\partial x_{n\!-\!1}}\!\!+\!\!2\!\sum\limits_{m=1}^n \!\!\Gamma^{n\!-\!1}_{\!n\!-\!1,m} K^{jm} &  \! \frac{\partial K^{jn}}{\partial x_{n\!-\!1}}\!\!+\!\!\sum\limits_{m=1}^n\!\! \big(\Gamma^{n\!-\!1}_{\!\!nm}\! \!+\!\!\Gamma_{\!\!n\!-\!1,m}^n \!\big) K^{\!jm} \!\!
   \\
  \! \frac{\partial K^{jn}}{\partial x_1} \!+\!\!\sum\limits_{m=1}^n \!\!\big( \Gamma_{\!1m}^n\! +\!\Gamma^1_{\!nm}\big) K^{jm}   \! &\!\!\! \cdots  \! \! & 
\frac{\partial K^{jn}}{\partial x_{n\!-\!1}}\! +\!\big(\Gamma^n_{n-1,m}+\Gamma_{\!nm}^{n\!-\!1} \big)K^{jm}   & 
    \;2 \sum\limits_{m=1}^n 
 \Gamma^n_{nm}  K^{jm}\!\!
 \end{pmatrix}\end{small} \!\!\begin{pmatrix} \nu_1 \\ \vdots \\ \nu_n \end{pmatrix}\\  [3mm]
\!\!\!\!\!\!&\!\!\!&\!\!\! \;\;\; + \, 2\lambda\bigg( \frac{\partial K^{j1}}{\partial x_1} + \sum_{m=1}^n \Gamma^{1}_{1m} K^{jm}  + \cdots + \frac{\partial K^{j,n-1}}{\partial x_{n-1}} + \sum_{m=1}^n \Gamma^{n-1}_{n-1,m} K^{jm} +  \sum_{m=1}^n \Gamma^{n}_{nm} K^{jm} \bigg) \begin{pmatrix} \nu_1 \\ \vdots \\ \nu_n \end{pmatrix},\;\;\;  \; j=1,\cdots, n.\end{eqnarray*}
 It is easy to see that  $\boldsymbol{\Upsilon} (t,x,y)$  is a continuous (matrix-valued) function for all $t> 0$, $x\in \Omega$ and $y\in \partial \Omega$. Further, for any  fixed $x\in \Omega$ and any $y\in \partial \Omega$,  since $x\ne y$ we see that 
 \begin{eqnarray} \label{23.12.18-1}  \lim\limits_{t\to 0^+} \mathbf{K} (t,x,y)=0,\end{eqnarray} 
 which implies \begin{eqnarray} \label{23.12.18-2} \lim\limits_{t \to 0^{+}} \frac{\partial \mathbf{K}(t,x,y)}{\partial T} \big|_{\partial \Omega} =0\;\;\;\mbox{for}\;\;\,   x\in \Omega, \;\;\, y\in \partial \Omega,  \end{eqnarray} 
  where $\frac{\partial}{\partial T}$ denotes the tangent derivative along the boundary $\partial \Omega$ in variable $y$. 
 From (\ref{23.12.18-1})--(\ref{23.12.18-2}) we get 
 \begin{eqnarray*} \label{23.12.18-3}  \lim\limits_{t\to 0^+} \nabla_l \mathbf{K} (t,x,y) \big|_{\partial \Omega} =0 \;\; \; \mbox{for all}\;\; x \in \Omega, \; y\in \partial \Omega, \;\;\, 1\le l <n,\end{eqnarray*}
   so that  \begin{eqnarray} \label{23.12.18-4} \lim\limits_{t\to 0^{+}} \boldsymbol{\Upsilon}(t,x,y)=0  \;\;\,\mbox{for any } x\in \Omega, 
   \; y\in \partial \Omega,\end{eqnarray}
    where $\nabla_l \mathbf{K}= (\nabla_l \mathbf{K}_1, \cdots, \nabla_l \mathbf{K}_n)$ and $\nabla_l {K}^{jk}:= \frac{\partial 
       K^{jk} (t,x,y)}{\partial y_l}  +\sum_{m=1}^n\Gamma^k_{lm} K^{jm} (t,x,y)$,  $\,(1\le j,k,l\le  n)$.
     Let $\mathbf{H}(t,x,y)$ be the solution of  
    \begin{eqnarray*} \left\{ \begin{array}{ll}  \frac{\partial \mathbf{u} (t,x,y)}{\partial t} = L_{\lambda,\mu} \mathbf{u} (t,x,y) \; \;\;\mbox{for}\;\,  t>0, \; x,y\in \Omega,\\
2  \mu  \big(\mbox{Def}\, \mathbf{u}(t,x,y)\big)^\#  \,\boldsymbol{\nu}  +\lambda \big(\mbox{div}\, (\mathbf{u}(t,x,y))\big)\, \boldsymbol{\nu} =\boldsymbol{\Upsilon}(t,x,y)\,  \;\;\mbox{for} \,\;  t>0, \, x\in \Omega, \, y\in \partial \Omega,\\
  \mathbf{u}(0, x,y) = \mathbf{0} \,\;\;\mbox{for} \,\;  x, y \in \Omega.\end{array} \right. \end{eqnarray*}
  From (\ref{23.12.18-4}), we get that the above parabolic system satisfy the compatibility condition.   
 Thus,  the  matrix-valued solution $\mathbf{H}(t,x,y)$, is  smooth in $(0, \infty) \times  \Omega \times  \Omega$ and  continuous on  $[0, \infty) \times \bar \Omega \times \bar \Omega$.   Then there exists a constant $C>0$ such that 
  \begin{eqnarray*}  | \mathbf{H} (t,x,y) |\le C \,\, \;\;\mbox{for all} \;\;   0\le t\le 1, \; x,\,y\in \bar \Omega, \;\;\,\end{eqnarray*}
Actually, we have  $\,\lim_{t \to 0^{+}} \!\int_{\Omega} \mbox{Tr}\; \mathbf{H} (t,x,x)\, dx=0 $,  
   and hence, for dimensions $n\ge 2$,   \begin{eqnarray}\label{23.12.9-1} \int_{\Omega} \mbox{Tr}\; \mathbf{H} (t,x,x)\, dx =o(nC(\mbox{vol}(\Omega))) =o(t^{-\frac{n-2}{2}})\,  \;\;\mbox{as}\;\; t\to 0^+.\end{eqnarray}
(\ref{23.12.9-1})  implies that we can precisely get the first three  terms asymptotic expansion. 
   It  can easily be verified  that 
\begin{eqnarray*} 
 && \mathbf{K}^{\begin{tiny}(D)\end{tiny}} (t,x,y)= \mathbf{K}(t,x,y)- \mathbf{K}(t,x,\overset{*}{y}), \\
&& \mathbf{K}^{\begin{tiny}(T)\end{tiny}} (t,x,y)= \mathbf{K}(t,x,y)+ \mathbf{K}(t,x,\overset{*}{y})-  \mathbf{H}(t,x,y) ,
\end{eqnarray*}   
 are the Green functions of \begin{eqnarray*}\left\{ \begin{array}{ll}\frac{\partial \mathbf{u}}{\partial t} +{L}_{\lambda,\mu}\mathbf{u}=0\;\;&\mbox{in}\;\, (0,+\infty) \times \Omega,\\
\mathbf{u}=\boldsymbol{\phi} \;\; &\mbox{on}\;\; \{0\}\times \Omega\end{array}\right.\end{eqnarray*} with zero Dirichlet and zero traction (i.e., free) boundary conditions, respectively,
where $y=(y',y_n)$, $y_n\ge 0$, and  $\overset{*}{y}:=\varsigma(y',y_n)=(y^{\prime},-y_n)$.  In other words,
 \begin{eqnarray*}\left\{    \begin{array}{ll} \frac{\partial \mathbf{K}^{\begin{tiny}(D)\end{tiny}}(t,x,y)}{\partial t} + {L}_{\lambda,\mu}\mathbf{K}^{\begin{tiny}(D)\end{tiny}}(t,x,y)=0,\;\;\; t>0,\, x,\, y\in\Omega,\\
 \mathbf{K}^{\begin{tiny}(D)\end{tiny}} (t, x,y)=0, \;\;\;\;  t>0,\; x\in \Omega, \,\;  y\in \partial \Omega,\\
 \mathbf{K}^{\begin{tiny}(D)\end{tiny}}(0, x,y)=\boldsymbol{\delta}(x-y), \;\;\;  x,y\in \Omega\end{array} \right. \end{eqnarray*}
and
 \begin{eqnarray*}\left\{  \begin{array}{ll}   \frac{\partial \mathbf{K}^{\begin{tiny}(T)\end{tiny}}(t,x,y)}{\partial t} + {L}\mathbf{K}^{\begin{tiny}(T)\end{tiny}}(t,x,y)=0, \,\;\;\; t>0, \;x, \,y\in \Omega,\\
 \mathcal {T}_{\lambda,\mu} (\mathbf{K}^{\begin{tiny}(T)\end{tiny}} (t, x,y))=0, \;\; \,\,t>0,\;  x\in \Omega, \;\;  y\in \partial \Omega,\\
 \mathbf{K}^{\begin{tiny}(T)\end{tiny}}(0, x,y)=\boldsymbol{\delta}(x-y),\;\;\;  x,y\in \Omega,\end{array} \right. \end{eqnarray*}
 where $ \mathcal {T}_{\lambda,\mu} (\mathbf{K}^{\begin{tiny}(T)\end{tiny}} (t, x,y))$ is the traction of  $\mathbf{K}^{\begin{tiny}(T)\end{tiny}} (t, x,y) $ on $\partial \Omega$.
   By combining  the fact  $\frac{\partial (\mathbf{K} (t,x,y)+ \mathbf{K}(t,x,\overset{*}{y}))}{\partial \boldsymbol{\nu}_y}\big|_{\partial \Omega} =0$  and $\mathcal{T}_{\lambda,\mu}  \big(\mathbf{K} (t,x,y)+\mathbf{K}(t, x, \overset{*}{y})\big)=\boldsymbol{\Upsilon} (t,x,y) =\mathcal{T}_{\lambda,\mu}  \mathbf{H} (t,x,y)$ for $t>0$,  $x\in \Omega$ and $y\in \partial \Omega$,  we  get $\mathcal{T}_{\lambda,\mu}  \mathbf{K}^{\begin{tiny}(T)\end{tiny}} (t,x,y) =0$ for all $t>0$, $x\in \Omega$ and $y\in \partial \Omega$. 
  Now, we will show  $\mathbf{K}^{\begin{tiny}(D)\end{tiny}}(t,x,y)$ and $\mathbf{K}^{\begin{tiny}(T)\end{tiny}}(t,x,y)$ satisfy the associated parabolic system.  
  In fact, for any $t>0$, $x,y \in \Omega$, we have
  $L_{\lambda,\mu}\mathbf{K}(t,x, y)=\mathcal{L}_{\lambda,\mu}\mathbf{K}(t,x, y)$,  so that
  \begin{eqnarray}\label{2022.10.29-8}\left\{ \begin{array}{ll}
  \Big( \frac{\partial}{\partial t}+{L}_{\lambda,\mu}\Big) \mathbf{K} (t,x, y)=
  \Big( \frac{\partial}{\partial t}+\mathcal{L}_{\lambda,\mu}\Big) \mathbf{K} (t,x, y)=0,\\
  \mathbf{K}(0, x,y)=\boldsymbol{\delta}(x-y) \end{array}\right.\end{eqnarray}
  by (\ref{2020.10.28-10}).
Noting that the  Jacobian matrix of the reflection $\varpi$ is $J$ (see (\ref{2022.11.8-1})), it follows from chain rule that for any fixed $t>0$ and $x\in \Omega$, and any $y=(y',y_n)\in \Omega$,
\begin{align*} &\left[{L}_{\lambda,\mu} ( \mathbf{K}(t,x,\overset{*}{y}))\right]\bigg|_{\text{\normalsize evaluated at the point $y$}}\\
&=\left[{L}_{\lambda,\mu} ( \mathbf{K}(t,x,\varsigma(y',y_n)))\right]\bigg|_{\text{\normalsize evaluated at the point $(y',y_n)$}}
 \\
&= \left[{L}_{\lambda,\mu} (\mathbf{K}(t, x, (y',-y_n))\big)\right]\bigg|_{\text{\normalsize evaluated at the point $(y',y_n)$}}\\
& \! =\!
\left\{\left[L_{\lambda,\mu} \big(g^{\alpha\beta}(y), g^{\alpha n}(y), g^{n\beta} (x), g^{nn}(x),
\big\{\Gamma^j_{kl}(x)\big\}_{1\le j,k,l\le n}, \Big\{\frac{\partial \Gamma^s_{jk}}{\partial  x_l}(x)\Big\}_{1\le s,j,k,l\le n}, \right.\right.\\
&
\left.\left. \quad\quad  \;\big\{R^j_k(x)\big\}_{1\le j,k\le n}, \,\frac{\partial}{\partial y_1}, \cdots, \frac{\partial}{\partial y_{n-1}}, \frac{\partial}{\partial y_n}\big)\right]\!  \mathbf{K} (t, x, (y'\!,-y_n))\!\right\}\!\Bigg|_{\text{\normalsize evaluated at $\!(y'\!,y_n)$}}\\
\!& \! =
 \left\{\left[L_{\lambda,\mu} \big(g^{\alpha\beta}(\overset{*}{y}), -g^{\alpha n}(\overset{*}{y}), -g^{n\beta}(\overset{*}{y}), g^{nn}(\overset{*}{y}), 
\big\{a_{jkl}\Gamma^{j}_{kl}(\overset{*}{x})\big\}_{1\le j,k, l\le n},\Big\{b_{sjkl}\frac{\partial \Gamma^s_{jk}}
{\partial {x}_l}(\overset{*}{x})\Big\}_{1\le s,j,k,l\le n},\right. \right.\\
& \left.\left.\quad \quad \; \;  \big\{c_{jk}R^j_k(\overset{*}{x})\big\}_{1\le j,k\le n},\,
 \frac{\partial}{\partial y_1}, \cdots, \frac{\partial}{\partial y_{n-1}}, -\frac{\partial}{\partial y_n}\big)\right]\!  \mathbf{K}(t,x,y)\!\right\}\Bigg|_{\text{\normalsize evaluated  at $\overset{*}{y}=(y',-y_n)$}} \\
&= L^\star_{\lambda,\mu} ( \mathbf{K}(t, x, \overset{*}{y}))\Big|_{\text{\normalsize evaluated at the point $\overset{*}{y}=(y',-y_n)$}}.  \end{align*}
 That is, the action of $L_{\lambda,\mu}$ to $\mathbf{K}(t,x,\overset{*}{y})$ (regarded as a vector-valued function of $y$) at the point ${y}=(y',y_n)$ is just the action of $L_{\lambda,\mu}^\star$
 to $\mathbf{K}(t,x,\overset{*}{y})$ (regarded as a vector-valued of $\overset{*}{y}$) at the point $\overset{*}{y}=(y',-y_n)$. Because of  $\varsigma(y',y_n)=(y',-y_n)\in \Omega^*$, we see $$L^\star_{\lambda,\mu}( \mathbf{K}(t,x,\overset{*}{y}))\big|_{\text{\normalsize evaluated at the point $\overset{*}{y}=(y'\!,-y_n)$}}\!=\!\mathcal{L}_{\lambda,\mu}(\mathbf{K} (t,x,\overset{*}{y}))\big|_{\text{\normalsize evaluated at the point $\overset{*}{y}\!=\!(y'\!,-y_n)$}}.$$
For any $t>0$, $x\in \Omega$ and $(y',-y_n)\in \Omega^*$, we have $$\Big(\frac{\partial}{\partial t}+ {\mathcal{L}_{\lambda,\mu}}\Big) \big(\mathbf{K} (t, x,  (y',-y_n))\big)=0.$$
In addition, $\mathbf{K} (t, x,  (y',-y_n))= \mathbf{K} (t, x,  \varsigma(y))$ for any $t>0$, $x,y\in \Omega$. By virtue of $x\ne (y',-y_n)$,
 this leads to $\mathbf{K} (0, x,  (y',-y_n))=0$ and
 \begin{eqnarray*} \label{2020.10.29-1} \Big(\frac{\partial}{\partial t}+ {L}_{\lambda,\mu}\Big) \big(\mathbf{K} (t, x,  (y',-y_n))\big)=0\;\; \mbox{for any}\;\,t>0, x\in \Omega \;\,\mbox{and}\;\,(y',-y_n)\in \Omega^*,\end{eqnarray*}
 i.e.,
  \begin{eqnarray} \label{2020.10.29-2}\left\{\! \begin{array}{ll}  \big(\frac{\partial}{\partial t}+ {{L}_{\lambda,\mu}}\big) \mathbf{K} (t, x,  \overset{*}{y})=0\;\; \mbox{for any}\;\,t>0, x\in \Omega \;\,\mbox{and}\;\,\overset{*}{y}\in \Omega^*,\\
   \mathbf{K} (0, x,  \overset{*}{y})=0 \;\;\, \mbox{for any}\,\; x,y\in \Omega.\end{array} \right.\end{eqnarray}
 Combining (\ref{2022.10.29-8}) and (\ref{2020.10.29-1}), we obtain that
 \begin{eqnarray} \label{2020.10.29-3} \left\{ \begin{array}{ll} \big(\frac{\partial}{\partial t}+ {{L}_{\lambda,\mu}}\big) \,\Big(\mathbf{K} (t, x, {y}) \mp \mathbf{K} (t, x,  \overset{*}{y})\Big)=0\;\;\mbox{for any}\;\,t>0, \; x,y\in \Omega,\\
 \mathbf{K} (0, x, {y}) \mp \mathbf{K} (0, x,  \overset{*}{y})=\boldsymbol{\delta} (x-y)\;\;\mbox{for any}\;\, x,y\in \Omega.\end{array}\right.\end{eqnarray}
  $\mathbf{K}(t,x, y)$ is $C^1$-smooth with respect to $y$ in $\mathcal{M}$ for any fixed $t>0$ and $x\in \Omega$ (see below),   so does it on the hypersurface $\partial \Omega$. Therefore, we get that  $\mathbf{K}^{\begin{tiny}(D)\end{tiny}} (t,x,y)$ (respectively $\mathbf{K}^{\begin{tiny}(T)\end{tiny}}(t,x,y)$) is the Green function in $\Omega$ with the Dirichlet (respectively, free) boundary condition on $\partial \Omega$.

\vskip 0.12 true cm

To show $C^{1}$-regularity of the fundamental solution $\mathbf{K}(t,x,y)$, it suffices to prove  $C^{1,1+\alpha}_{loc}$-regularity for a $W^{1,2}_2$ strong solution  $\mathbf{u}$ of the  parabilic system $\big(\frac{\partial}{\partial t} +\mathcal{L}_{\lambda,\mu}\big)\mathbf{u}=0$ in $(0,+\infty)\times \mathcal{M}$, where $W^{1,2}_p:= \{ \mathbf{u}\big| \mathbf{u},\frac{\partial \mathbf{u}}{\partial t}, D\mathbf{u}, D^2 \mathbf{u} \in L^p\} $.
 When the coefficients of an elliptic system are smooth on both sides of an $(n-1)$-dimensional hypersurface (may be discontinuous  crossing this hypersurface), the corresponding $C^{1,1+\alpha}$-regularity for solutions of a parabolic equation system is a special case of Dong's result (see Theorem 4 of p.$\,$141 in \cite{Don}).
In fact, Dong in \cite{Don} has given regularity results to the strong solutions for parabolic equation with  more general coefficients. 
This type of system arises from the problems of linearly elastic laminates and composite materials (see, for example, \cite{CKVC}, \cite{LiVo}, \cite{LiNi}, \cite{Don} and \cite{Xi-11}).  

Let us discuss the parabolic (elastic) system in more detail. In fact, on p.$\,$10177 in \cite{Liu2},  we have expressed the Lam\'{e} operator $L_{\lambda, \mu}$ in $\Omega$  as the form of components relative to local coordinates:
 \begin{eqnarray} \label{2023.50.22-12} \\
 L_{\lambda, \mu}\mathbf{u}=\!\!\!\!\!\!&&\!\!\!\left\{ -\mu\Big( \sum_{m,l=1}^n g^{ml} \frac{\partial^2 }{\partial x_m\partial x_l}\Big){\mathbf{I}}_n -(\mu+\lambda)
\begin{bmatrix} \sum\limits_{m=1}^n g^{1m} \frac{\partial^2}{ \partial x_m\partial x_1} & \cdots & \sum\limits_{m=1}^n g^{1m} \frac{\partial^2}{ \partial x_m\partial x_n}
  \\    \vdots & {} & \vdots  \\
 \sum\limits_{m=1}^n g^{nm} \frac{\partial^2}{ \partial x_m\partial x_1} & \cdots & \sum\limits_{m=1}^n g^{nm} \frac{\partial^2}{ \partial x_m\partial x_n}   \end{bmatrix} \right.\nonumber\\
[1.5mm]  && +\mu \Big( \sum\limits_{m,l,s=1}^n g^{ml} \Gamma_{ml}^s \frac{\partial}{\partial x_s} \Big){\mathbf{I}}_n - \mu  \begin{bmatrix}  \sum\limits_{m,l=1}^n 2g^{ml} \Gamma_{1m}^1 \frac{\partial }{\partial x_l} & \cdots &  \sum\limits_{m,l=1}^n 2g^{ml} \Gamma_{nm}^1 \frac{\partial }{\partial x_l}\\
\vdots & {}& \vdots \\
\sum\limits_{m,l=1}^n 2g^{ml} \Gamma_{1m}^n \frac{\partial }{\partial x_l} & \cdots & \sum\limits_{m,l=1}^n 2g^{ml} \Gamma_{nm}^n \frac{\partial }{\partial x_l}
\end{bmatrix}\nonumber \\
[1.5mm]  && -(\mu+\lambda) \begin{bmatrix}\sum\limits_{m,l=1}^n g^{1m} \Gamma_{1l}^l  \frac{\partial }{\partial x_m} & \cdots & \sum\limits_{m,l=1}^n g^{1m} \Gamma_{nl}^l  \frac{\partial }{\partial x_m} \\       \vdots & {} & \vdots  \\
  \sum\limits_{m,l=1}^n g^{nm} \Gamma_{1l}^l  \frac{\partial }{\partial x_m} & \cdots & \sum\limits_{m,l=1}^n g^{nm} \Gamma_{nl}^l  \frac{\partial }{\partial x_m} \end{bmatrix}\nonumber
 \\
[1.5mm]   &&  -\mu \begin{bmatrix} \sum\limits_{l,m=1}^n g^{ml} \big(  \frac{\partial \Gamma^1_{1l}}{\partial x_m} +  \Gamma_{hl}^1 \Gamma_{1m}^h - \Gamma_{1h}^1 \Gamma_{ml}^h \big) & \cdots &  \sum\limits_{l,m=1}^n g^{ml} \big(  \frac{\partial \Gamma^1_{nl}}{\partial x_m} +  \Gamma_{hl}^1 \Gamma_{nm}^h - \Gamma_{nh}^1 \Gamma_{ml}^h \big)
\\ \vdots & {} & \vdots \\
\sum\limits_{l,m=1}^n g^{ml} \big(  \frac{\partial \Gamma^n_{1l}}{\partial x_m} +  \Gamma_{hl}^n \Gamma_{1m}^h - \Gamma_{1h}^n \Gamma_{ml}^h \big) & \cdots &  \sum\limits_{l,m=1}^n g^{ml} \big(  \frac{\partial \Gamma^n_{nl}}{\partial x_m} +  \Gamma_{hl}^n \Gamma_{nm}^h - \Gamma_{nh}^n \Gamma_{ml}^h \big)
\end{bmatrix} \nonumber\\
[1.5mm] && \left.- (\mu+\lambda) \begin{bmatrix} \sum_{l,m=1}^n g^{1m} \frac{\partial \Gamma^l_{1l}}{\partial x_m}& \cdots & \sum_{l,m=1}^n g^{1m} \frac{\partial \Gamma^l_{nl}}{\partial x_m}\\
\vdots & {} & \vdots \\
\sum_{l,m=1}^n g^{nm} \frac{\partial \Gamma^l_{1l}}{\partial x_m}& \cdots & \sum_{l,m=1}^ng^{nm} \frac{\partial \Gamma^l_{nl}}{\partial x_m}\end{bmatrix} - \mu \begin{bmatrix}
R^1_1 & \cdots & R^1_n\\
\vdots & {} & \vdots \\
R^n_1 & \cdots & R^n_n\end{bmatrix} \right\}\begin{bmatrix} u^1\\
\vdots\\
u^n\end{bmatrix},  \nonumber \end{eqnarray}
where ${\mathbf{I}}_n$ is the $n\times n$ identity matrix.
The above expression $L_{\lambda,\mu} \mathbf{u}$ will played a key role in our proof.
From the local expression  (\ref{2023.50.22-12}) of $L_{\lambda,\mu}$, we see that the top-order coefficients of $\mathcal{L}_{\lambda,\mu}$ are not ``too bad'' since only the first $(n-1)$ coefficients of the $n$-th column in the second matrix in $L_{\lambda,\mu}$ defined on $\Omega$ are changed their signs in $L^\star_{\lambda,\mu}$ at the reflection points of $\Omega^*$. In other words, only  $(\sum_{m=1}^n g^{1m} (x)
\frac{\partial^2 u^n}{\partial x_m\partial x_n}$,$\cdots$, $\sum_{m=1}^n g^{n-1,m} (x)
\frac{\partial^2 u^n}{\partial x_m\partial x_n})^T$ in $\Omega$ is changed into $(-\sum_{m=1}^n g^{1m} (\overset{*}{x})
\frac{\partial^2 u^n}{\partial x_m\partial x_n}$,$\cdots$, $-\sum_{m=1}^n g^{n-1,m} (\overset{*}{x})
\frac{\partial^2 u^n}{\partial x_m\partial x_n})^T$ in $\Omega^*$ in the second term in $\mathcal{L}_{\lambda,\mu}$ (see  (\ref{24.4.9-1}) of $L_{\lambda,\mu}$).
  For any small coordinate chart $V\subset \mathcal{M}$, if  $V\subset \mathcal{M}\setminus (\partial \Omega)$, then the solution $\mathbf{u}$ belongs to $[C^\infty((0, +\infty)\times V)]^n \cap [W^{1,2}_2 ((0, +\infty)\times V)]^n$ since the coefficients of parabolic system 
\begin{eqnarray} \label{2023.6.13-3} \Big(\frac{\partial }{\partial t}+\mathcal{L}_{\lambda,\mu}\Big)\mathbf{u}=0\end{eqnarray} are smooth in $V$.
If the coordinate chart $V\subset \mathcal{M}$ and $V\cap \partial \Omega \ne \varnothing$, then we can find a (local) diffeomorphism $\Psi$ such that $\Psi(V)= U\subset \mathbb{R}^n$ and $\partial \Omega$ is mapped onto $U\cap \{x\in \mathbb{R}^n \big| x_n=0\}$ (i.e., $V\cap \partial \Omega$ is flatten into hyperplane  $\{x\in \mathbb{R}^n\big|x_n=0\}$ by $\Psi$).  By this coordinate transformation $\Psi$, the parabolic system (\ref{2023.6.13-3}) is changed into another parabolic system $(\frac{\partial}{\partial t}- L_{\lambda,\mu})\mathbf{v}=0$, whose coefficients are smooth on both sides of $(n-1)$-dimensional hyperplane  $\{x\in \mathbb{R}^n\big|x_n=0\}$ (may be discontinuous crossing this hyperplane, i.e., the coefficients of $L_{\lambda,\mu}$ have jump only on this hyperplane). It follows from Dong's regularity result (Theorem 4 of \cite{Don}) that $\mathbf{v}\in [C^{1,1+\alpha}((0,+\infty)\times U)]^n$, so that $\mathbf{u}\in  [C^{1,1+\alpha}((0,+\infty)\times V)]^n$. 
 Note that $\mathcal{M}$ is a compact closed Riemannian manifold. Thus we find from the above discussion  and (global) geometric analysis technique that for any initial value $\boldsymbol{\phi}\in [C^\infty (\mathcal{M})]^n$, there exists a (global) strong solution in $[W^{1,2}_2((0,+\infty)\times \mathcal{M})]^n$ which is in  $[C^{1,1+\alpha}((0,+\infty)\times \mathcal{M})]^n$.  Of course, this result  also holds to our case for $C^{1,1+\alpha}$-regularity of the fundamental solution $\mathbf{K} (t,x,y)$ on $(0,+\infty)\times \mathcal{M}\times \mathcal{M}$ (see \cite{Liu2}).

\vskip 0.20 true cm

Therefore, the integral kernels $\mathbf{K}^{\begin{tiny}(D)/(T)\end{tiny}}(t,x,y)$
 of $\frac{\partial \mathbf{u}}{\partial t}+L_{\lambda,\mu}^{\begin{tiny}(D)/(T)\end{tiny}}\mathbf{u}=0$ can be expressed on $(0,\infty)\times \Omega\times \Omega$ as
 \begin{eqnarray} \label{c4-23} && {\mathbf{K}}^{\begin{tiny}(D)\end{tiny}} (t,x,y) =\mathbf{K}(t,x,y)- \mathbf{K}(t,x,\overset{\ast} {y}),\\
\label{24.6.4-2} && {\mathbf{K}}^{\begin{tiny}(T)\end{tiny}} (t,x,y) =\mathbf{K}(t,x,y)+ \mathbf{K}(t,x,\overset{\ast} {y}) -\mathbf{H} (t,x,y)
 \end{eqnarray}
 $\overset{*} {y}$ being the double of $y\in \Omega$ (see, p.$\,$53 of \cite{MS-67}).
 Since the strongly continuous semigroup $(e^{-t\mathcal{L}_{\lambda,\mu}})_{t\ge 0}$ can also be represented as  \begin{eqnarray*} e^{-t\mathcal{L}_{\lambda,\mu} }=\frac{1}{2\pi i} \int_{\mathcal{C}} e^{-t\tau} (\tau I- \mathcal{L}_{\lambda,\mu})^{-1} d\tau,\end{eqnarray*}
where $\mathcal{C}$ is a suitable curve in the complex plane in the positive direction around the spectrum of $\mathcal{L}_{\lambda,\mu}$ (i.e., a contour around the positive real axis). It follows that
 \begin{eqnarray*} {\mathbf{K}} (t,x,y) \!=\! e^{-t\mathcal{L}_{\lambda,\mu}}\delta(x\!-\!y)\! =\! \frac{1}{(2\pi)^n}\!\! \int_{{\Bbb R}^n}\!\! e^{i(x\!-\!y)\cdot\xi} \bigg(\!\frac{1}{2\pi i} \! \int_{\mathcal{C}}\!\! e^{-t\tau}\; \iota \big((\tau {I}\! -\!\mathcal{L}_{\lambda,\mu})^{\!-1}\big) d\tau\!\bigg) d\xi, \;\; \forall\, t>0,\,  x,y\in \mathcal{M},\end{eqnarray*}
where  $\iota(A)$  denotes the full symbol of a pseudodifferential operator $A$.

 We claim that \begin{eqnarray}\; \label{2022.11.11-1}\;\;\, \frac{1}{2\pi i}\! \int_{\mathcal{C}} (\tau I\!-\! \mathcal{L}_{\lambda,\mu} )^{-1} e^{-t\tau}\,\delta(x\!-\!y)\,d\tau  = \frac{1}{2\pi i} \int_{\mathcal{C}} \Big( \int_{\mathbb{R}^n} e^{i(x\!-\!y)\cdot \xi} \sum_{j\le -2} \mathbf{q}_j ( x,\xi, \tau) \,d\xi \Big) e^{-t \tau} d\tau.   \;\;\quad \end{eqnarray}
 In fact, for any smooth vector-valued function $\boldsymbol{\phi}$ with compact support we have \begin{eqnarray*} \big(e^{-t\mathcal{L}_{\lambda,\mu}} \boldsymbol{\phi} \big)(x) \!\!&\!\!=\!\!&\! \! \Big( \frac{1}{2\pi i} \int_{\mathcal{C}} e^{-t\tau} (\tau I- \mathcal{L}_{\lambda,\mu} )^{-1}  d\tau \Big) \, \boldsymbol{\phi}(x)\\
 \!&=\!& \frac{1}{2\pi i} \int_{\mathcal{C}} e^{-t\tau} \Big( \int_{\mathbb{R}^n} e^{i x\cdot \xi} \sum_{j\le -2} \mathbf{q}_j (x, \xi, \tau) \hat{\boldsymbol{\phi}} (\xi) \,d\xi \Big) d\tau.\end{eqnarray*}
On the one hand, from the left-hand side of (\ref{2022.11.11-1}), we get 
\begin{eqnarray} \label{2022.11.11-2}
&& \int_{\mathbb{R}^n}\Big[ \Big( \frac{1}{2\pi i} \int_{\mathcal{C}} (\tau I- \mathcal{L}_{\lambda,\mu} )^{-1}  e^{-t\tau} d\tau\Big) (\delta (x-y))\Big]\boldsymbol{\phi}(y)dy \\
&&   \;\quad\quad \quad =  \Big( \frac{1}{2\pi i} \int_{\mathcal{C}} (\tau- \mathcal{L}_{\lambda,\mu})^{-1} e^{-t\tau} d\tau\Big) \boldsymbol{\phi}(x)
= e^{-t\mathcal{L}_{\lambda,\mu} } \boldsymbol{\phi} (x).\nonumber\end{eqnarray}
On the other hand, from the right-hand side of (\ref{2022.11.11-1}) we obtain
\begin{eqnarray} \label{2022.11.11-3} && \int_{\mathbb{R}^n} \Big[ \frac{1}{2\pi i} \int_C \Big(\int_{\mathbb{R}^n} e^{i(x-y)\cdot \xi} \sum_{j\le -2} \mathbf{q}_j ( x, \xi,\tau) d\xi\Big) e^{-t\tau} d\tau \Big] \boldsymbol{\phi}(y)dy \\
&& \quad \quad \quad =\frac{1}{2\pi i} \int_C \Big(\int_{\mathbb{R}^n} e^{ix\cdot \xi} \sum_{j\le -2} \mathbf{q}_j ( x, \xi,\tau) d\xi\Big) e^{-t\tau} d\tau \int_{\mathbb{R}^n} e^{-y\cdot \xi} \boldsymbol{\phi}(y) dy \nonumber\\
  && \quad \quad \quad = \frac{1}{2\pi i} \int_C \Big(\int_{\mathbb{R}^n} e^{ix\cdot \xi} \sum_{j\le -2} \mathbf{q}_j (x, \xi,\tau) \hat {\boldsymbol{\phi}} (\xi) d\xi\Big) e^{-t\tau} d\tau
= e^{-t\mathcal{L}_{\lambda,\mu}} \boldsymbol{\phi} (x).\qquad \; \nonumber\end{eqnarray}
 Thus, the desired identity (\ref{2022.11.11-1}) is asserted by (\ref{2022.11.11-2}) and (\ref{2022.11.11-3}).

  In particular, for every $t>0$ and $x\in \Omega$, \begin{eqnarray}\label{2020.7.5-1}&&  {\mathbf{K}} (t,x,x) = e^{-t\mathcal{L}_{\lambda,\mu}}\delta(x-x) = \frac{1}{(2\pi)^n} \int_{{\Bbb R}^n}\bigg(\frac{1}{2\pi i} \int_{\mathcal{C}} e^{-t\tau}\; \iota \big((\tau {I} -\mathcal{L}_{\lambda,\mu})^{-1}\big) d\tau\bigg) d\xi\\
    && \qquad\,\; \qquad = \frac{1}{(2\pi)^n} \int_{{\Bbb R}^n}\bigg(\frac{1}{2\pi i} \int_{\mathcal{C}} e^{-t\tau}\; \iota \big((\tau {I} -{L}_{\lambda,\mu})^{-1}\big) d\tau\bigg) d\xi\nonumber\\
    && \qquad\,\;\qquad = \frac{1}{(2\pi)^n} \int_{{\mathbb{R}}^n} \Big( \frac{1}{2\pi i} \int_{\mathcal{C}} e^{-t\tau} \sum_{l\ge 0} q_{-2-l} (x, \xi, \tau) \,d\tau \Big) d\xi,\nonumber\\
\label{2020.7.5-2} &&  {\mathbf{K}} (t,x,\overset{*}{x}) = e^{-t\mathcal{L}_{\lambda,\mu}}\delta(x-\overset{*}{x}) = \frac{1}{(2\pi)^n} \int_{{\Bbb R}^n} e^{i(x-\overset{*}{x})\cdot\xi} \bigg(\frac{1}{2\pi i} \int_{\mathcal{C}} e^{-t\tau}\; \iota \big((\tau {I} -\mathcal{L}_{\lambda,\mu})^{-1}\big) d\tau\bigg) d\xi \qquad\; \\
&& \qquad \qquad\;\;\; = \frac{1}{(2\pi)^n} \int_{{\Bbb R}^n} e^{i(x-\overset{*}{x})\cdot\xi} \bigg(\frac{1}{2\pi i} \int_{\mathcal{C}} e^{-t\tau}\; \iota \big((\tau {I} -{\mathcal{L}_{\lambda,\mu}})^{-1}\big) d\tau\bigg) d\xi\nonumber\\
 && \qquad\;\;\;\qquad = \frac{1}{(2\pi)^n} \int_{{\mathbb{R}}^n} e^{i(x-\overset{*}{x})\cdot \xi} \Big( \frac{1}{2\pi i} \int_{\mathcal{C}} e^{-t\tau} \sum_{l\ge 0} q_{-2-l} (x, \xi, \tau) \,d\tau \Big) d\xi,\nonumber
 \end{eqnarray}
 where $\sum_{l\ge 0} {\mathbf{q}}_{-2-l} (x,\xi,\tau) $ is the full symbol of $(\tau I-L_{\lambda,\mu})^{-1}$.

 Firstly, from the discussion on p.$\,$10182 of \cite{Liu2}, we know that
   \begin{eqnarray}\label{2022.11.9-1}
   && \;\;\;\;{\mathbf{q}}_{-2} (x,\xi,\tau) =\frac{1}{\tau- \mu \sum\limits_{l,m=1}^n g^{lm}\xi_l \xi_m }\,{\mathbf{I}}_n\\
 && \qquad \;\;\;\; \; \; +\frac{\mu+\lambda}{  \big(\tau- \mu \sum\limits_{l,m=1}^n g^{lm}\xi_l \xi_m \big)\big(\tau- (2\mu+\lambda)  \sum\limits_{l,m=1}^n g^{lm}\xi_l \xi_m \big)}\begin{bmatrix} \sum\limits_{r=1}^n g^{1r} \xi_r \xi_1 &\cdots &  \sum\limits_{r=1}^n g^{1r} \xi_r \xi_n\\
\vdots& {} &\vdots \\
 \sum\limits_{r=1}^n g^{nr} \xi_r \xi_1 &\cdots &  \sum\limits_{r=1}^n g^{nr} \xi_r \xi_n\end{bmatrix}\nonumber\end{eqnarray}
 and   \begin{eqnarray}\label{2020.6.6-1}  && \mbox{Tr}\, \big({\mathbf{q}}_{-2} (x,\xi,\tau)\big)=
   \frac{n}{\big(\tau \!-\!\mu \!\sum_{l,m=1}^n \!g^{lm} \xi_l\xi_m\big)} \\
 && \qquad \qquad \qquad \qquad \;\;  + \frac{(\mu+\lambda)\sum_{l,m=1}^n g^{lm} \xi_l\xi_m}{
  \big(\!\tau\! -\!\mu\! \sum_{l,m=1}^n \!g^{lm} \xi_l\xi_m\!\big)
 \big(\!\tau \!-\!(2\mu\!+\!\lambda)  \!\sum_{l,m=1}^n \!g^{lm} \xi_l\xi_m\!\big)}
.\nonumber \end{eqnarray}
       For each $x\in \Omega$, we use a geodesic normal coordinate system centered at this $x$. It follows from \S11 of Chap.1 in \cite{Ta1}
        that in such a coordinate system, $g_{jk}(x)=\delta_{jk}$ and $\Gamma_{jk}^l(x)=0$. Then (\ref{2020.6.6-1}) reduces to
        \begin{eqnarray}\label{2020.7.6-3} \quad \quad \,\;\mbox{Tr} \,\big({\mathbf{q}}_{-2} (x, \xi,\tau)\big)=
     \frac{n}{(\tau -\mu |\xi|^2)}+
      \frac{(\mu+\lambda)|\xi|^2}{(\tau -\mu |\xi|^2) (\tau-(2\mu+\lambda)|\xi|^2)},\end{eqnarray}
  where $|\xi|=\sqrt{\sum_{k=1}^n \xi^2_k}$ for any $\xi\in {\mathbb{R}}^n$.
   By applying the residue theorem (see, for example, Chap.$\,$4, \S5 in \cite{Ahl}) we get
 \begin{eqnarray} \label{3.10} \frac{1}{2\pi i}\! \int_{\mathcal{C}}\! e^{-t\tau} \bigg(\frac{n}{(\tau -\mu |\xi|^2)}\!+\!
      \frac{(\mu+\lambda)|\xi|^2}{(\tau -\mu |\xi|^2) (\tau-(2\mu\!+\!\lambda)|\xi|^2)} \bigg) d\tau\!=\! (n\!-\!1) e^{-t\mu|\xi|^2} \!+\!  e^{-t(2\mu+\lambda)|\xi|^2}.\end{eqnarray}
       It follows that
      \begin{eqnarray}\label{2020.7.10-1} \frac{1}{(2\pi)^n}\!\!\!&\!\!\!&\!\!\!\! \!\!\! \! \int_{{\mathbb{R}}^n}\Big( \frac{1}{2\pi i} \int_{\mathcal{C}} e^{-t\tau}\, \mbox{Tr}\,({\mathbf{q}}_{-2} (x, \xi,\tau) ) d\tau \Big) d\xi \\
      \!\!\! &=\!\!\!&
    \frac{1}{(2\pi)^n} \int_{{\Bbb R}^n} \bigg((n-1) e^{-t\mu|\xi|^2} +  e^{-t(2\mu+\lambda)|\xi|^2}\bigg)
         d\xi \nonumber\\
               \!\!\! &=\!\!\!&\frac{n-1}{(4\pi \mu t)^{n/2}} +  \frac{1}{(4\pi (2\mu+\lambda) t)^{n/2}},
                  \nonumber\end{eqnarray}
                  and hence
        \begin{eqnarray}\label{2020.7.12-3}  && \int_{\Omega}\! \left\{\! \frac{1}{(2\pi)^n}\! \int_{{\mathbb{R}}^n}\!\Big( \frac{1}{2\pi i} \int_{\mathcal{C}} e^{-t\tau}\, \mbox{Tr}\,({\mathbf{q}}_{-2} (x, \xi,\tau) ) d\tau\! \Big) d\xi\!\right\}\! dV\\
        && \qquad \;\;=\Big(\frac{n-1}{(4\pi \mu t)^{n/2}} \! +\!  \frac{1}{(4\pi (2\mu\!+\!\lambda) t)^{n/2}}\!\Big){\mbox{Vol}}(\Omega).\nonumber\end{eqnarray}

In the above discussion, if we replace $x\in \Omega$ by $\overset{*}{x} \in \Omega^*$, then (\ref{2022.11.9-1}) will become
\begin{eqnarray*}\label{2022.11.9-2}
  \!\! \!\!\!&& \!\!\!\!\!{\mathbf{q}}_{-2} (\overset{*}{x},\xi,\tau) =\frac{1}{\tau- \mu \sum\limits_{l,m=1}^n (g^{lm}(\overset{*}{x}))\xi_l \xi_m }\,{\mathbf{I}}_n +\frac{\mu\!+\!\lambda}{  \big(\tau\!-\! \mu \sum\limits_{l,m=1}^n\! (g^{lm}(\overset{*}{x}))\xi_l \xi_m \big)\big(\tau\!- \!(2\mu\!+\!\lambda)  \sum\limits_{l,m=1}^n \!(g^{lm}(\overset{*}{x}))\xi_l \xi_m \big)}\\
 && \;\;\;\;\times \begin{bmatrix}\! \sum\limits_{r=1}^n \!(g^{1r}(\overset{*}{x})) \xi_r \xi_1 \!&\! \cdots \!&\! \sum\limits_{r=1}^n( g^{1r}(\overset{*}{x}))\xi_r \xi_{n-1}\!& \!\sum\limits_{r=1}^n \!(-g^{1r}(\overset{*}{x})) \xi_r \xi_n\\
\vdots\!& {}\! &\vdots &\vdots \\
\sum\limits_{r=1}^n (g^{n\!-\!1,r}(\overset{*}{x})) \xi_r \xi_1 \!\!&\!\cdots \!&\! \sum\limits_{r=1}^n (g^{n-1,r}(\overset{*}{x})) \xi_r \xi_{n\!-\!1}\!& \! \sum\limits_{r=1}^n (-g^{n\!-\!1,r}(\overset{*}{x})) \xi_r \xi_n
\\
 \sum\limits_{r=1}^n\! (-g^{nr}(\overset{*}{x})) \xi_r \xi_1 \!&\!\cdots \!& \!\sum\limits_{r=1}^n \!(-g^{nr}(\overset{*}{x})) \xi_r \xi_{n-1}\!& \! \sum\limits_{r=1}^n (g^{nr} (\overset{*}{x}))\xi_r \xi_n\end{bmatrix}\nonumber\end{eqnarray*}
 and   \begin{eqnarray*}\label{2022.11.9-4}  \mbox{Tr}\, \big({\mathbf{q}}_{-2} (\overset{*}{x},\xi,\tau)\big)=
   \frac{n}{\big(\tau \!-\!\mu \!\sum_{l,m=1}^n \!(g^{lm}(\overset{*}{x})) \xi_l\xi_m\big)} \!+\! \frac{(\mu+\lambda)\sum_{l,m=1}^n (g^{lm}(\overset{*}{x})) \xi_l\xi_m}{
  \big(\!\tau\! -\!\mu\! \sum_{l,m=1}^n \!(g^{lm}(\overset{*}{x})) \xi_l\xi_m\!\big)
 \big(\!\tau \!-\!(2\mu\!+\!\lambda)  \!\sum_{l,m=1}^n \!(g^{lm}(\overset{*}{x})) \xi_l\xi_m\!\big)}
.\end{eqnarray*}
This implies that all expressions (\ref{2020.6.6-1})--(\ref{2020.7.12-3}) of the above  trace symbols  have the same form either in $\Omega$ or in $\Omega^*$.

      For given (small) $\epsilon>0$ , denote by $U_\epsilon(\partial \Omega)=\{z\in {\mathcal{M}}\big| \mbox{dist}\, (z, \partial \Omega)<\epsilon\}$ the $\epsilon$-neighborhood of $\partial \Omega$ in $\mathcal{M}$.
             When $x\in \Omega\setminus U_\epsilon (\partial \Omega)$,
we see by taking geodesic normal coordinate system at $x$ that (\ref{2020.7.6-3}) still holds at this $x$. According to (\ref{3.10}) we have
  that
     \begin{eqnarray} \label{3.11}  \mbox{Tr}\,({\mathbf{q}}_{-2}(t, x, \overset{*}{x}))\!\!\! &=\!\!\!&
     \frac{1}{(2\pi)^n} \int_{{\Bbb R}^n} e^{i(x-\overset{*}{x})\cdot\xi}
    \bigg((n-1) e^{-t\mu|\xi|^2} +  e^{-t(2\mu+\lambda)|\xi|^2}\bigg)
         d\xi \nonumber\\
               \!\!\! &=\!\!\!&\frac{n-1}{(4\pi \mu t)^{n/2}}e^{-\frac{|x-\overset{*}{x}|^2}{4t\mu}} +  \frac{1}{(4\pi (2\mu+\lambda) t)^{n/2}}e^{-\frac{|x-\overset{*}{x}|^2}{4t(2\mu+\lambda)}} \quad \, \mbox{for any}\;\, x\in \Omega \setminus U_\epsilon (\partial \Omega),
                  \nonumber\end{eqnarray}
     which exponentially tends to zero as $t\to 0^+$ because $|x-\overset{*}{x}|\ge \epsilon$.  Hence
     \begin{eqnarray} \label{2020.7.6-4}  \int_{\Omega \setminus U_\epsilon (\partial \Omega)} \left(\mbox{Tr}\,({\mathbf{q}}_{-2}(t, x, \overset{*}{x}))\right)\,dV\!\!\! &=\!\!\!&
     O(t^{1-\frac{n}{2}}) \quad \;\mbox{as}\;\; t\to 0^+.
                 \end{eqnarray}

                Secondly, for $l\ge 1$, it can be verified that
                  $\mbox{Tr}\, ({\mathbf{q}}_{-2-l} (x, \xi, \tau))$ is a sum of finitely many terms, each of which  has the following form:
                  $$\frac{r_k(x, \xi)}{(\tau-\mu \sum_{l,m=1}^n g^{lm} \xi_l \xi_m )^s (\tau -(2\mu +\lambda)\sum_{l,m=1}^n g^{lm} \xi_l\xi_m )^j },$$
where $k-2s-2j=-2-l$, and $r_k(x,\xi)$ is the symbol independent of $\tau$ and homogeneous of degree $k$.
 Again we take the geodesic normal coordinate systems center at $x$ (i.e., $g_{jk}(x)=\delta_{jk}$ and $\Gamma_{jk}^l(x)=0$), by applying residue theorem we see  that, for $l\ge 1$ , \begin{eqnarray*}\label{2020.7.12-1} && \frac{1}{(2\pi)^n} \int_{{\mathbb{R}}^n}\Big( \frac{1}{2\pi i} \int_{\mathcal{C}} e^{-t\tau}\, \mbox{Tr}\,({\mathbf{q}}_{-2-l} (x, \xi,\tau) ) d\tau \Big) d\xi=O(t^{l-\frac{n}{2}}) \;\;\mbox{as}\;\, t\to 0^+ \;\;\, \mbox{uniformly for} \;\, x\in \Omega,
 \end{eqnarray*} and
\begin{eqnarray} \!\! \label{2020.7.14-1}  \\
 \frac{1}{(2\pi)^n} \!\int_{{\mathbb{R}}^n}\! e^{i(x-\overset{*}{x})\cdot \xi} \Big(\! \frac{1}{2\pi i} \!\int_{\mathcal{C}}\! e^{-t\tau}\, \mbox{Tr}\,({\mathbf{q}}_{-2-l} (x, \xi,\tau) ) d\tau \!\Big) d\xi\!=\!O(t^{l-\frac{n}{2}}\!) \;\,\mbox{as}\;\, t\to 0^+ \,\; \mbox{uniformly for} \;\, x\in \Omega.\nonumber
      \end{eqnarray}
Therefore
\begin{eqnarray}\label{2020.7.12-10}  \!\!  \int_{\Omega} \bigg\{\frac{1}{(2\pi)^n} \int_{{\mathbb{R}}^n}\Big( \frac{1}{2\pi i} \int_{\mathcal{C}} e^{-t\tau}\, \sum_{l\ge 1}\mbox{Tr}\,({\mathbf{q}}_{-2-l} (x, \xi,\tau) ) d\tau \Big) d\xi\bigg\} dV =O(t^{1-\frac{n}{2}}) \;\;\mbox{as}\;\, t\to 0^+,
 \end{eqnarray} and
\begin{eqnarray}\label{2020.7.13-5}  \!\!\! \!\!
\int_{\Omega} \!\bigg\{\! \frac{1}{(2\pi)^n} \!\int_{{\mathbb{R}}^n}\!\!\! e^{i(x-\overset{*}{x})\cdot \xi} \Big( \frac{1}{2\pi i} \!\int_{\mathcal{C}} e^{-t\tau} \sum_{l\ge 1}\mbox{Tr}\,({\mathbf{q}}_{-2-l} (x, \xi,\tau) ) d\tau \!\Big) d\xi\bigg\}dV\!=\!O(t^{1-\frac{n}{2}})\;\;\mbox{as}\,\, t\to 0^+. 
      \end{eqnarray}
Combining (\ref{2020.7.5-1}), (\ref{2020.7.12-3}) and (\ref{2020.7.12-10}), we have
       \begin{eqnarray} \label{3.11}  \!\! \!\! \int_{\Omega}\mbox{Tr}\,({\mathbf{K}}(t, x, x)) \, dV =    \bigg[  \frac{n-1}{(4\pi \mu t)^{n/2}} +  \frac{1}{(4\pi (2\mu+\lambda) t)^{n/2}}\bigg]{\mbox{Vol}}(\Omega)+O(t^{1-\frac{n}{2}})\;\;\mbox{as}\,\; t\to 0^+. \end{eqnarray}

     Now,  we will consider the case of $\int_{\Omega\cap U_\epsilon(\partial \Omega)} \left\{\frac{1}{(2\pi)^n} \int_{{\mathbb{R}}^n} e^{i(x-\overset{*}{x})\cdot \xi} \Big(  \frac{1}{2\pi i} \int_{\mathcal{C}} e^{-t\tau}\,\mbox{Tr}\,({\mathbf{q}}_{-2} (x, \xi,\tau) ) d\tau \Big) d\xi\right\} dV$.
We pick a self-double patch $W$ of $\mathcal{M}$ (such that $W\subset U_\epsilon (\partial \Omega)$) covering a patch $W\cap \partial \Omega$ of $\partial \Omega$  endowed (see the diagram on p.$\,$54 of \cite{MS-67}) with local coordinates $x$ such that
\vskip 0.26 true cm

\begin{figure}[h]
\centering
\begin{tikzpicture}[scale=1,line width=0.8]

\draw (4,0) arc (70: 90: 20);
\draw (3.5,1.5) arc (70: 85: 20);
\draw (2.7,-1) arc (70: 86: 16);
\draw (3.5,1.52) -- (2.7,-1);
\draw (-1.57,2.65) -- (-1.63,-0.1);

\draw[fill=black] (1,1.8) circle (1pt);
\draw[fill=black] (0.82,0.87) circle (1pt);
\draw[fill=black] (0.69,0) circle (1pt);

\draw[->,>=stealth] (1,1.8) .. controls (0.8,0.87) and (0.69,0) .. (0.5,-1.5);

\node at (-0.2,2) {$\Omega^{*}$};
\node at (-0.5,0.5) {$\Omega$};
\node at (0.5,-1.8) {$x_n>0$};
\node at (1,0) {$x$};
\node at (1.4,1.8) {$x^*$};
\node at (-2.4,1.5) {$\partial \Omega$};
\node at (-2,-2){};

\end{tikzpicture}
\end{figure}

\vskip 0.1 true cm 

   $\epsilon>x_n>0$ in $W\cap \Omega$; $\,x_n=0$ on $W\cap \partial \Omega$;
  $\; x_n (\overset{*}{x})=-x_n(x)$; and the positive $x_n$-direction is perpendicular to $\partial \Omega$. This has the effect that (\ref{2021.2.6-3})--(\ref{2021.2.6-4})  and  \begin{eqnarray}
   \label{2021.2.6-5}  \sqrt{|g|/g_{nn}} \; dx_1\cdots dx_{n-1}\!\!\!&\!=\!&\!\!\! \mbox{the element of (Riemannian) surface area on} \,\, \partial \Omega.\end{eqnarray}
    We choose coordinates $x'=(x_1,\cdots, x_{n-1})$ on an open set in $\partial \Omega$ and then
 coordinates $(x', x_{n})$ on a neighborhood in $\bar\Omega$ such that
$x_{n}=0$ on $\partial \Omega$ and $|\nabla x_{n}|=1$ near $\partial \Omega$ while $x_{n}>0$ on $\Omega$ and such that $x'$ is constant
on each geodesic segment in $\bar\Omega$ normal to $\partial \Omega$.
   Then the metric tensor on $\bar \Omega$ has
the form (see \cite{LU} or p.$\,$532 of \cite{Ta2})
\begin{eqnarray} \label{a-1}  \big(g_{jk} (x',x_{n}) \big)_{n\times n} =\begin{pmatrix} ( g_{jk} (x',x_{n}))_{(n-1)\times (n-1)}& 0\\
      0& 1 \end{pmatrix}. \end{eqnarray}
       Furthermore, we can take a geodesic normal coordinate system for $(\partial \Omega, g)$ centered at $x_0=0$, with respect to $e_1, \cdots, e_{n-1}$, where  $e_1, \cdots, e_{n-1}$ are the principal curvature vectors. As Riemann showed, one has (see p.$\,$555 of \cite{Ta2})
            \begin{eqnarray} \label{7/14/1}& g_{jk}(x_0)= \delta_{jk}, \; \; \frac{\partial g_{jk}}{\partial x_l}(x_0)
 =0  \;\;  \mbox{for all} \;\; 1\le j,k,l \le n-1,\\
 & -\frac{1}{2}\frac{\partial g_{jk}}{\partial x_n} (x_0) =\kappa_k\delta_{jk}  \;\;  \mbox{for all} \;\; 1\le j,k \le n-1,\nonumber
 \end{eqnarray}
 where $\kappa_1\cdots, \kappa_{n-1}$ are the principal curvatures of $\partial \Omega$ at point $x_0=0$.
   Due to the special geometric normal coordinate system and (\ref{7/14/1})--(\ref{a-1}), we see that for any $x\in \{z\in \Omega\big|\mbox{dist}(z, \partial \Omega)< \epsilon\}$,
   \begin{eqnarray} \label{18-4-5-1} x-\overset{\ast} {x}=(0,\cdots, 0, x_n-(-x_n))=(0,\cdots, 0,2x_n).\end{eqnarray}
              By (3.17) of \cite{Liu2}, (\ref{7/14/1}), (\ref{3.10}) and (\ref{18-4-5-1}), we find that
      \begin{eqnarray*} \label{18-4-1-1} &&\int_{W\cap \Omega} \bigg\{\frac{1}{(2\pi)^n} \int_{{\Bbb R}^{n}} e^{i\langle x-\overset{*}{x}, \xi\rangle} \Big( \frac{1}{2\pi i} \int_{\mathcal{C}} e^{-t\tau} \,\mbox{Tr}\,\big({\mathbf{q}}_{-2} (x, \xi, \tau)\big) d\tau \Big)  \, d\xi\bigg\} dV \\
           &&   =\! \int_0^\epsilon dx_n\! \int_{W\cap \partial \Omega}\!
 \frac{dx'}{(2\pi)^n}\!\int_{{\Bbb R}^{n}} \!e^{i\langle 0, \xi'\rangle + i2x_n \xi_n} \bigg[\!\frac{1}{2\pi i}\! \int_{\mathcal{C}}\! e^{-t\tau} \!\Big(  \frac{n}{(\tau \!-\!\mu |\xi|^2)}\!+\!
      \frac{(\mu\!+\!\lambda)|\xi|^2}{(\tau \!-\!\mu |\xi|^2) (\tau\!-\!(2\mu\!+\!\lambda)|\xi|^2)}\Big) d\tau\!\bigg] d\xi\\
  && \,=  \int_0^\epsilon dx_n \int_{W\cap \partial \Omega}
 \frac{dx'}{(2\pi)^n}\int_{{\Bbb R}^{n}} e^{i2x_n \xi_n} \bigg((n-1) e^{-t\mu|\xi|^2} +  e^{-t(2\mu+\lambda)|\xi|^2} d\tau\bigg) d\xi\\
 &&\,= \int_0^\epsilon dx_n \int_{W\cap \partial \Omega}  \frac{dx'}{(2\pi)^n}\int_{-\infty}^\infty e^{2ix_n \xi_n}\bigg[ \int_{{\Bbb R}^{n-1}} \bigg( (n-1) e^{-t\mu(|\xi'|^2+\xi_n^2)} +  e^{-t(2\mu+\lambda)(|\xi'|^2+\xi_n^2)} \bigg) d\xi' \bigg]d\xi_n  \nonumber\\
    &&\,= \int_0^\epsilon dx_n \int_{W\cap \partial \Omega}  \frac{1}{(2\pi)^n}\bigg[\int_{-\infty}^\infty e^{2ix_n \xi_n}e^{-t\mu \xi_n^2} \bigg( \int_{{\Bbb R}^{n-1}}  (n-1) e^{-t\mu \sum_{j=1}^{n-1}\xi_j^2}  d\xi' \bigg)d\xi_n\bigg] dx'\nonumber\\
  &&\quad\, + \int_0^\epsilon dx_n \int_{W\cap \partial \Omega}  \frac{1}{(2\pi)^n}\bigg[\int_{-\infty}^\infty e^{2ix_n \xi_n}e^{-t(2\mu+\lambda) \xi_n^2} \bigg( \int_{{\Bbb R}^{n-1}} e^{-t(2\mu+\lambda) \sum_{j=1}^{n-1}\xi_j^2} d\xi' \bigg)d\xi_n\bigg]dx'\nonumber, \end{eqnarray*}
 where  $\xi=(\xi', \xi_n)\in {\Bbb R}^n$, $\xi'=(\xi_1, \cdots, \xi_{n-1})$.
   A direct calculation shows that
   \begin{eqnarray*}
 &&  \frac{1}{(2\pi)^n}\bigg[\int_{-\infty}^\infty e^{2ix_n \xi_n}e^{-t\mu \xi_n^2} \bigg( \int_{{\Bbb R}^{n-1}}  (n-1) e^{-t\mu \sum_{j=1}^{n-1}\xi_j^2}  d\xi' \bigg)d\xi_n\bigg]=\frac{n-1}{(4\pi \mu  t)^{n/2}}\, e^{-\frac{(2x_n)^2}{ 4\mu t} },\\
  &&    \frac{1}{(2\pi)^n}\bigg[\int_{-\infty}^\infty e^{2ix_n \xi_n}e^{-t(2\mu+\lambda) \xi_n^2} \bigg( \int_{{\Bbb R}^{n-1}} e^{-t(2\mu+\lambda) \sum_{j=1}^{n-1}\xi_j^2} d\xi' \bigg)d\xi_n= \frac{1}{(4\pi (2 \mu+\lambda)  t)^{n/2}}\, e^{-\frac{(2x_n)^2}{ 4(2\mu+\lambda) t} }.\end{eqnarray*}
   Hence
    \begin{eqnarray} \label{18-4-1-10.} \\
 && \int_{W\cap \Omega} \bigg\{\frac{1}{(2\pi)^n} \int_{{\Bbb R}^{n}} e^{i\langle x-\overset{*}{x}, \xi\rangle} \Big( \frac{1}{2\pi i} \int_{\mathcal{C}} e^{-t\tau} \,\mbox{Tr}\,\big({\mathbf{q}}_{-2} (x, \xi, \tau)\big) d\tau \Big)  \, d\xi\bigg\} dV\nonumber \\
  &&\;\;\quad = \int_0^\epsilon dx_n \int_{W\cap \partial \Omega} \left[\frac{n-1}{(4\pi \mu  t)^{n/2}}\, e^{-\frac{(2x_n)^2}{ 4\mu t} }+\frac{1}{(4\pi (2 \mu+\lambda)  t)^{n/2}}\, e^{-\frac{(2x_n)^2}{ 4(2\mu+\lambda) t} }\right]dx'
    \nonumber\\
    &&\;\; \quad= \int_0^\infty dx_n \int_{W\cap \partial \Omega} \left[ \frac{n-1}{(4\pi \mu  t)^{n/2}}\, e^{-\frac{(2x_n)^2}{ 4\mu t} }+\frac{1}{(4\pi (2\mu+\lambda)  t)^{n/2}}\, e^{-\frac{(2x_n)^2}{ 4(2\mu+\lambda) t} }\right]dx'\nonumber\\
    && \;\; \quad\;\quad - \int_\epsilon^\infty dx_n \int_{W\cap \partial \Omega} \left[\frac{n-1}{(4\pi \mu  t)^{n/2}}\, e^{-\frac{(2x_n)^2}{ 4\mu t} }+ \frac{1}{(4\pi (2\mu +\lambda) t)^{n/2}}\, e^{-\frac{(2x_n)^2}{ 4(2\mu+\lambda) t} }\right]dx'\nonumber\\
   &&\;\;\quad= \frac{n-1}{4} \cdot \frac{\mbox{Vol}(W\cap \partial \Omega)}{(4\pi \mu t)^{(n-1)/2}} + \frac{1}{4} \cdot \frac{\mbox{Vol}(W\cap \partial \Omega)}{(4\pi (2\mu +\lambda) t)^{(n-1)/2}}\nonumber\\
   &&\;\;\quad \quad \, -   \int_{W\cap \partial \Omega} \left\{\int_\epsilon^\infty\bigg[\frac{n-1}{(4\pi \mu  t)^{n/2}}\, e^{-\frac{(2x_n)^2}{ 4\mu t} }+ \frac{1}{(4\pi (2\mu +\lambda) t)^{n/2}}\, e^{-\frac{(2x_n)^2}{ 4(2\mu+\lambda) t}} \bigg]dx_n\right\}dx'. \nonumber\end{eqnarray}
      It is easy to verify that for any  fixed $\epsilon>0$,
   \begin{eqnarray} \label{18-4-1-11.}\begin{aligned} && \int_\epsilon^\infty \frac{1}{(4\pi \lambda t)^{\frac{n}{2}}} e^{-\frac{(2x_n)^2}{4\mu t}}dx_n = O(t^{1-n/2})\quad \; \, \mbox{as} \;\, t\to 0^+,\qquad \;\; \quad \\
  && \int_\epsilon^\infty \frac{1}{(4\pi (2\mu+\lambda) t)^{\frac{n}{2}}} e^{-\frac{(2x_n)^2}{4(2\mu+\lambda) t}}dx_n = O(t^{1-n/2})
   \quad \; \, \mbox{as} \;\, t\to 0^+.\end{aligned}\end{eqnarray}
   From (\ref{18-4-1-10.}) and (\ref{18-4-1-11.}), we get that
     \begin{eqnarray}\label{18-4-1-6}\!\!\!\! \!\!\!\!\! \!\!\!\!\! && \int_{W\cap \Omega} \bigg\{\frac{1}{(2\pi)^n} \int_{{\Bbb R}^{n}} e^{i\langle x-\overset{*}{x}, \xi\rangle} \Big( \frac{1}{2\pi i} \int_{\mathcal{C}} e^{-t\tau} \,\mbox{Tr}\,\big({\mathbf{q}}_{-2} (x, \xi, \tau)\big) d\tau \Big)  \, d\xi\bigg\} dV\\
\!\!\!\!\!\!\!\!\! \!\!\!\! && \,\quad \;\,\quad     =
     \frac{n\!-\!1}{4} \cdot \frac{\mbox{Vol}(W\cap \partial \Omega)}{(4\pi \mu t)^{(n-1)/2}} + \frac{1}{4} \cdot \frac{\mbox{Vol}(W\cap \partial \Omega)}{(4\pi (2\mu +\lambda) t)^{(n-1)/2}}
 + O(t^{1-n/2}) \quad \; \mbox{as} \, \; t\to 0^+.\nonumber\end{eqnarray}
  For any $x\in \Omega\cap U_\epsilon (\partial \Omega)$, we have
  \begin{eqnarray} \label{18-4-3-1}\;\;\,\;\, \;\;\;\;\; \mbox{Tr}\,(K(t, x, \overset{*}{x}))
 \!\!\!&\!=&\!\!\!\! \frac{1}{(2\pi)^n}\int_{{\Bbb R}^n} \!  e^{i\langle x-\overset{*}{x}, \xi\rangle}\big( \frac{1}{2\pi i} \!\int_{\mathcal{C}}\! e^{-t\tau}\, \mbox{Tr}\,\big({\mathbf{q}}_{-2} (x, \xi, \tau)\big) d\tau\big) d\xi\\
  \!\!\!&\!&\!\!\!\!+ \frac{1}{(2\pi)^n}\int_{{\Bbb R}^n}   e^{i\langle x-\overset{*}{x}, \xi\rangle}\big(\sum_{l\ge 1} \frac{1}{2\pi i} \int_{\mathcal{C}} e^{-t\tau} \mbox{Tr}\,\big({\mathbf{q}}_{-2-l} (x, \xi, \tau)\big) d\tau\big) d\xi\nonumber\\
  \!\!\!&\!=&\!\!\!\! \!\frac{1}{(2\pi)^n}\!\!\int_{{\Bbb R}^n} \! \! e^{i\langle x-\overset{*}{x}, \xi\rangle}\big( \frac{1}{2\pi i} \!\int_{\mathcal{C}} e^{-t\tau} \mbox{Tr}\,\big({\mathbf{q}}_{-2} (x, \xi, \tau)\big) d\tau\big) d\xi\!+\!O(t^{1+\frac{n}{2}}) \;\; \mbox{as}\;\, t\to 0^+,\nonumber\end{eqnarray}  where the second equality used (\ref{2020.7.14-1}).
  Combining (\ref{18-4-1-6}) and (\ref{18-4-3-1}), we have
     \begin{eqnarray} \label{3..15}
     &&\int_{W\cap \Omega}\mbox{Tr}\big( \mathbf{K}(t, x, \overset{*}{x})\big)dx=  \frac{n-1}{4} \cdot \frac{\mbox{Vol}(W\cap \partial \Omega)}{(4\pi \mu t)^{(n-1)/2}} \\
  && \quad\,\quad \;\,\quad \;+ \frac{1}{4} \cdot \frac{\mbox{Vol}(W\cap \partial \Omega)}{(4\pi (2\mu +\lambda) t)^{(n-1)/2}}
 + O(t^{1-n/2}) \quad \; \mbox{as} \, \; t\to 0^+.\nonumber\end{eqnarray}
     It follows from  (\ref{c4-23}), (\ref{2020.7.5-2}), (\ref{2020.7.6-4}), (\ref{2020.7.13-5}), (\ref{3.11}) and (\ref{3..15}) that
 \begin{eqnarray} \label{a-4-1-3} && \int_{W\cap \Omega}\mbox{Tr}\big({\mathbf{K}}(t, x, x) \big) dx  \mp \int_{W\cap \Omega}\mbox{Tr}\big( {\mathbf{K}}(t, x, \overset{*}{x})\big)dx \\
 &&\qquad\;\; =\bigg[  \frac{n-1}{(4\pi \mu t)^{n/2}}  +  \frac{1}{(4\pi (2\mu+\lambda) t)^{n/2}}\bigg]\mbox{Vol}(W\cap\Omega) \nonumber\\
&&\qquad \;\, \;\,\;  \;\mp \frac{1}{4}\bigg[(n-1) \frac{\mbox{Vol}(W\cap \partial \Omega)}{(4\pi \mu t)^{(n-1)/2}}+ \frac{\mbox{Vol}(W\cap \partial \Omega)}{(4\pi (2\mu +\lambda) t)^{(n-1)/2}}\bigg]\nonumber \\
 &&\quad \qquad \; \; +O(t^{1-n/2})\quad \; \mbox{as} \, \; t\to 0^+.\nonumber\end{eqnarray}
  Note that \begin{eqnarray}  \label {2.12.10-6} \sum_{k=1}^\infty e^{- t \tau_k^{\begin{tiny}(D)\end{tiny}} } = \int_\Omega \mbox{Tr}\, \big( \mathbf{K}^{\begin{tiny}(D)\end{tiny}} (t, x,x) \big) \,dx  
 = \int_\Omega  \mbox{Tr}\, \Big( \mathbf{K} (t,x,x) -\mathbf{K} (t, x, \overset{*}{x})\Big) \,dx    \end{eqnarray} 
and   \begin{eqnarray} \label{23.12.10-7}
\sum_{k=1}^\infty e^{- t \tau_k^{\begin{tiny}(T)\end{tiny}} } 
\!\!\!\!  &&  \!\!\!= \int_\Omega \mbox{Tr}\, \big( \mathbf{K}^{\begin{tiny}(T)\end{tiny}} (t, x,x) \big) \,dx  \\
 && \!\!\! = \int_\Omega  \mbox{Tr} \,\Big( \mathbf{K} (t,x,x) +\mathbf{K} (t, x, \overset{*}{x})\Big) \,dx  - \int_\Omega  \mbox{Tr}\, \Big( \mathbf{H}  (t, x, x) \Big) \,dx. \nonumber \end{eqnarray}

Finally,  by applying Lemma 2.1, we get the third  coefficient $a_2$  in the asymptotic expansions for the Dirichlet and traction boundary conditions \begin{eqnarray} \label{24.5.22-01}&a_2 = 
\frac{1}{6 (4\pi)^{n/2}} \Big[\Big( \frac{1}{ (\lambda+2\mu)^{(n-2)/2}} +\frac{n-7}{ \mu^{(n-2)/2}} + \frac{12\mu}{n} \big( \frac{1}{(\lambda+2\mu)^{n/2}} +\frac{n-1}{\mu^{n/2}} \big) \Big)\! \int_{\Omega} \! R(x)\, dV\\ 
 &\quad  + 2 \,\Big(\frac{1}{ (\lambda+2\mu)^{(n-2)/2}} +\frac{n-7}{ \mu^{(n-2)/2}}\Big) \int_{\partial \Omega}H(x) \,ds\Big]. \nonumber \end{eqnarray} 
 From  (\ref{a-4-1-3})--(\ref{23.12.10-7}),  (\ref{23.12.9-1}) and (\ref{24.5.22-01}),  we get  the conclusion of Theorem 1.4.
   $\square$

  \vskip 0.38  true cm 

   \noindent{\bf  Remark 5.1.} \ {\it  \ (i) \ \ For the two-term heat trace asymptotic formula for the Lam\'{e} operator with Dirichlet or traction boundary conditions, we refer the reader to  \cite{ Liu2}, \cite{Liu-23} or \cite{PiOr},  in which  an  open problem mentioned by Avramidi in \cite{Avr} was affirmatively answered. }

\vskip 0.82 true cm

\section{the Stokes eigenvalues and buckling eigenvalues}
\vskip 0.39 true cm

 The buckling eigenvalue problem with Dirichlet boundary condition  has interpretations in physics, that is, it describes the critical buckling load of a clamped plate subjected to a uniform compressive force around its boundary.  In 1986, Girault and Raviart (see \cite{GiR} or \cite{ChL}) proved that all the Stokes
eigenvalues with Dirichlet boundary condition coincide with all the buckling eigenvalues in the two-dimensional case. 
In this section, we shall also show that in two dimensions, the eigenvalue problem for the  Stokes operator with Cauchy force boundary condition is equivalent to the eigenvalue problem for the buckling plate problem with Cauchy force boundary condition.

\vskip 0.28  true cm  

    \noindent  {\it Proof of Theorem 1.8.} \    As pointed out above,  we shall only prove the case of the  Cauchy force boundary condition, i.e.,  we shall show  that
    \begin{eqnarray} \label{24.4.24-01} \varsigma_{k+1}^{\begin{tiny} (C)\end{tiny}}= \Lambda_k^{\begin{tiny} (C)\end{tiny}}\;\;\,\;\mbox{for all}\;\; k\ge 1.\end{eqnarray}    

 For a vector field $\mathbf{u}=(u^1,u^2)$ we define $\mathbf{u}^{\perp}=(-u^2, u^1)$ and for a scalar $\psi$ we define $\nabla^{\perp} \psi=(-\partial_2 \psi, \partial_1 \psi)$. Clearly,   $(\mathbf{u}^{\perp})^{\perp}=-\mathbf{u}$ and $(\nabla^{\perp})^{\perp}=-\nabla \psi$. 
 By $\omega(\mathbf{u})$ we mean the vorticity (scalar curl) of $\mathbf{u}$:  
 $\;\omega(\mathbf{u})= \partial_1 u^2-\partial_2u^1$. Note that $\omega(\mathbf{u}) =-\mbox{div}\, \mathbf{u}^{\perp}$. In particular,  if  $\mbox{div}\; \mathbf{u}=0$,  then $\nabla^{\perp} \omega(\mathbf{u}) =\Delta \mathbf{u}$.

  Clearly, $(\varsigma_1, \mathbf{u}_1)=(0, c)$ is the first eigenpair of the Lam\'{e} operator with the Cauchy force boundary condition.  
  Let $(\varsigma_k, \mathbf{u}_k)$, ($k\ge 2$),  be a Stokes eigenpair with Cauchy force boundary condition.  Then there is a scalar function $p_k$ in $H^{1}(\Omega)$ such that 
 \begin{eqnarray}\label{24.4.23-10} \left\{ \begin{array}{ll}  -\mu\, \Delta \mathbf{u}_k + \nabla p_k =  \varsigma_k \mathbf{u}_k \;\, &\mbox{in}\;\;\Omega,\\
 \mbox{div}\; \mathbf{u}_k =0 \;\,  &\mbox{in}\;\; \Omega,\\
\Big(2\mu \,(\mbox{Def}\; \mathbf{u}_k)-p_k \mathbf{I}_2\Big) \,\boldsymbol{\nu}=0\;\,  &\mbox{on} \;\;\partial \Omega. \end{array} \right.\end{eqnarray}  
 Here the above pressure function $p_k$ is uniquely determined by letting $\,\inf\limits_{x\in \Omega} p_k(x)=0$, since $p_k$ is uniquely determined up to an additive arbitrary constant.
Note that the Cauchy force boundary condition in eigenvalue problem (\ref{24.4.23-10})  can exactly be expressed as 
\begin{eqnarray} \label{24.4.23-11} \begin{pmatrix} 2\mu \,\partial_1 u^1_k +p_k  &  \mu \big( \partial_2 u^1_k +\partial_1 u^2_k\big) \\
\mu \big( \partial_1 u^2_k + \partial_2 u^1_k\big) &  2\mu \,\partial_2 u^2_k +p_k\end{pmatrix} \begin{pmatrix} \nu^1 \\
\nu^2 \end{pmatrix} =0  \;\;\;\; \; \mbox{on}\;\;\partial \Omega.\end{eqnarray}
Since $\boldsymbol{\nu}=(\nu^1, \nu^2)$ is a non-zero solution of linear equation (\ref{24.4.23-11}),  we immediately get 
\begin{eqnarray} \label{24.4.23-012}  \mbox{det} \;\begin{pmatrix} 2\mu \,\partial_1 u^1_k +p_k  &  \mu \big( \partial_2 u^1_k +\partial_1 u^2_k\big) \\
\,\mu \big( \partial_1 u^2_k + \partial_2 u^1_k\big) &  2\mu \,\partial_2 u^2_k +p_k\end{pmatrix}\, =0  \;\;\;\; \; \mbox{on}\;\;\partial \Omega,\end{eqnarray} 
Thus,  by (\ref{24.4.23-012}) and $\mbox{div}\, \mathbf{u}_k=0$ on $\Omega$ we have  
\begin{eqnarray}  \label{24.4.2-12}  p_k=\mu \sqrt{\big(\partial_1 u^2_k  + \partial_2 u_k^1 \big)^2 +4\big(\partial_1 u_k^1)^2} \;\;\;\mbox{on}\;\; \partial \Omega.\end{eqnarray}  
Taking the curl of Equation (\ref{24.4.23-10}), we see that the vorticity $\omega_k = \mbox{curl}\,(\mathbf{u}_k)$
satisfies 
\begin{eqnarray} \label{24.4.23.11}  \left\{ \begin{array}{ll}  \mu\,\Delta \omega_k + \varsigma_k \omega_k=0 \;   &\mbox{in}\;\; \Omega,\\
\Big(2\mu \,(\mbox{Def}\; \mathbf{u}_k)-\mu \sqrt{\big(\partial_1 u^2_k  + \partial_2 u_k^1 \big)^2 +4\big(\partial_1 u_k^1)^2}\, \;\mathbf{I}_2\Big) \,\boldsymbol{\nu}=0\;\,  &\mbox{on} \;\;\partial \Omega, \end{array} \right.\end{eqnarray} 
 Here we have used the fact that  $\omega (\nabla p_k)=0$ in $\Omega$, and 
\begin{eqnarray*}  \label{24.4.2-12}  p_k=\mu \sqrt{\big(\partial_1 u_k^2 +\partial_2 u_k^1\big)^2 +4\big(\partial_1 u_k^1)^2} \;\;\;\mbox{on}\;\; \partial \Omega.\end{eqnarray*}  
That is, $\omega_k$ is an eigenfunction of $-\mu \Delta$, but with Cauchy force boundary
conditions on the velocity $\mathbf{u}_k$.

Let  $\psi_k$ be the stream function for $\mathbf{u}_k$. That is,  
\begin{eqnarray} \label{24.6.23-2} \mathbf{u}_k = \nabla^{\perp} \psi_k \;\; \mbox{in}\;\; \Omega.\end{eqnarray} 
Then  $\omega_k =\Delta \psi_k $  in $\Omega$. 
For the known function $\omega_k$ defined in $\Omega$,  let $\psi_k$ is the solution of  the following  Poisson equation:
\begin{eqnarray} \label{24.6.23-3} \left\{ \begin{array} {ll} \Delta \psi_k=\omega_k \;\; &\mbox{in}\;\;\Omega,\\
\psi_k=0\;\;& \mbox{on}\;\; \partial \Omega.\end{array}\right.\end{eqnarray} 
  From (\ref{24.6.23-3}), we have 
\begin{eqnarray*} \left. \begin{array}{ll}  & \partial_1 u^1_k= - \partial_1\partial_2  \psi_k,  \;\;\; \partial_1u_k^2= \partial_1^2 \psi_k,\\          & \partial_2 u_k^1= -\partial_2^2 \psi_k, \;\;\;   \partial_2u_k^2 =\partial_1\partial_2 \psi_k\end{array} \right. \;\; \,\mbox{in}\;\; \Omega.\end{eqnarray*}  
Thus,  $\psi_k$  satisfies 
 \begin{eqnarray}\label{24.4.23-12} \left\{ \begin{array}{ll}  \mu \,\Delta \Delta \psi_k + \varsigma_k \Delta \psi_k = 0 \;\, &\mbox{in}\;\;\Omega,\\
 \psi_k =0 \;\,  &\mbox{on}\;\; \partial \Omega,\\
(B {\psi})\,\boldsymbol{\nu}=0\;\,  &\mbox{on} \;\;\partial \Omega, \end{array} \right.\end{eqnarray} 
where \begin{eqnarray*} B\psi = \mu \begin{pmatrix}  -2 \partial_1\partial_2 \psi_k  &  \partial_1^2 \psi_k -\partial_2^2 \psi_k\\
\partial_1^2 \psi_k -\partial_2^2 \psi_k &  2\partial_1\partial_2 \psi_k \end{pmatrix} - \mu \sqrt{ (\partial_1^2 \psi_k -\partial_2^2\psi_k)^2 +4(\partial_1\partial_2 \psi_k)^2} \;\, \mathbf{I}_2\;\;\;\, \mbox{on}\;\;\partial \Omega. \end{eqnarray*}
This is the corresponding eigenvalue problem for the buckling plate with Cauchy force boundary condition.

 Conversely, given $\psi_k$
satisfying Equation (\ref{24.4.23-12}),  $\omega_k= \Delta\psi_k$  and $\mathbf{u}_k =\nabla^{\perp} \psi_k$ satisfy Equation (\ref{24.4.23.11})
and one can show, at least for sufficiently smooth boundaries, that $\mathbf{u}_k$ satisfies Equation (\ref{24.4.23-10}). Note that $(0, c)$ is not a buckling eigenpair with Cauchy force boundary condition. Thus, except the first Stokes eigenpair, the eigenvalue problem for the Stokes operator and the buckling plate are equivalent for the Cauchy force boundary condition, which implies (\ref{24.4.24-01}).
 \qed

   \vskip 0.3  true cm 

   \noindent{\bf  Remark 6.1.} \ {\it For a bounded domain $\Omega$ with smooth boundary in the Euclidean plane $\mathbb{R}^2$,   the previous  formulas  (\ref{24.5.23-1}) in Corollary 1.10 only hold upto two-term asymptotics. It is easy to see that the third coefficients of the asymptotic expansions   for  $\sum_{k=1}^\infty e^{-t \Xi_k}$,      $\;  \sum_{k=1}^\infty e^{-t \Gamma_k}$ and $ \sum_{k=1}^\infty e^{-t \Lambda_k^{(D)}}$ are different.  Note also that  $\mu \Delta$ is the usual Laplacian defined on $0$-form and  $L_{-\mu, \mu}=\mu \Delta$ is the  (vector form)  Laplacian defined on $1$-form.}

\vskip 1.16 true cm

\section*{Acknowledgments}

 This research was supported by the NNSF of China (12271031)
 and  NNSF of China (11671033/A010802).

\vskip 0.99 true cm

\end{document}